\pdfoutput=1
\documentclass[reprint,amsmath,amssymb,aps,pre,superscriptaddress,nofootinbib,floatfix]{revtex4-2}

\usepackage{graphicx}
\usepackage{dcolumn}
\usepackage{bm}
\usepackage{mathtools}
\usepackage{amsthm}
\usepackage{xcolor}
\usepackage[colorlinks=true,hypertexnames=false,linkcolor=blue!55!black,citecolor=blue!55!black,urlcolor=blue!55!black]{hyperref}
\usepackage{booktabs}
\usepackage{enumitem}
\usepackage[nameinlink,capitalise,noabbrev]{cleveref}

\hypersetup{pdftitle={Predicting phase memory and slow relaxation in synchronized networks},pdfauthor={Francesco Paolo Maiale}}

\newif\ifarxiv
\arxivtrue
\newcommand{\SMsee}{\ifarxiv See the Supplemental Material, available as an ancillary file with this arXiv
  submission,\else See Supplemental Material at [URL will be inserted by publisher]\fi}

\newcommand{\R}{\mathbb R}
\newcommand{\Z}{\mathbb Z}
\newcommand{\one}{\bm 1}
\newcommand{\trans}{\mathsf T}
\newcommand{\diag}{\operatorname{diag}}
\newcommand{\osc}{\operatorname{osc}}
\newcommand{\circu}{\operatorname{circ}}
\newcommand{\norm}[1]{\left\lVert#1\right\rVert}
\newcommand{\abs}[1]{\left\lvert#1\right\rvert}

\theoremstyle{plain}
\newtheorem{theorem}{Theorem}
\newtheorem{proposition}{Proposition}

\newtheorem{appthm}{Theorem}[section]
\newtheorem{appprop}[appthm]{Proposition}
\newtheorem{appcor}[appthm]{Corollary}
\newtheorem{applem}[appthm]{Lemma}
\theoremstyle{definition}
\newtheorem{appdefn}[appthm]{Definition}
\newtheorem{appassum}[appthm]{Assumption}
\theoremstyle{remark}
\newtheorem{apprem}[appthm]{Remark}
\theoremstyle{plain}
\crefname{appthm}{Theorem}{Theorems}\crefname{appprop}{Proposition}{Propositions}
\crefname{appcor}{Corollary}{Corollaries}\crefname{applem}{Lemma}{Lemmas}
\crefname{appdefn}{Definition}{Definitions}\crefname{appassum}{Assumption}{Assumptions}
\crefname{apprem}{Remark}{Remarks}

\newcommand{\proofin}[1]{(see \cref{a:#1})}
\newcommand{\numcheck}[2][]{\emph{Numerical check\if\relax\detokenize{#1}\relax\else~[#1]\fi.} #2}

\begin{document}

\title{Predicting phase memory and slow relaxation in synchronized networks}

\author{Francesco Paolo Maiale}
\email{francescopaolo.maiale@gssi.it}
\affiliation{Gran Sasso Science Institute, Viale Francesco Crispi 7, 67100 L'Aquila, Italy}

\date{\today}

\begin{abstract}
After a weak local kick, synchronized oscillators retain a common phase shift. With positive effective rates, the fraction retained in linear response is the stationary probability of the kicked node in the linearized random walk.
On a semi-infinite chain of identical Sakaguchi oscillators, these weights decay away from the free end with exact length $\xi=-1/\ln\cos 2\alpha$ for $0<|\alpha|<\pi/4$.
On the symmetric locked branch of a long even chain, equal kicks at opposite ends leave the same permanent shift. The sign of the difference between the mean phase displacements of the two halves identifies the kicked end, and its relaxation time grows as $e^{N/(2\xi)}$ at fixed $\alpha$.
On untwisted networks with small edge lags and positive coupling slope at zero, we bound the error in predicting log weights from lifted phases with slope $2\Gamma''(0)/\Gamma'(0)$.
On a fixed finite connected graph with positive rates, the self-consistent weighted potential is unique up to a constant and its log-weight error is $O(\eta^3)$ as the largest absolute edge residual $\eta$ tends to zero. A separate stationarity-residual bound applies to computed candidates under explicit smallness conditions.
Paired kicks detect stationary currents on twisted rings. In three full oscillator models with non-odd coupling, the predicted preference for phase leaders or laggers holds in every configuration eligible for comparison. In fixed-voltage grid models with equal damping and uniform balancing, relative memory equals one minus the marginal loss sensitivity.
\end{abstract}

\maketitle

\section{Introduction}\label{sec:intro}

\begin{figure*}[t]
\includegraphics[width=\textwidth]{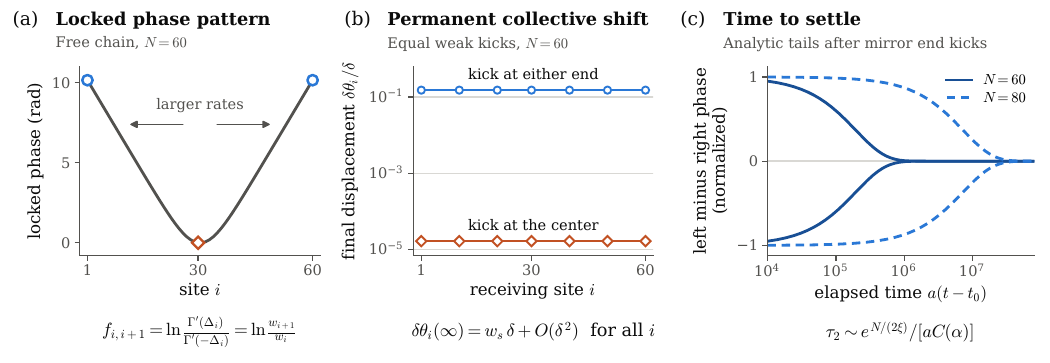}
\caption{\label{fig:concept}From the locked phase pattern to the permanent shift and its
relaxation on the positive symmetric Sakaguchi branch at $\alpha=0.4$.
(a) Lifted phase
relative to the center of a free chain of $N=60$ sites. Arrows indicate the direction of
larger effective rates. We reconstruct the phases from the exact ratios of adjacent weights.
(b) In linear response to a kick of size $\delta$, every receiving site retains the same
displacement $w_s\delta$. A kick at either end retains $15\%$ of the kick, versus
$1.7\times10^{-5}$ of the kick at site $30$. (c) Kicks at opposite ends excite the
antisymmetric slow mode with opposite signs. The curves show its analytic late-time contribution
to the difference between the mean phase displacements of the left and right halves,
normalized by its magnitude at a late reference time $t_0$. The displayed asymptote holds
for long even chains at fixed $0<\alpha<\pi/4$, with $\xi=-1/\ln\cos2\alpha$ and
$C(\alpha)$ from Theorem~\ref{thm:sink}. Increasing the chain from $60$ to $80$ sites
increases the relaxation time $37$-fold, from $a\tau_2\approx1.95\times10^5$ to $7.25\times10^6$.}
\end{figure*}

What can a synchronized network remember about where it was kicked?
After a weak local kick, the oscillators recover their locked phase differences and retain a
common phase shift. On the 60-site chain in Fig.~\ref{fig:concept}, a kick at either end
leaves about $15\%$ of its amplitude in that shift. A kick at a central site leaves only
$1.7\times10^{-5}$ of its amplitude. The permanent shift distinguishes an end from the
center. The transient also distinguishes the two ends. Kicks of the same sign at opposite
ends give opposite signs for the difference between the mean phase displacements of the
left and right halves. Light pulses in circadian systems and current pulses in electrochemical
arrays can probe collective phase responses~\cite{Hannay2015,Kiss2002}.

In linear response the collective shift after a kick $\delta$ at node $s$ is $w_s\delta$, where
$w$ is the normalized left null vector of the locked state's Jacobian~\cite{Kori2009}. With
positive effective rates, that Jacobian generates a random walk and $w$ is its stationary
distribution~\cite{TonjesBlasius2009b}. Odd coupling gives uniform weights on a symmetric graph.
A phase lag can bias the walk toward phase leaders or laggers~\cite{Kori2009}. The same
weights describe influence in directed networks~\cite{Masuda2009}, and localization near
a pacemaker changes collective noise~\cite{Masuda2010}. Related work addresses entrainment and signal
routing~\cite{KoriMikhailov2004,Kirst2016} and the transient spread of disturbances in
grids~\cite{Kettemann2016}, where resistance distances govern certain time-integrated
responses with the mean phase removed~\cite{TylooColettaJacquod2018,TylooPagnierJacquod2019}.

We ask how the locked phases determine the source weights $w_s$ and how long the
approach to the permanent shift retains information about the source. For an untwisted
state with small edge lags and positive coupling slope at zero, we bound the error of
predicting $\ln w_s$ from the lifted phase at $s$. This gives a discrete error estimate
for the exponential density suggested by the continuum theory~\cite{TonjesBlasius2009a}.
A self-consistent potential uses the directional rates and has log-weight error
$O(\eta^3)$ at fixed graph size as the largest absolute edge residual $\eta$ tends to zero.
A separate bound uses the stationarity residual of a computed candidate. On the symmetric
Sakaguchi chain, the locked lags fix both the decay length $\xi$ of the weights near a
free end and the conductances that govern exchange between the two halves. In long even
chains at fixed phase lag $0<|\alpha|<\pi/4$, the exchange time grows as $e^{N/(2\xi)}$.
We calculate its prefactor and give a length-independent lower bound on all faster decay rates. This discrete calculation
builds on the continuum free-end condition, monotone chain profiles and competing-pacemaker
models~\cite{KopellErmentrout1986,RenErmentrout1998,TonjesBlasius2009a}.

Sections~\ref{sec:setup}--\ref{sec:identical} develop the weight estimates, test them on a
graph atlas and derive a degree law for identical oscillators. Sections~\ref{sec:geometry}
and~\ref{sec:topology} trace the response on chains, boxes and rings. Sections~\ref{sec:oscillators}
and~\ref{sec:networks} compare phase predictions with full oscillator dynamics and apply
the weight laws to lossy grids and measured connectome wiring. Section~\ref{sec:outlook}
discusses collective noise and possible kick experiments. Methods and further checks are
in the Supplemental Material~\cite{SM}.

\section{Perturbation memory as a stationary distribution}\label{sec:setup}

We consider $N\ge2$ phase oscillators on a finite connected undirected graph $G=(V,E)$ with symmetric weights
$a_{ij}=a_{ji}\ge0$ satisfying $a_{ii}=0$ and $a_{ij}>0$ on edges $\{i,j\}\in E$. More precisely, our model is
\begin{equation}\label{eq:model}
  \dot\theta_i=\omega_i+\sum_j a_{ij}\,\Gamma(\theta_j-\theta_i),
\end{equation}
where $\Gamma$ is a smooth $2\pi$-periodic coupling function, for instance obtained by phase
reduction of a limit-cycle oscillator. We refer to a solution
$\theta(t)=\theta^*+\Omega t\,\one$ (with $\one=(1,\dots,1)^\trans$) in which every oscillator runs at the common frequency $\Omega$ as a \emph{locked state}.
We write $\Delta_{ij}=\theta^*_j-\theta^*_i$ for the phase lags across the edges.

\paragraph*{Linear response.}
Let $\theta$ be a locked state. A small perturbation $\delta\theta$ evolves as
$\dot{\delta\theta}=J\,\delta\theta$, which means
\[
\dot{\delta\theta}_i=\sum_j W_{ij}(\delta\theta_j-\delta\theta_i),
\]
with
\begin{equation}\label{eq:J}
  J=W-\diag(W\one),\qquad W_{ij}=a_{ij}\,\Gamma'(\Delta_{ij}).
\end{equation}
The matrix $J$ has zero row sums, because a uniform phase shift is neutral. Throughout, we assume that every
\emph{effective weight} is positive, $W_{ij}>0$ on every edge in both directions, which holds
at sufficiently small lags when $\Gamma'(0)>0$. Oscillator $i$ is then pulled toward neighbor $j$ at rate $W_{ij}$, so $J$ is
the generator of a continuous-time random walk that jumps from $i$ to $j$ at rate $W_{ij}$
\proofin{prop:generator}. The eigenvalue $0$ is simple and all other eigenvalues have negative real
part, recovering the positive-rate case of Ermentrout's stability theorem
\cite{Ermentrout1992,RenErmentrout1998} (the same generator argument covers nonnegative rates with
strongly connected positive-rate support). The stochastic matrix $e^{Jt}$ maps each perturbation to
an average of its initial values, so both the supremum norm and the diameter are nonincreasing.
The unique stationary vector satisfies $w>0$, $w^\trans J=0$ and $\sum_iw_i=1$; see
Ref.~\cite{TonjesBlasius2009b}. The projection $w^\trans\delta\theta(t)$ is conserved, and
\[
\delta\theta(t)\to(w^\trans\delta\theta(0))\,\one.
\]

\paragraph*{Influence.}
Kick oscillator $s$ by a small $\delta$. The network relaxes back to the locked orbit, leaving every oscillator permanently shifted
by
\begin{equation}\label{eq:shift}
  \delta\Psi=w_s\,\delta+O(\delta^2)
\end{equation}
\proofin{prop:shift}. We refer to $w_s$ as the \emph{influence} (or \emph{perturbation memory}) of node $s$. This is the collective-phase gradient of Kori
\emph{et al.}~\cite{Kori2009} and the influence of Masuda, Kawamura and Kori~\cite{Masuda2009,Masuda2010}.
As a left zero mode, $w$ is paired with the uniform right zero mode. The permanent shift measures
the projection on this biorthogonal pair, as does the driven response of non-reciprocal mechanical
metamaterials~\cite{Schomerus2020}.

\paragraph*{Reciprocity and the gauge field.}
Because the graph is undirected, the rates are symmetric whenever $\Gamma'$ is even at the occupied lags.
In particular, an odd coupling function up to an additive constant gives a symmetric walk and
uniform $w$. An even, nonconstant component of $\Gamma$ can make the effective rates asymmetric.
Their asymmetry is carried by the antisymmetric edge function
\begin{equation}\label{eq:gauge}
  f_{ij}:=\ln\frac{W_{ij}}{W_{ji}}=F(\Delta_{ij}),\qquad F(x)=\ln\frac{\Gamma'(x)}{\Gamma'(-x)},
\end{equation}
which does not depend on $a_{ij}$. By analogy with the continuum phase-diffusion equation, where the locked phase profile $\theta^0$ acts as a gauge field inducing a drift $2\gamma\nabla\theta^0$~\cite{TonjesBlasius2009a}, we call the discrete edge bias $f$ the \emph{gauge field}. The function $F$ is odd on its domain. When $\Gamma'(0)>0$,
its expansion near zero is $F(x)=2\gamma x+O(x^3)$, with
\begin{equation}\label{eq:gamma}
  \gamma:=\frac{\Gamma''(0)}{\Gamma'(0)}.
\end{equation}
The coefficient $\gamma$ measures the leading asymmetry at small lags. For sufficiently small positive $\Delta_{ij}$, $\gamma>0$ gives $W_{ij}>W_{ji}$ so the walk drifts toward phase leaders, whereas when $\gamma<0$ it drifts toward phase laggers. For Sakaguchi coupling~\cite{Sakaguchi1986},
$\Gamma(x)=\sin(x-\alpha)$, $\gamma=\tan\alpha$; for the diffusive-type family
$\Gamma(x)=\sin x+c(1-\cos x)$, $\gamma=c$. Writing $S_{ij}=(W_{ij}W_{ji})^{1/2}$, the rates factor as
$W_{ij}=S_{ij}e^{f_{ij}/2}$.

Outside the positive-rate regime, effective weights can change sign. We track where this occurs in the numerical atlas of Sec.~\ref{sec:exact} and discuss signed influence and transient growth in Sec.~\ref{sec:outlook}.

\section{Stationary current and controlled approximations}\label{sec:exact}

The locked state supplies the directional rates. Balancing the stationary currents determines
the weights governing the permanent kick response. For any candidate log-weight potential
$\Phi\in\R^V$, define on each edge the residual and the conductance
\begin{equation}\label{eq:hc}
  h_{ij}(\Phi)=f_{ij}-(\Phi_j-\Phi_i),\quad c_{ij}(\Phi)=S_{ij}\,e^{(\Phi_i+\Phi_j)/2}.
\end{equation}

\begin{proposition}\label{thm:hodge}
The potential $\Phi=\ln w$ is the unique solution, up to an additive constant, of
\begin{equation}\label{eq:H}
  \sum_{i\sim j}c_{ij}(\Phi)\,\sinh\!\big(\tfrac12h_{ij}(\Phi)\big)=0\qquad\text{for all }j.
\end{equation}
The stationary probability current on each edge is
\begin{equation}\label{eq:current}
\begin{aligned}
  \jmath_{ij} & :=w_iW_{ij}-w_jW_{ji} \\ & =2\,c_{ij}(\ln w)\,\sinh\!\big(\tfrac12h_{ij}(\ln w)\big).
  \end{aligned}
\end{equation}
\end{proposition}

Equations~\eqref{eq:H} and \eqref{eq:current} specialize the force--current representation of
Kaiser, Jack and Zimmer~\cite[Eqs.~(12)--(13)]{KaiserJackZimmer2018} to oscillator influence.
For $v=e^\Phi$, $v_iW_{ij}-v_jW_{ji}=2c_{ij}\sinh(h_{ij}/2)$, so $v^\trans J=0$ is equivalent
to Eq.~\eqref{eq:H} \proofin{thm:hodge}. The edge affinity $h$ has the same cycle sums as
$f$; refer to~\cite{Schnakenberg1976}. The decomposition separates the influence potential from the stationary
current, whose divergence vanishes.

When $f$ is exact, $f_{ij}=\Phi_j-\Phi_i$, the residual vanishes, the current vanishes, and the walk is
in detailed balance (Kolmogorov's criterion~\cite{Kolmogorov1936,Kelly1979}). Then
\begin{equation}\label{eq:exactgauge}
  w_i=\frac{e^{\Phi_i}}{\sum_ke^{\Phi_k}}.
\end{equation}
This is always the case on a tree, where $w_j/w_i=\exp\sum_{e\in\mathrm{path}(i,j)}f_e$.

\paragraph*{First order in the edge lags.}
Use principal edge lags ($\abs{\Delta_e}\le d<\pi$) and assume $\Gamma'(x)>0$ on
$[-d,d]$. On a cycle $C$ with winding number $q_C=(2\pi)^{-1}\sum_{e\in C}\Delta_e$,
oddness of $F$ gives for the cycle sum $\circu_C f:=\sum_{e\in C}f_e$ the following identity:
\begin{equation}\label{eq:circ}
  \circu_Cf=4\pi\gamma\,q_C+\rho_C,\qquad\abs{\rho_C}\le\frac{M_3}6\sum_{e\in C}\abs{\Delta_e}^3,
\end{equation}
with $M_3=\sup_{\abs x\le d}\abs{F'''(x)}$ \proofin{prop:circ}. In an untwisted state ($q_C=0$ on
every cycle), the phases admit a real lift $\tilde\theta$ with $\tilde\theta_j-\tilde\theta_i=\Delta_{ij}$
on every edge, and the gauge field is a gradient up to third order. Then
\begin{equation}\label{eq:firstorder}
  \ln w_i=2\gamma\,\tilde\theta_i+\text{const}+O(d^3),
\end{equation}
with the explicit bound
\begin{equation}\label{eq:firstbound}
  \max_i\Big|\ln w_i-2\gamma\tilde\theta_i+\ln\sum_ke^{2\gamma\tilde\theta_k}\Big|\le
  \frac{M_3}6\max_T\sum_{e\in T}\abs{\Delta_e}^3,
\end{equation}
the maximum taken over spanning trees $T$ \proofin{thm:firstorder}. The proof uses the Markov chain
tree theorem~\cite{LeightonRivest1986,AnantharamTsoucas1989}, the spanning-tree (cofactor) form of $w$
also used by Kori \emph{et al.}~\cite{Kori2009}.

Note that Eq.~\eqref{eq:firstorder} turns the sign of $\gamma$ into a quantitative ranking. Indeed, to leading
order, the log-weight contrast is $2|\gamma|\osc\tilde\theta$ (with $\osc x:=\max_i x_i-\min_i x_i$),
with an error at most the spanning-tree bound in Eq.~\eqref{eq:firstbound} \proofin{thm:firstorder}. A phase ranking therefore needs differences large enough to resolve
against the remainder. For two Sakaguchi oscillators the sign rule follows from the exact weights of
Kori \emph{et al.}, who also observed it in globally coupled populations~\cite{Kori2009}. In continuum
phase diffusion the Cole--Hopf substitution linearizes the nonlinear phase-diffusion
equation~\cite[Sec.~6.3]{Kuramoto1984}, including the continuum limit of Sakaguchi oscillator
lattices~\cite[Sec.~3]{SakaguchiShinomotoKuramoto1988} (see also Ref.~\cite{TonjesBlasius2009a}); with a global phase lift and periodic or reflecting
boundaries, the adjoint operator then has stationary density proportional to
$e^{2\gamma\tilde\theta}$; see~\cite{Masuda2010,TonjesBlasius2009a}. For the nonperiodic exponential
surrogate $\Gamma(x)=(e^{\gamma x}-1)/\gamma$ used on the lift in the lattice Cole--Hopf
analysis~\cite{BlasiusTonjes2005}, $F(x)=2\gamma x$ exactly. In an untwisted state the walk is then
reversible and $w\propto e^{2\gamma\tilde\theta}$. This also holds for a periodic coupling that agrees
with the surrogate on the interval of lags used, including both signs. Equation~\eqref{eq:firstbound}
gives a discrete error estimate for the original periodic coupling.

\paragraph*{A self-consistent third-order potential.}
The natural approximation to Eq.~\eqref{eq:H} linearizes $\sinh$ but keeps the conductances
self-consistent. In particular, let $\hat\Phi$ solve
\begin{equation}\label{eq:Hhat}
  \sum_{i\sim j}c_{ij}(\hat\Phi)\,h_{ij}(\hat\Phi)=0\qquad\text{for all }j,
\end{equation}
the weighted least-squares potential of $f$ with conductances $c(\hat\Phi)$. For fixed conductances
this is the standard graph Hodge projection~\cite{JiangLimYaoYe2011}. On every finite connected
graph with positive rates, Eq.~\eqref{eq:Hhat} has a unique solution up to an additive constant
\proofin{prop:globalpotential}.
We use Newton's method to compute the self-consistent potential. An exact solution of
Eq.~\eqref{eq:Hhat} cancels the linear term in the current balance. For a numerical candidate,
the stationarity residual enters the error bound directly.

\begin{theorem}[A posteriori error bound]\label{thm:weighted}
Let $\hat\Phi$ be any potential, $\hat c=c(\hat\Phi)$, $\hat h=h(\hat\Phi)$ and
$\eta=\max_e\abs{\hat h_e}$. Let $g^\trans=(e^{\hat\Phi})^\trans J$ be the stationarity residual of
$e^{\hat\Phi}$, $\hat L$ the Laplacian with conductances $\hat c$, $\hat R$ its largest effective
resistance, $\norm g_*=(g^\trans\hat L^+g)^{1/2}$, $\vartheta=\frac12\hat R\norm g_1$, and
$Q=\hat R\sum_e\hat c_e\ge N-1$.
\begin{enumerate}[label=\textup{(\alph*)}]
\item If $\vartheta<1$ and $\mathcal B:=\sqrt{\hat R}\norm g_*/(1-\vartheta)<1$, then
  \begin{equation}\label{eq:certified}
    \osc(\ln w-\hat\Phi)\le-\ln(1-\mathcal B).
  \end{equation}
\item Under the conditions of \textup{(a)}, the error of the correction $\hat\Phi+\hat L^+g$ is at most
  $\mathcal B\bigl(\sqrt Q\,\eta e^{\eta/2}+\mathcal B/(1-\mathcal B)\bigr)$ in oscillation.
\item If $\hat\Phi$ solves Eq.~\eqref{eq:Hhat} and $\cosh(\eta/2)\eta^3Q\le6$, then
  \begin{equation}\label{eq:third}
    \osc(\ln w-\hat\Phi)\le\frac{\sqrt Q}{12}\cosh(\eta/2)\,\eta^3 ,
  \end{equation}
  and the one-step correction $\hat\Phi+\hat L^+g$ is accurate to $O(Q\eta^4)$.
\end{enumerate}
\end{theorem}

An energy estimate eliminates the skew part of the drift
\proofin{thm:weighted}. In part (a), twice the shortest-path eccentricity of any node, with
edge lengths $1/\hat c_e$, can replace $\hat R$, because an effective resistance is at most a path length. We evaluate that upper bound, the residual $g$
and the grounded solve for $\norm g_*$ in double precision; a numerical certificate would enclose
those quantities and verify the smallness conditions with controlled rounding. Part (c) uses
self-consistency to cancel the linear residual, whereas an unweighted Hodge projection leaves a linear
residual when the conductances vary (on a Barab\'asi--Albert graph with $N=1000$ and $c=1$, its error
has local slope $\approx1.36$ in $\eta$, against $3.0$ for $\hat\Phi$ \proofin{thm:weighted}).

In localized states, weak conductances can make $\hat R$ large and amplify a nonzero residual, although
contrast alone does not decide whether the smallness conditions hold (for instance, the exact
path-integrated potential on a positive-rate tree has $g=0$ and satisfies both conditions at arbitrary
contrast). A separate spanning-tree comparison, applying O'Cinneide's entrywise perturbation
principle~\cite[Thm.~1]{OCinneide1993} to a comparison chain with stationary weights proportional to $e^{\hat\Phi}$, gives cubic order at
fixed $N$ for the exact solution of Eq.~\eqref{eq:Hhat} whenever $e^{\eta/2}\cosh(\eta/2)\eta^3/48<1$,
with a prefactor proportional to $N$ and independent of $\hat R$ \proofin{prop:treeweighted}.

In the 37 states used to develop the bound, the a posteriori bound was within $1.6$--$22$ times
the measured error where it applied (31 states); details are in the Supplemental Material~\cite{SM}.
The correction $\hat\Phi+\hat L^+g$ is a first-order stationary-distribution perturbation
\cite{Schweitzer1968}. Classical normwise bounds use condition numbers from the group inverse,
first-passage times or ergodicity coefficients~\cite{ChoMeyer2001}, which can grow exponentially
in metastable chains. Entrywise relative bounds avoid that growth
\cite{OCinneide1993,ThiedeVanKotenWeare2015}. Theorem~\ref{thm:weighted} controls the
oscillation of $\ln w$ using a candidate's residual and symmetric-network resistances.

\begin{figure*}[t]
\includegraphics[width=\textwidth]{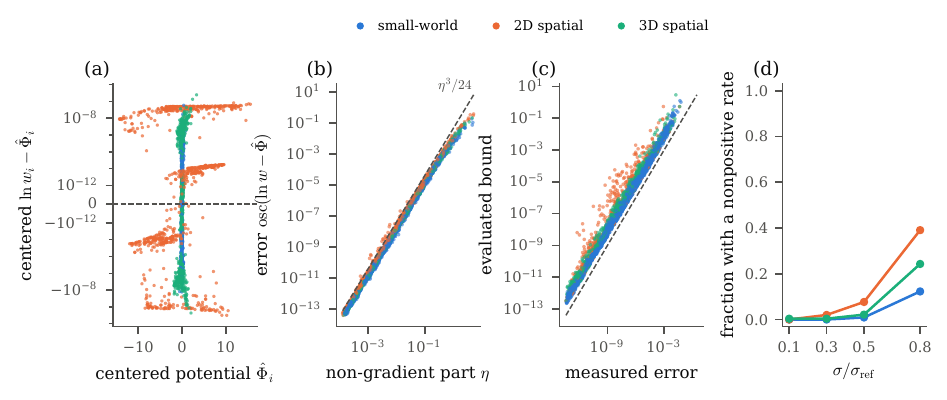}
\caption{\label{fig:atlas}Disorder atlas across eight graph families. Colors group small-world (Barab\'asi--Albert, Erd\H{o}s--R\'enyi, Watts--Strogatz and block-model),
two-dimensional spatial and three-dimensional spatial families.
(a) Centered residual $\ln w-\hat\Phi$ against centered $\hat\Phi$ for one eligible state per family
($c=1$, $\sigma=0.5\sigma_{\rm ref}$, largest available $N$, $400$ sampled nodes each).
The vertical axis uses a symmetric logarithmic scale.
(b) Error against the largest edge residual $\eta$ (axis label ``non-gradient part'') for eligible states. Dashed, the heuristic $\eta^3/24$.
(c) Finite double-precision bound against measured error. Dashed, equality.
(d) Fraction of saved locks with a nonpositive rate against $\sigma/\sigma_{\rm ref}$.}
\end{figure*}

\paragraph*{Testing the approximation on an atlas.}
Across the atlas, the measured error in the self-consistent log weights typically scaled
cubically at small edge residuals. The path-length test for the error bound usually gave
a finite bound at low weight contrast and rarely did so in strongly localized states
(Fig.~\ref{fig:atlas}). We tested eight graph families with up to $N=10^5$ nodes:
Barab\'asi--Albert, Erd\H{o}s--R\'enyi, Watts--Strogatz and two-block stochastic block models; random
geometric graphs in two and three dimensions; and square and cubic lattices. The six couplings are
Sakaguchi $\alpha\in\{0.1,0.2,0.3\}$ and diffusive type $c\in\{0.5,1,2\}$. In each cell (one graph realization with one coupling), frequencies are sampled at
$\sigma/\sigma_{\rm ref}\in\{0.1,0.3,0.5,0.8\}$ using an operational disorder reference from
continuation. The largest three-dimensional samples use one realization; the full design is in the
Supplemental Material~\cite{SM}.
The design and five predictions were fixed in advance, after pilot calculations; the reported numbers
come from a later analysis with repaired solvers, a recomputed disorder reference and the accuracy
gate defined below. Two predictions failed: the cap $\osc(\ln w-\hat\Phi)/\eta^3\le1/24$ in every state and the
$95\%$ pass rate of the degree law (Sec.~\ref{sec:identical})~\cite{SM}.
In $304$ of the $1552$ cells, continuation found no positive disorder
amplitude $\sigma_{\rm ref}$, so those cells were not sampled.
All of these have Sakaguchi coupling: $96$ of the $99$ Sakaguchi cells on two-dimensional and $84$ of $91$ on
three-dimensional random geometric graphs, $60$ and $58$ of $99$ on Erd\H{o}s--R\'enyi and block-model graphs,
and $6$ of $99$ on Barab\'asi--Albert graphs, consistent with the continuation stops of Sec.~\ref{sec:identical}.
On random geometric graphs the disorder atlas therefore rests almost entirely on diffusive-type coupling.

The reference $\sigma_{\rm ref}$ is the largest frequency spread at which our continuation from
the in-phase state still finds a locked state, located by doubling and bisection to $2\%$. It
describes the numerical procedure and is not a locking threshold: another continuation route can
find a lock at a larger spread or on another branch. A lock with positive rates is stable
(Sec.~\ref{sec:setup}). Where a rate is negative we test stability numerically, with the full
transverse spectrum on small graphs and sparse eigenvalue searches on large ones; such a search
can miss an unstable mode. The error $e=\osc(\ln w-\hat\Phi)$ is often close to round-off. We
therefore compute it from a correction to the density $e^{\hat\Phi}$, not as a difference of two
logarithms, and we use a state only when $e$ exceeds ten times its estimated numerical
floor~\cite{SM}. We call this test the accuracy gate and the states that pass it eligible.
Of $4988$ locks that pass the stability test, $4680$
have positive effective rates and $4195$ are eligible (Fig.~\ref{fig:atlas}). The double-precision bound is finite in $3882$,
with $0$ observed violations and median bound-to-error ratio
$8.1$ (range $3.1$ to $1.1\times10^{4}$). In $313$ further
eligible states the test using the path-length upper bound on resistance did not establish
the smallness conditions. This test gives a finite bound in
$99.6\%$ of eligible states with $\osc\ln w<3$, in $36\%$ of those
with $10\le\osc\ln w<20$, and in $0.6\%$ of the $155$ states with
$\osc\ln w\ge20$ (for example, the two-dimensional random geometric state in Fig.~\ref{fig:atlas}(a) has
$\osc\ln w\approx32$ and $\vartheta\approx1.69\times10^{8}$).
\ifnum0>0\relax The bound evaluation remains unresolved in a further
$0$ eligible states.\fi

Among $888$ cells with at least three eligible points at $\eta\le0.1$, the median
log--log slope is $2.98$, and $96\%$ lie in $[2.7,3.3]$ (range $1.83$ to $3.63$).
The median cubic ratio $e/\eta^3$ is $0.015$,
with $90\%$ below $0.022$. While, to leading order, a single-edge cubic flow moves the potential
by at most $\abs{\hat h_e}^3/24$ ($1/24\approx0.042$), multiple edges can reinforce one another:
$21$ atlas states meet or exceed $1/24$, and the largest eligible measured ratio is
about $0.182$.

\section{Degree heterogeneity and phase memory}\label{sec:identical}

The first-order law~\eqref{eq:firstorder} makes the locked profile the link between graph structure
and memory. Even identical oscillators can acquire phase differences when the zero-lag input
varies from node to node.

\paragraph*{The degree law.}
Take identical oscillators with Sakaguchi coupling $\Gamma(x)=\sin(x-\alpha)$ and uniform edge weights
$a_{ij}=KA_{ij}$, where $A$ is the unweighted adjacency matrix, $k_i=\sum_j A_{ij}$ are the degrees,
$L=\diag(k)-A$ is the graph Laplacian, and $L^+$ is its pseudoinverse. The zero-lag value $\Gamma(0)=-\sin\alpha\neq0$ shifts the effective
frequency of node $i$ by $-Kk_i\sin\alpha$. The locked phase depends on the response of the whole
graph to those offsets. Continuing the in-phase state in
$\beta=\tan\alpha$ gives $\theta=-\beta L^+k+O(\beta^3)$, the identical-oscillator case of the
linearized phases of Ref.~\cite{Skardal2015}, and with Eq.~\eqref{eq:H}:

\begin{theorem}[Degree law]\label{thm:degree}
For $\alpha$ near zero with $\abs{\tan\alpha}$ below a graph-dependent constant, the
in-phase-continued state exists, all its
effective weights are positive, and
\begin{equation}\label{eq:degree}
  \ln w_i=-2\tan^2\!\alpha\,(L^+k)_i+\text{const}+O(\tan^4\alpha).
\end{equation}
$\ln w$ is an even function of $\alpha$.
\end{theorem}

The proof is an implicit-function argument in $\beta$, matched order by order to Eq.~\eqref{eq:H}
\proofin{thm:degree}. The evenness is exact, not perturbative: physically, flipping the sign of $\alpha$
reverses both the locked phase ordering ($\theta^*\to-\theta^*$) and the leader/lagger drift bias
($\gamma\to-\gamma$), so the two sign changes cancel. Formally, the reflected coupling
$\tilde\Gamma(x)=-\Gamma(-x)$ yields the reflected locked state $-\theta^*$ with the \emph{same}
Jacobian, so $w$ and the whole relaxation spectrum coincide; for Sakaguchi coupling
$\tilde\Gamma_\alpha=\Gamma_{-\alpha}$ \proofin{prop:reflect}. Numerically the two branches agree to the
last bit.

Equation~\eqref{eq:degree} ranks influence by the degree potential $-L^+k$ for either sign of the
phase lag. In terms of the effective resistance distances $R_{\rm eff}(i,j)$ of the unweighted graph,
$(L^+k)_i=-\tfrac12\sum_j(k_j-\bar k)\,R_{\rm eff}(i,j)+\text{const}$, so at this order
$\ln w_i=\tan^2\!\alpha\sum_j(k_j-\bar k)\,R_{\rm eff}(i,j)+\text{const}$.
The sum weights each resistance distance by the signed degree deviation $k_j-\bar k$, so
nodes of the same degree can have different influence. Resistance
distances also rank the transient vulnerability of nodes in lossless networks~\cite{TylooPagnierJacquod2019},
a different quantity. On Barab\'asi--Albert and Watts--Strogatz graphs
with $N=1000$ and $\alpha=0.1$--$0.3$, the correlation between $\ln w$ and the prediction is
$0.9985$--$0.99999$. The error of Eq.~\eqref{eq:degree} scales as $\tan^{4.02}\alpha$ and
$\tan^{4.03}\alpha$ on the two families, and the error of the phases as $\tan^{3.01}\alpha$ and
$\tan^{3.02}\alpha$, as predicted.

\paragraph*{Heterogeneous frequencies.}
For couplings with $\Gamma(0)=0$ and $\Gamma'(0)>0$, such as the diffusive-type family, the
in-phase state of identical oscillators is stable and has uniform influence. Along the branch
from this state with uniform coupling $a_{ij}=KA_{ij}$ and frequency spread $\sigma$, to first order in $\sigma$ \proofin{rem:freqlaw},
\begin{equation}\label{eq:freqlaw}
  \ln w=\frac{2\gamma}{K\Gamma'(0)}\,L^+\omega+\text{const}+O(\sigma^2).
\end{equation}
The sign of $\gamma$ determines whether influence increases or decreases with the phase potential
$L^+\omega$ at this order. Spatial graphs can sustain a large spread of that potential when frequency
heterogeneity excites their slowly varying Laplacian modes.

\paragraph*{Validity domain.}
The degree law is a small-$\alpha$ expansion on the in-phase-continued branch. Summing the locking
equations over a set $S$ gives an exact identity and a necessary condition for a locked state with
principal lags bounded by $d\in[0,\pi]$ \proofin{prop:cut}:
\begin{equation}\label{eq:cutcrit}
  \abs{\sin\alpha}\Big(\cos d\,\abs S\bar k-2e_{\mathrm{int}}(S)\Big)\le\abs{\partial S},
\end{equation}
where $\abs{\partial S}$ counts the edges leaving $S$, $e_{\mathrm{int}}(S)$ counts the internal
edges and $\bar k$ is the mean degree. This is the
cut-set form of the classical single-node condition~\cite[Lemma~4.2 and Sec.~4.3]{DorflerBullo2014},
applied to the effective frequencies $-Kk_i\sin\alpha$. For $\alpha>0$, low-degree structures have less
negative frequency offsets $-Kk_i\sin\alpha$ and thus run faster than the bulk; a subgraph $S$ attached
to the rest of the network by few edges cannot be held back once the coupling capacity $\abs{\partial S}$
of its boundary is exceeded. Dropping the finite-lag
correction gives a cut estimate $\alpha_{\rm cut}$ for positive phase lag, minimized over single nodes,
dangling trees outside the 2-core, and maximal degree-2 paths in the 2-core~\cite{SM}. The resulting conditions
are $\sin\alpha\lesssim k_i/\bar k$ for a single node and
$\sin\alpha\lesssim1/(n(\bar k-2)+2)$ for a dangling tree of $n$ nodes when its denominator is
positive.

\begin{figure*}[t]
\includegraphics[width=\textwidth]{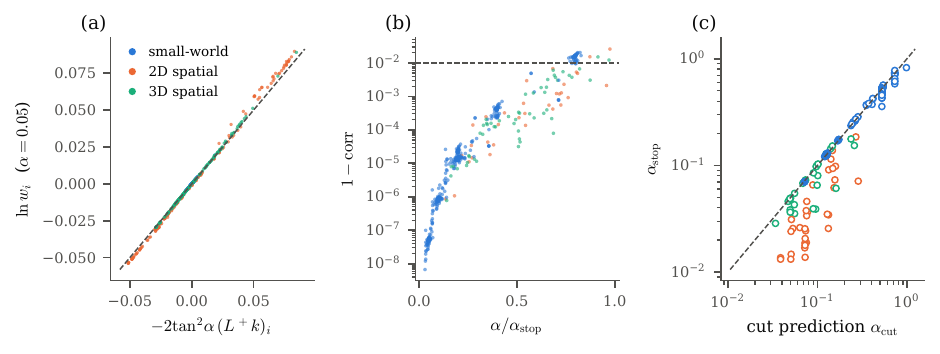}
\caption{\label{fig:degree}Identical oscillators across graph families. Colors as in
Fig.~\ref{fig:atlas}. (a) $\ln w$ against the degree law~\eqref{eq:degree} at $\alpha=0.05$, both
centered, for two realizations per group with $N\le10^4$; dashed, equality. (b) $1-{\rm corr}$ between
$\ln w$ and $-2\tan^2\!\alpha\,L^+k$, against $\alpha/\alpha_{\rm stop}$, in non-box locks that pass the stability test
with $\alpha\le0.1$ and positive rates; dashed, the threshold ${\rm corr}=0.99$.
(c) Non-box continuation stops against the leading cut estimate $\alpha_{\rm cut}$;
filled symbols mark small-world realizations whose predicted limiting cut lies outside the 2-core.
Dashed line, equality.}
\end{figure*}

The leading cut estimate can be inaccurate at large lags or when a different part of the
graph limits continuation. In the identical-oscillator atlas, we denote by $\alpha_{\rm stop}$ the largest
$\alpha$ reached by continuation on each graph and use it to normalize the fixed-lag samples.
Free boxes have a positive stable in-phase-continued state throughout $0\le\alpha<\pi/4$
(Sec.~\ref{sec:geometry}), so box stops in that interval are numerical. We compare the cut estimate,
the numerical stops, and the degree-law accuracy on the non-box graphs (Fig.~\ref{fig:degree}):
\begin{itemize}
\item \emph{Dangling-tree cuts.} On small-world graphs whose critical set predicted by Eq.~\eqref{eq:cutcrit} is a dangling tree,
  $\alpha_{\rm stop}$ agrees with the prediction to within $5\%$ in all $55$ realizations, the largest
  deviation being $2.0\%$. Excluding the 15 realizations seen before the protocol
  was fixed leaves 40, all within $2.0\%$ of the prediction.
\item \emph{2-core cuts.} On the other small-world graphs, whose critical set lies in the 2-core, the median of
  $\alpha_{\rm stop}/\alpha_{\rm cut}$ falls from $1.00$--$1.03$ at $N=10^2$ to $0.88$ (Barab\'asi--Albert) and
  $0.82$ (Watts--Strogatz) at $N=10^5$.
\item \emph{Random geometric graphs.} The median of $\alpha_{\rm stop}$ falls from $0.092$
  at $N=10^2$ to about $0.013$ at $10^5$ in two dimensions (three realizations; two of these stops are conservative
  values from a coarser search), and from $0.10$ to
  $0.044$ at $10^4$ in three. The median of $\alpha_{\rm stop}/\alpha_{\rm cut}$
  ranges from $0.29$ to $0.66$ and from $0.89$ to $0.99$, respectively.
  On four two-dimensional realizations ($N=10^3$ and $10^4$),
  pseudo-arclength continuation confirmed local folds $0.26$--$1.10\%$ above the stops,
  while an alternate continuation route on one $N=10^4$ realization reached a coexisting branch that continued $15\%$ above
  the local fold~\cite{SM}. These checks concern the traced solution branches;
  the first loss of transverse stability along each branch is a separate question.
\item \emph{Degree-law accuracy.} In the non-box population with $\alpha\le0.1$ and positive rates, the correlation between $\ln w$
  and $-2\tan^2\!\alpha\,L^+k$ is at least $0.99$ in $417$ of $417$ states below
  $0.75\,\alpha_{\rm stop}$ and in $11$ of $60$ closer to the stop.
  The median two-point log--log slope of the error against $\tan\alpha$, from
  $\alpha=0.025$ to $0.05$, is $4.06$ over $159$ non-box realizations.
\end{itemize}

\section{Boundary localization and slow relaxation}\label{sec:geometry}

The degree law applies near the in-phase state. On trees, the exact path formula of
Sec.~\ref{sec:exact} also applies at large influence contrast as long as the rates are positive.
We use free chains to relate this contrast to the time scale of a kick response. Products of
the same chains give the spatial profile on a box.

\begin{figure*}[t]
\includegraphics[width=\textwidth]{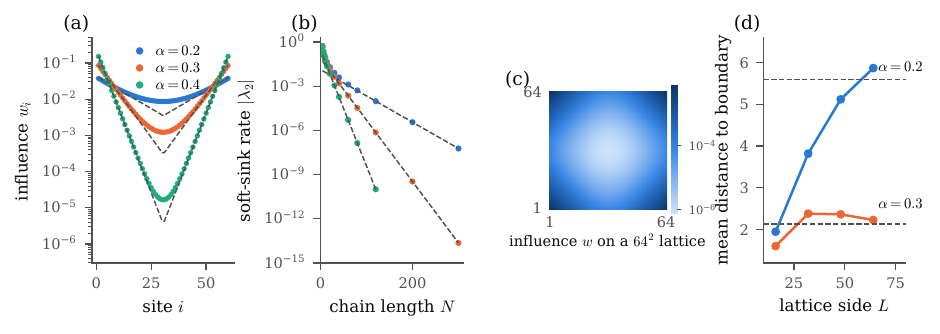}
\caption{\label{fig:geometry}Geometry. (a) Influence on an identical Sakaguchi chain of $N=60$ oscillators
(markers) against the long-chain end profile $\frac12(1-r)\,r^{k-1}$, with $k=1$ at the nearer end and
$r=\cos2\alpha$ (dashed). (b) Soft-sink rate $\abs{\lambda_2}$ against $N$, computed in 160-digit
arithmetic (markers), and its large-$N$ asymptote $C(\alpha)\,r^{N/2}$ with $a=1$ (dashed). (c) Influence on a
$64\times64$ lattice at $\alpha=0.3$, on a logarithmic scale: the memory sits in the corners. The four $6\times6$
corner boxes hold $45\%$ of it, and on this lattice $w=w_1\otimes w_1$ holds to $4.6\times10^{-13}$. (d) Mean
distance of the memory to the boundary against the lattice side, with the asymptotes $r^2/(1-r^2)=2.14$
($\alpha=0.3$) and $5.59$ ($\alpha=0.2$) (dashed). Both curves overshoot before converging, because at
intermediate sides the center of each chain, where the lags vanish, still holds memory. In follow-up
calculations outside the plotted range, at $\alpha=0.2$ the mean distance is about $6.2$ at $L=96$
and $5.60$ at $L=256$.}
\end{figure*}

\paragraph*{Exact skin length.}
On a tree the gauge field is exact, so every positive locked state of a chain obeys
$w_{i+1}/w_i=\Gamma'(\Delta_i)/\Gamma'(-\Delta_i)$ exactly (Sec.~\ref{sec:exact}). Take identical
oscillators with uniform coupling $a$, and let $x_0$ be a zero of $\Gamma$ with $\Gamma'(\pm x_0)>0$.
The lag $-x_0$ is the free-end boundary condition of Kopell and Ermentrout's continuum
limit~\cite{KopellErmentrout1986}. Bounds on the lags of finite chains for monotone profiles
are given by Ren and Ermentrout~\cite{RenErmentrout1998}.
On a semi-infinite chain the lag $-x_0$ holds exactly on every edge: the profile
$\theta_i=-(i-1)x_0$ is an exact locked state, and its generator is a uniform Hatano--Nelson
chain~\cite{HatanoNelson1996,Bergholtz2021} with rates $a\Gamma'(-x_0)$ toward increasing $i$ and
$a\Gamma'(x_0)$ toward the free end. Hence
\begin{equation}\label{eq:skin}
  w_i=(1-r)\,r^{\,i-1},\qquad r=\frac{\Gamma'(-x_0)}{\Gamma'(x_0)},\qquad \xi=-\frac1{\ln r},
\end{equation}
whenever $r<1$ \proofin{prop:skin}. For Sakaguchi coupling $x_0=\alpha$ and
\begin{equation}\label{eq:xi}
  \xi=-\frac1{\ln\cos2\alpha},\qquad0<\abs\alpha<\pi/4 .
\end{equation}
The memory of a perturbation is stored at the free end, the pacemaker of the chain, and decays into
the bulk by a factor $\cos2\alpha$ per site. A geometric profile is elementary on a chain whose links
carry a fixed bias~\cite{Masuda2009}, and in active mechanical metamaterials such a bias is engineered
by feedback, with skin modes and unidirectional amplification observed~\cite{Brandenbourger2019,Ghatak2020};
here the locked state itself sets the bias. For diffusive coupling $x_0=0$ and $r=1$: an identical
diffusive chain on its in-phase branch has uniform influence $w_i=1/N$ at finite $N$. On the semi-infinite chain the constant profile at $r=1$ has infinite mass.
In direct continuations of finite chains the fitted decay rate agrees with $-\ln\cos2\alpha$ to within
$0.2\%$ once the wave is established [Fig.~\ref{fig:geometry}(a)].

\paragraph*{The soft sink.}
A finite chain has two free ends, so two pacemakers. Their waves meet at a sink in the middle, which we call soft because the conductances there are small but not zero. For
Sakaguchi coupling with $0<\alpha<\pi/4$ and even $N=2n$, there is exactly one symmetric positive
locked state up to a common phase, and it is stable. The lags approach $-\alpha$ near a free end
in the long-chain limit and tend to $0$ towards the center along an explicit
``center orbit'' \proofin{lem:twosource}; existence, monotonicity and stability of such profiles were
proved for a broad class of couplings by Ren and Ermentrout~\cite{RenErmentrout1998}. The stationary
weight sits at both ends, and the small conductances near the sink make the end-to-end
resistance large. The slowest mode is antisymmetric and governs exchange between the two halves
of the walk. Exponentially close
rates were found for competing pacemakers in the continuum~\cite{TonjesBlasius2009a}; on the
discrete chain the end-to-end resistance gives two-sided bounds and the asymptotic prefactor.

\begin{theorem}[Soft sink]\label{thm:sink}
Let $c_k=w_kW_{k,k+1}$ be the edge conductances of the two-source state,
$R_N=\sum_kc_k^{-1}$ the end-to-end resistance, $h_N(i)=R_N^{-1}\sum_{j<i}c_j^{-1}$ the equilibrium
potential between the ends and $V_N=\sum_iw_i\,(1-2h_N(i))^2$. For every even $N$ the second
eigenvalue of $J$ is real, its eigenvector is antisymmetric, and
\begin{equation}\label{eq:sinkbounds}
  \frac4{R_N}\le\abs{\lambda_2}\le\frac4{R_N\,V_N},\qquad 0\le1-V_N\le\frac{2n(n-1)\,r^{\,n-1}}P,
\end{equation}
with $r=\cos2\alpha$ and a constant $0<P\le1$ given by a product along the center orbit.
For even $N\ge4$ the third relaxation rate satisfies
\begin{equation}\label{eq:thirdrate}
  \abs{\lambda_3}\ge aP\cos\alpha\,(1-\sqrt r)^2>0.
\end{equation}
At fixed $\alpha$ as $N\to\infty$,
\begin{equation}\label{eq:sink}
  \lambda_2=-aC(\alpha)r^{N/2}(1+o(1))=-aC(\alpha)e^{-N/(2\xi)}(1+o(1)),
\end{equation}
where $C(\alpha)=4/A(\alpha)$ and $A(\alpha)=\lim_{n\to\infty}a\,r^{n}R_N$ is given explicitly by a
convergent series along the center orbit, and
$C(\alpha)=8\alpha^4\,(1+O(\alpha^2))$ as $\alpha\to0$.
\end{theorem}

The bounds~\eqref{eq:sinkbounds} apply the hitting-time and capacity estimates of the
potential-theoretic approach to metastability~\cite{BEGK2002} to the odd sector, using the explicit
birth--death forms of the equilibrium potential and the capacity~\cite[Secs.~7.1.4 and 7.3]{BovierdenHollander2015};
they hold for every even $N$ \proofin{thm:sink}. The small-$\alpha$ limit follows from a tanh kink: the
center orbit satisfies $\abs{\Delta^*_m+\alpha\tanh(\alpha^2m)}\le K\alpha^3e^{-\alpha^2m}$ for $\alpha$ below a constant $\alpha_0$, which gives
$C=8\alpha^4(1+O(\alpha^2))$ for even Sakaguchi chains \proofin{prop:continuum}.
More generally, let $\Gamma(x_0)=0$ with $x_0>0$, $\Gamma'>0$ on $[-x_0,x_0]$,
$\Gamma_e(x)<\Gamma_e(x_0)$ for $0\le x<x_0$, and $r<1$, where $\Gamma_e$ is the even part.
For either parity the symmetric lock exists. At fixed coupling $\abs{\lambda_2}$ scales as
$r^{\lfloor N/2\rfloor}$ with a parity-dependent prefactor, and $\abs{\lambda_3}\ge c/N$
for $N\ge4$ with $c>0$ independent of $N$
\proofin{thm:general}. For families approaching $x_0=0$ with a quadratic leading variation
of $\Gamma_e$ across the occupied lags and nearly constant odd-part slope (see \cref{a:rem:weakskin}),
the formal weak-skin approximation is
$\abs{\lambda_2}\approx2a\Gamma'(0)(1-r)^2r^{N/2}$ (for odd $N$ the exact exponent is $\lfloor N/2\rfloor$,
a difference of $\sqrt r\to1$ in this limit). Its prefactor can be too small by an order
of magnitude for strongly asymmetric couplings (Sec.~\ref{sec:oscillators}).

The relaxation time of this mode, $1/\abs{\lambda_2}\propto e^{N/(2\xi)}$, is exponentially long.
At fixed $\alpha$, Eq.~\eqref{eq:thirdrate} bounds the remaining rates away from zero uniformly
in even $N$. At $N=160$, $\abs{\lambda_3}$ lies within $5\%$ of the uniform-chain band edge
$a(1-\sqrt r)^2$ for $\alpha=0.3$ and $0.5$. The nonuniform center can keep the third rate
strictly below that edge for arbitrarily long chains \proofin{rem:subband}.
In linear response, identical kicks at mirror-related sites leave the same permanent shift but
excite the slow antisymmetric mode with opposite signs. After the faster modes subside, the
remaining phase difference between the halves identifies the kicked side. The predicted slope of $\ln\abs{\lambda_2}$ against $N$,
$\frac12\ln\cos2\alpha$, is $-0.1807$ at $\alpha=0.4$ and $-0.3078$ at $\alpha=0.5$; the measured slopes
are $-0.1807$ and $-0.3075$. With 160-digit arithmetic the constant $C(\alpha)$ is reproduced to
5--11 digits for $\alpha=0.2$--$0.6$, with $N$ up to $300$ ($\alpha=0.2$, $0.3$) or $120$ [Fig.~\ref{fig:geometry}(b)].

\paragraph*{Corners.}
Rectangular lattices with free boundaries are Cartesian products of chains. Ren and Ermentrout
showed that on such arrays, with nearest-neighbor coupling and frequencies that are sums of row and
column profiles, the locked state reduces to independent chain problems, two in two dimensions and
three in three~\cite[Lemma~3.3, Thm.~3.4 and Remark~4]{RenErmentrout1998}. The next proposition states
this reduction for general Cartesian products and adds the factorization of the Jacobian and of the
influence.

\begin{proposition}[Product law]\label{thm:product}
Take identical oscillators on $G\square H$ with the weights of each factor repeated on every
copy of that factor. Let $u$ and $v$ be positive locked states on $G$ and $H$ with
influences $w_G$ and $w_H$. Then $\theta_{(x,y)}=u_x+v_y$ is a positive locked state, its Jacobian is
the Kronecker sum $J_G\otimes I+I\otimes J_H$, and
\begin{equation}\label{eq:product}
  w=w_G\otimes w_H .
\end{equation}
If $u$ and $v$ are the branches continued from in-phase locking along a common parameter path, so is $\theta$.
\end{proposition}

The product of the two-source chains gives the positive in-phase-continued Sakaguchi state on
every finite box for $|\alpha|<\pi/4$ \proofin{thm:product}. Stable rotating waves can coexist with
this branch. Examples occur on even square arrays for certain odd couplings~\cite[Thm.~1]{PaulletErmentrout1994},
with persistence at small phase lag in the setting of Ref.~\cite[Sec.~3.1]{RenErmentrout1998}.
Corner localization refers to the
in-phase-continued branch. On a $D$-dimensional box, the exact finite-box influence is
$w(x)=\prod_{m=1}^D w^{(m)}(x_m)$. As every side grows at fixed $0<|\alpha|<\pi/4$, the geometric profiles
near the ends give corner localization, with decay $r^{\sum_{m=1}^D d_m(x)}$ at fixed distances $d_m$
from the nearest face along each axis. The mean distance to the boundary tends to
$r^D/(1-r^D)$, and the share in corner boxes of side $\ell$ tends to $(1-r^\ell)^D$
\proofin{cor:corner} [Fig.~\ref{fig:geometry}(c,d)]. On a $64\times64$ lattice at $\alpha=0.3$, the four
$6\times6$ corner boxes, $3.5\%$ of the nodes, hold $45\%$ of the influence, and Eq.~\eqref{eq:product}
holds to $4\times10^{-12}$ on lattices from $16^2$ to $64^2$, $8^3$ and $12^3$.

The selection mechanism is a competition of pacemakers. With $\Gamma(0)\ne0$, boundary oscillators
receive fewer inputs and the boundary acts as a pacemaker region~\cite{RenErmentrout1998}; the original
lattice Sakaguchi model subtracts the zero-lag value, setting $\Gamma(0)=0$ and removing this
effect~\cite{SakaguchiShinomotoKuramoto1988}. In two dimensions, the limiting frequency shift
$\Omega-\omega$ is $2a\Gamma(-x_0)$ for a diagonal wave from a corner and
$a\Gamma(-x_0)+2a\Gamma(0)$ for a wave from a straight edge. These shifts explain which boundary
feature entrains the lattice, as for competing pacemakers in disordered
media~\cite{BlasiusTonjes2005}. For Sakaguchi coupling the corner is faster for every
$0<\alpha<\pi/4$. The corner profile is a product of end profiles of a single stationary mode.
For the distinction from an intrinsic second-order skin effect~\cite{Kawabata2020}, see
Appendix~\ref{a:sec:product}.

\section{Stationary currents and kick reciprocity}\label{sec:topology}

In the free chain of Fig.~\ref{fig:concept}, end and center kicks leave very different
permanent shifts, although detailed balance makes the stationary current zero on every
edge. A uniformly twisted ring can have equal source weights and a stationary current
around the ring. Every site then leaves the same permanent shift, while exchanging the
kick and observation sites can reveal the current during relaxation. A nonzero reciprocity
defect for a single pair at one time suffices. Zero defect characterizes detailed balance
when it holds for every pair at every positive time.

\begin{figure}[tbp]
\includegraphics[width=\columnwidth]{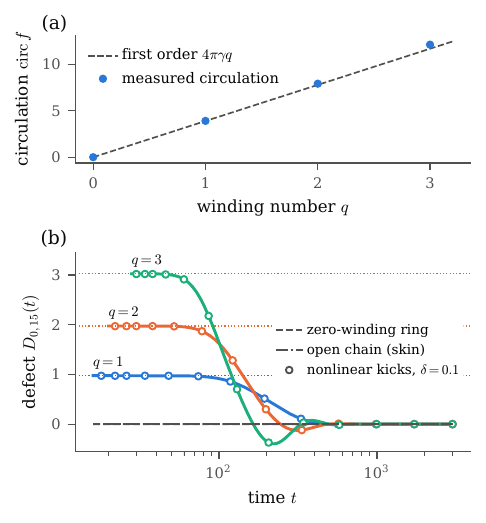}
\caption{\label{fig:topology}Topology, on rings of $N=60$ Sakaguchi oscillators with $\alpha=0.3$. (a) Circulation
of the gauge field around $q$-twisted states, computed from the locked state (markers), against the first-order
value $4\pi\gamma q$ (dashed); the exact value $NF(2\pi q/N)$ coincides with the markers. (b) Reciprocity defect $D_{0,15}(t)$ from the linear propagator (lines)
and from direct nonlinear kicks of size $\delta=0.1$ (open circles), on the window where both responses exceed
$10^{-3}$ of the kick; time is in units of $1/a$. Dotted, the plateau value $15F(2\pi q/N)$. The zero-winding ring has
a small defect and the open chain has zero defect in linear response, although the chain's raw
responses differ by factors up to $e^{2.9}$.}
\end{figure}

\paragraph*{Twisted rings.}
On a ring of $N\ge3$ identical oscillators with uniform weight $a$, the $q$-twisted state
$\theta_i=i\Delta$, $\Delta=2\pi q/N$, has positive rates $A_\pm=a\Gamma'(\pm\Delta)$ for steps
$i\to i\pm1$ in the regime considered here \proofin{prop:ring}. Its generator is circulant, so
the influence is exactly uniform. The cycle affinity is $NF(\Delta)$ and the current on each
oriented edge is $\jmath=(A_+-A_-)/N$. Both are nonzero exactly when $A_+\ne A_-$. When the
Taylor interval in Eq.~\eqref{eq:circ} applies, $NF(\Delta)=4\pi\gamma q+O(q^3/N^2)$.
For $N=60$ and Sakaguchi $\alpha=0.3$ the computed
circulations are $3.90$, $7.90$ and $12.10$ for $q=1,2,3$, against the first-order values $3.89$,
$7.77$ and $11.66$ [Fig.~\ref{fig:topology}(a)]. With unequal directional rates the relaxation spectrum has complex modes,
$\operatorname{Im}\lambda_k=(A_+-A_-)\sin(2\pi k/N)$, as noted before~\cite{TonjesBlasius2009a,Mihara2018,KimKori2025}:
the slowest Fourier modes rotate by $0.62$, $1.25$ and $1.92$ radians per $e$-fold for $q=1,2,3$.
Winding also labels electrical loop-flow states in the meshed-grid examples of
Ref.~\cite{Coletta2016}, including cases with moderate ohmic losses. Electrical power flow can
circulate in the lossless case, where the auxiliary walk is reversible. For cycles of ordinary
reciprocal lines at small lags the leading affinity of the walk is $\sum_{e\in C}2\gamma_e\Delta_e$
(Sec.~\ref{sec:networks}). With a common line gain it reduces to $4\pi\gamma q_C$ at this order.

\paragraph*{An observable for the flux.}
To detect stationary probability current, compare weak-kick responses with source and receiver
interchanged. Let $P_{ij}(t)=(e^{Jt})_{ij}$ be the phase displacement of node $i$ from its
unperturbed trajectory at time $t$, per unit kick at node $j$ in linear response. For $t>0$,
define the \emph{reciprocity defect}
\begin{equation}\label{eq:defect}
  D_{ij}(t):=\ln\frac{w_i\,P_{ij}(t)}{w_j\,P_{ji}(t)} .
\end{equation}
Its exponential is the stationary two-time detailed-balance ratio
of Ref.~\cite[Sec.~IV A and Fig.~1]{BenAvrahamDoroszPleimling2011}.

\begin{proposition}[Reciprocity]\label{thm:recip}
The following are equivalent: (i) the gauge field has zero circulation on every cycle; (ii) the walk
is in detailed balance; (iii) $D_{ij}(t)=0$ for all $i,j$ and all $t>0$. If they hold, the
relaxation spectrum is real.
\end{proposition}

The proof is standard reversibility~\cite{Kelly1979} \proofin{thm:recip}; what matters is what $D$
separates. Non-reciprocity of the raw responses, $P_{ij}\neq P_{ji}$, is generic and says little: in
a skin chain it is a pure gauge effect, absorbed by the potential, and $D\equiv0$ even where
$P_{ij}/P_{ji}$ is large. A nonzero defect certifies that the walk carries a circulating current,
that is, a cycle with nonzero flux. Time-reversal asymmetry of two-time joint probabilities is the
standard signature of broken detailed balance in stationary systems~\cite{Gnesotto2018}; $D$ is its
response-based form, measured with kicks rather than inferred from fluctuations. On a twisted ring the defect is explicit: from the Skellam form of
the propagator~\cite{Skellam1946}, $D_{0m}(t)\to mF(\Delta)$ as $t\downarrow0$ and $D_{0m}(t)\to0$ as $t\to\infty$, for
$0<m<N/2$ \proofin{prop:ring}. On the infinite line, the bare Skellam kernel for the unwrapped biased walk $X_t$ has ratio
$\mathbb P(X_t=m)/\mathbb P(X_t=-m)=(A_+/A_-)^m$ at every $t>0$. On the ring the stationary
weights cancel, and wrap images contribute at every positive time. The defect stays close to
$mF(\Delta)$ while those images are small relative to the unwrapped kernel. For $N=60$,
$\alpha=0.3$, and $m=15$, the finite-ring defect stays within $1\%$ of $15F(2\pi q/N)$ through
$at=50$ for $q=1,2,3$. It then departs from this plateau, with the $q=3$ defect about $46\%$
below the plateau at $at=100$.

Direct simulation of the nonlinear model tests the kick response [Fig.~\ref{fig:topology}(b) and
Table~\ref{tab:defect}]. For $N=60$,
$\alpha=0.3$ and a source--receiver distance $m=15$, the sampled linear maxima agree with the
plateau value $mF(2\pi q/N)$ for $q=1,2,3$; nonlinear kicks of size $10^{-3}$ reproduce the linear values
to $2\times10^{-4}$, and kicks of size $0.1$ to within $0.3\%$. In linear response and for small kicks,
the control states give $|D|\le1.4\times10^{-4}$ at the same distance, including an open skin chain with
$70\times$ influence contrast and raw response ratios up to $e^{2.9}$. At kick size $0.1$ the controls
rise to $5.9\times10^{-4}$ (ring) and $4.7\times10^{-3}$ (chain), still at least two hundred times below
the $q=1$ signal. Swapping source and receiver at
distance 15 on a $q=1$ ring changes the $w$-normalized response by a factor $e^{0.98}\approx2.7$.

\begin{table}[tbp]
\caption{\label{tab:defect}Maximal reciprocity defect $\abs{D_{0,15}}$ on rings and an open chain of $N=60$ Sakaguchi
oscillators ($\alpha=0.3$), from the linear propagator and from direct nonlinear kicks of size
$\delta$ over the window where both responses exceed $10^{-3}$ of the kick, against the plateau value $15F(\Delta)$, $\Delta=2\pi q/N$. Controls: a zero-winding ring with
heterogeneous frequencies (influence contrast $2.1$) and an open chain (contrast $70$).}
\footnotesize
\begin{ruledtabular}
\begin{tabular}{lcccc}
State & linear & $\delta=10^{-3}$ & $\delta=0.1$ & $15F(\Delta)$\\
\colrule
ring, $q=1$ & 0.9757 & 0.9757 & 0.9744 & 0.9757\\
ring, $q=2$ & 1.9754 & 1.9754 & 1.9725 & 1.9754\\
ring, $q=3$ & 3.0255 & 3.0254 & 3.0205 & 3.0255\\
$q=0$ ring, heterog. & $1.4\times10^{-4}$ & $1.3\times10^{-4}$ & $5.9\times10^{-4}$ & ---\\
open chain & $4.7\times10^{-15}$ & $4.8\times10^{-5}$ & $4.7\times10^{-3}$ & ---\\
\end{tabular}
\end{ruledtabular}
\end{table}

Operationally, the same two kick experiments supply both the permanent shifts ($w_i,w_j$) and the
transient responses ($P_{ij}(t),P_{ji}(t)$), yielding $D_{ij}(t)$ directly and complementing methods
that infer stationary probability flux from fluctuations~\cite{Battle2016}.

\section{Testing the phase reduction}\label{sec:oscillators}

For the three full oscillator models with non-odd coupling, the sign of $\gamma$ predicted
whether leaders or laggers carried greater memory. This held at both parameter values in
every configuration eligible for the sign comparison. The errors in individual memory
shares varied across models, graphs and disorder amplitudes.
The gain $\gamma=\Gamma''(0)/\Gamma'(0)$ follows by averaging the phase-response vector
against the coupling term over one cycle. We computed it for nineteen parameter sets of
eight models, with electrical, mechanical, chemical and fast synaptic channels. Table~\ref{tab:gamma}
gives a selection; the full catalog, coupling functions and checks are in the Supplemental
Material~\cite{SM}. The computation reproduces the Stuart--Landau closed forms to
$4\times10^{-12}$, and the adjoints agree with direct kicks to $2.3\times10^{-6}$.

\begin{table}[tbp]
\caption{\label{tab:gamma}Memory gain $\gamma=\Gamma''(0)/\Gamma'(0)$ for standard oscillator models under
diffusive (electrical, mechanical, or species-diffusion), input, and fast synaptic coupling. For $\gamma>0$ perturbation
memory favors phase leaders at small lags in an untwisted state, and $\gamma<0$ favors laggers. Rows marked ``no'' have
$\Gamma'(0)<0$ and lie outside the expansion about stable in-phase locking.}
\footnotesize
\begin{ruledtabular}
\begin{tabular}{llrc}
Model & coupling & $\gamma$ & $\Gamma'(0)>0$\\
\colrule
Stuart--Landau ($c$) & diffusive & $c$ & yes\\
van der Pol ($\mu=1,5$) & velocity & $0.40/4.44$ & yes\\
van der Pol ($\mu=1,5$) & position & $-4.30/{-2.87}$ & no\\
FitzHugh--Nagumo & electrical & $2.66$ & yes\\
Morris--Lecar I/II & electrical & $0.66/0.94$ & yes\\
Hodgkin--Huxley & electrical & $2.59$ & yes\\
FHN, ML, HH & exc.\ synapse & $5.5$--$9.7$ & yes\\
FHN, ML, HH & inh.\ synapse & $4.0$--$7.9$ & no\\
Wilson--Cowan & E$\to$I & $-0.36$ & yes\\
Goodwin & diffusive & $4.78$ & yes\\
Brusselator & act./inh. & $3.26/{-0.80}$ & yes\\
\end{tabular}
\end{ruledtabular}
\end{table}

Where the small-lag remainder is controlled, the contrast depends on both the gain and the
phase spread: $\osc\ln w\approx2|\gamma|\osc\tilde\theta$. The sign selects leaders or laggers.
The coupling channel matters: a single-harmonic fit reverses the sign of $\gamma$ for the
sampled electrically coupled FitzHugh--Nagumo, Hodgkin--Huxley and Morris--Lecar~II models. In some sampled
neuron couplings, the weak-skin prefactor underestimates the even-chain constant by up to
a factor of $19$~\cite{SM}.
Models with $\Gamma'(0)<0$ require checks at their actual locked lags.
In a position-coupled van der Pol chain, the cycle-averaged adjoint projection gives signed
influences from $-14$ to $+35$ times the mean. Position kicks at eight nodes check the
instantaneous adjoint response at one phase; the five resolved ones agree within $0.4\%$~\cite{SM}. The seven parameter sets outside the Stuart--Landau and van der Pol families were recomputed at a second
parameter value: $\gamma$ and $\Gamma'(0)$ keep their signs in every physical channel, while
$\gamma$ changes by up to $46\%$~\cite{SM}.

\paragraph*{Full networks.}
We compare four full oscillator models with their phase reductions: Morris--Lecar type-I
and Hodgkin--Huxley neurons with gap junctions, inhibitor-coupled Brusselators ($\gamma=-0.80$),
and Stuart--Landau oscillators without shear as an odd control. Each is placed on a $40$-node
chain and on $100$-node Barab\'asi--Albert and random geometric graphs. The coupling levels are
$\epsilon=K\Gamma'(0)/\Omega_0\in\{0.001,0.003,0.01,0.03\}$, with $\Omega_0$ the nominal
angular frequency. Node parameters realize a frequency spread of half the positive-rate
continuation limit, scaled with coupling to leave the leading phase model and its weights
$w$ unchanged. Deviations therefore combine finite-coupling and node-parameter effects.
Of $48$ designs, $46$ have recorded phase locks and $45$ pass
the full-model checks (the three excluded $\epsilon=0.03$ designs are detailed below).

In the full network, $v_i=T^{-1}\int_0^Tz_i\cdot F_i\,dt$ is the response to a local displacement
along the $T$-periodic locked orbit's tangent at node $i$, averaged over the starting phase. Here $z$ is the periodic adjoint
normalized by $z\cdot F=1$, $F$ is the network vector field, and $z_i,F_i$ denote their node-$i$ components. Then $\sum_iv_i=1$, and $v=w$ for a phase
model. To check the adjoint, we kick the first state component in both directions at
$130$ sources across $13$ configurations and compare the asymptotic
time shifts with the corresponding instantaneous adjoint components at the kick phase. The comparison gives $0$
violations of the prescribed tolerance~\cite{SM}. The adjoint check and the two results listed below were
predictions fixed before the run; the coefficients and crossing levels reported after them are exploratory. Usable configurations have a feasible parameter
design, a transversely stable locked orbit and an adjoint that passes the normalization and
periodicity checks. The sign comparison also requires $v>0$.
\begin{itemize}
\item \emph{The two smallest couplings give near-linear error scaling.} Between $\epsilon=0.001$
  and $0.003$, the error $L_1=\sum_i\abs{v_i-w_i}$ grows with a fitted exponent of
  $0.98$ to $1.03$ over $12$ model--graph pairs.
  At $\epsilon=0.001$, $L_1\le0.018$.
\item \emph{The phase ranking persists in the sampled full models.} In the $33$ usable configurations of the
  three non-odd models, the Spearman correlation $\rho$ of $\ln v$ with the node phases of the full
  network has the sign of $\gamma$, with minimum magnitude about $0.986$. Morris--Lecar and Hodgkin--Huxley
  networks remember their leaders, and Brusselator networks remember their laggers.
\end{itemize}

The leading error coefficient $L_1/\epsilon$ varies by a factor of about
$135$ among the non-odd models, and by $5\times10^{4}$ when the odd control is included (Table~\ref{tab:fullnet}). For nonnegative normalized influences,
$L_1/2$ is the misplaced share of memory. Log--log interpolation puts a $5\%$ misplaced share
at $\epsilon=0.0058$ for Hodgkin--Huxley on the scale-free graph and
$\epsilon=0.013$ for Brusselators on the geometric graph; the other ten pairs
remain below it over their usable ranges.

Three designs at $\epsilon=0.03$ are unusable. One Hodgkin--Huxley design exceeded the
frequency-spread limit. In another, the isolated cycle of one cell had a nontrivial Floquet
multiplier of $0.739$, above the protocol cutoff of $0.5$ that marks the edge of the oscillatory range. The full Brusselator network on
the scale-free graph drifted although its phase model locked.
The three non-odd models are repeated at a second parameter value (Morris--Lecar $I=50$,
Brusselator $B=2.5$, Hodgkin--Huxley $I=15$; lower rows of Table~\ref{tab:fullnet}). The three results
above hold there as well: the adjoint check has no violation in $100$ sources, the fitted exponents lie
between $0.98$ and $1.06$ in $9$ model--graph pairs, and $\rho$ has the sign of $\gamma$, with magnitude at least $0.985$, in all $29$ of the $34$
usable configurations where $v>0$. In the other five, all on the chain, $v$ has entries at round-off level
(magnitude below $2\times10^{-13}$) that are not positive, and the comparison is not made. In those five the rank
correlation with the node phases still has the sign of $\gamma$: it is $0.93$ to $0.999$ for $\ln\abs v$ and $0.68$ to $1.00$
for $v$ itself.
Between the two parameter values, $L_1/\epsilon$ changes by factors of $0.5$ to $6.9$. The largest
change occurs on the Hodgkin--Huxley chain, where the frequency spread in units of $K\Gamma'(0)$ also rises
from $0.0648$ to $1.875$ under the half-range rule. The comparison therefore includes changes
in disorder as well as in the node parameter. Two $\epsilon=0.03$ designs are unusable, both
of them also unusable at the first value.
Each model uses one coupling channel on one realization per graph; full design and exclusion details are
in the Supplemental Material~\cite{SM}.

\begin{table}[tbp]
\caption{\label{tab:fullnet}Comparison of full oscillator networks with their phase reductions. Coefficient
$L_1/\epsilon$ of the phase-reduction error $L_1=\sum_i\abs{v_i-w_i}$ at $\epsilon=0.001$ for Stuart--Landau
(SL), Morris--Lecar type-I (ML I) and Hodgkin--Huxley (HH) neurons with gap junctions and Brusselators
coupled through the inhibitor, on the chain
($N=40$), the Barab\'asi--Albert graph (BA) and the random geometric graph (RGG) ($N=100$). In parentheses, the
coupling estimated by log--log interpolation where $L_1$ crosses $0.1$ within the usable range. $\rho$: the value of
$\rho(\ln v,\theta)$ closest to zero over the usable configurations of the model with $v>0$. The lower rows repeat the three non-odd models at a second parameter value ($I=50$, $B=2.5$ and $I=15$).}
\footnotesize
\begin{ruledtabular}
\begin{tabular}{lrlllr}
Model & $\gamma$ & chain & BA & RGG & $\rho$\\
\colrule
SL ($c=0$) & $0$ & $0.00034$ & $0.07$ & $0.048$ & -- \\
ML I & $0.66$ & $0.19$ & $0.15$ & $0.76$ & $0.999$ \\
Brusselator & $-0.80$ & $0.13$ & $4.5$ & $5.5$ ($0.013$) & $-0.998$ \\
HH & $2.59$ & $0.15$ & $17.6$ ($0.0058$) & $7.5$ & $0.986$ \\
\colrule
ML I, $50$ & $0.68$ & $0.28$ & $0.26$ & $0.39$ & $0.994$ \\
Bruss., $2.5$ & $-0.92$ & $0.46$ & $8.8$ & $12.9$ ($0.0061$) & $-0.995$ \\
HH, $15$ & $2.36$ & $1.0$ & $18.6$ ($0.0053$) & $6.0$ ($0.019$) & $0.985$ \\
\end{tabular}
\end{ruledtabular}
\end{table}

\section{Tests on empirical networks}\label{sec:networks}

In the two lossy grid models, retaining the bias of each line improves the prediction
of bus memory at nominal loading. With equal damping and uniform balancing, the relative
memory at a bus also equals one minus its marginal loss sensitivity. On measured
connectome wiring, assigned oscillator dynamics test the weight approximation in states
with much stronger localization.

\paragraph*{Line-dependent rates.}
With the voltage magnitudes $v_i$ fixed at their values in the case files and equal damping
coefficients $D_i\equiv D$, the first-order (overdamped) grid model is
\begin{equation}\label{eq:grid}
  \dot\theta_i=\lambda P_i-P_i^{\rm el}(\theta)=\omega_i+\sum_{j\ne i}A_{ij}\sin(\theta_j-\theta_i+\beta_{ij}).
\end{equation}
Here $t$ is physical time divided by $D$, $P_i$ is the net injection (generation minus load),
$\lambda$ scales the loading, and
$P_i^{\rm el}=\operatorname{Re}[V_i\overline{(YV)_i}]$ is the active power leaving bus $i$, with
$V_i=v_ie^{\mathrm i\theta_i}$ and $Y$ the bus admittance matrix. This gives $A_{ij}=v_iv_j\abs{Y_{ij}}$,
$\beta_{ij}=\arg Y_{ij}-\pi/2$ and $\omega_i=\lambda P_i-v_i^2\operatorname{Re}Y_{ii}$. This is the standard
phase-oscillator form of lossy grid models~\cite{DorflerBullo2012,NishikawaMotter2015,Hellmann2020};
$\beta_{ij}$ is the phase shift of Ref.~\cite{NishikawaMotter2015}, with the same sign, and
Eq.~\eqref{eq:grid} is their model without inertia and with equal damping. For an ordinary reciprocal inductive line
of resistance $r$ and reactance $x$, $\beta=\arctan(r/x)$. In the convention of Eq.~\eqref{eq:model}
the gain of this Sakaguchi coupling is
\begin{equation}\label{eq:gammaline}
  \gamma_e=-r_e/x_e .
\end{equation}
Off-nominal taps and phase shifters can alter the amplitudes and lags used in the directional
rates. In a lossless grid the rates are symmetric even in the presence of taps and phase shifters,
so the memory is uniform.

Proposition~\ref{thm:hodge} and Theorem~\ref{thm:weighted} use only the effective rates and apply when the
coupling varies between lines or directions. On ordinary reciprocal lines the small-lag bias is
$f_e\approx2\gamma_e\Delta_e$. Different gains can give nonzero cycle affinity already at first
order, even in an untwisted state. The sign of each line gain is therefore insufficient to order
the buses by phase. Direction-dependent lags can also give a bias at zero phase difference.
For small biases the leading potential $\Phi_{\rm lin}$ is the weighted least-squares potential
of the actual $f$, with fixed conductances $S$: it satisfies
$\sum_{i\sim j}S_{ij}\,[f_{ij}-(\Phi_{{\rm lin},j}-\Phi_{{\rm lin},i})]=0$ for all $j$.

\paragraph*{Memory and marginal loss sensitivity.}
Differentiating the locked equations with respect to $\omega_k$ along the local branch and
multiplying by $w^\trans$ gives $\partial\Omega/\partial\omega_k=w_k$. This is the local frequency
response in Refs.~\cite{TonjesBlasius2009b,Kori2009,Skardal2016}. Summing Eq.~\eqref{eq:grid} over the buses gives
$N\Omega=\lambda\sum_iP_i-P_{\rm loss}$, where $P_{\rm loss}=\sum_iP^{\rm el}_i$ is the total active loss, including consumption in shunt conductances; this
frequency--loss balance appears in Ref.~\cite{Hellmann2020}.
Hence, at $\lambda=1$,
\begin{equation}\label{eq:mlf}
  Nw_k=1-\frac{\partial P_{\rm loss}}{\partial P_k}.
\end{equation}
Here $Nw_k$ is the influence relative to its uniform value, and the derivative is the marginal
loss sensitivity with distributed slack~\cite{Litvinov2004}. In Eq.~\eqref{eq:grid} every bus
has the same slack weight, because the common frequency absorbs the imbalance uniformly
\cite{Dhople2020}. Left null weights already describe how local fluctuations excite lossy
grids~\cite{Plietzsch2022}. A related static relation at a marginal power-flow state connects
left-null multipliers with incremental losses~\cite[Eq.~(19)]{Ayuev2015}. Equation~\eqref{eq:mlf}
connects that type of sensitivity to dynamical phase-kick memory under the balancing convention
of Eq.~\eqref{eq:grid}. Within this model, a bus where incremental injection reduces losses
holds more than the uniform share. This identification of relative phase-kick memory with marginal loss sensitivity requires
equal damping. The Supplemental Material derives the weights for unequal damping and the
distinct responses to phase and frequency kicks in the swing model, with an IEEE~118 check~\cite{SM}.

\begin{figure*}[t]
\includegraphics[width=\textwidth]{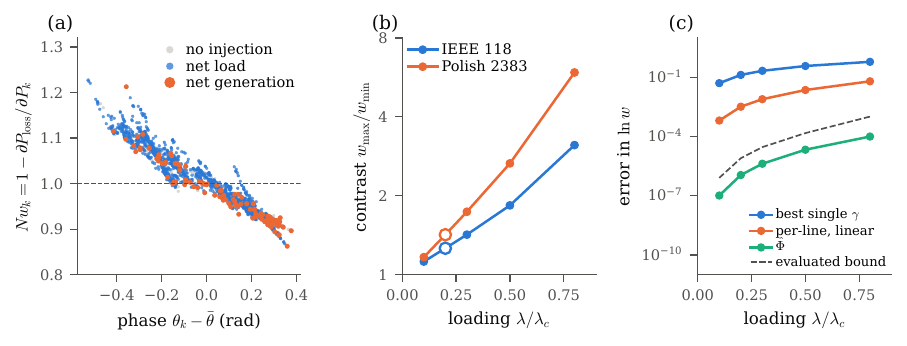}
\caption{\label{fig:networks}Real networks. (a) Polish transmission grid (2383 buses) at its operating point:
memory $Nw_k$, which by Eq.~\eqref{eq:mlf} is one minus the marginal loss sensitivity of bus $k$, against the
phase. The phase laggers hold the most memory, and the ten most influential buses are net loads except one
$0.7$\,MW generator; dashed line, uniform memory. (b) Memory contrast along the loading series of both grids; open markers mark the
operating point $\lambda=1$. The operational loading reference is $\lambda_c\approx5.0$ for both grids. (c) Error in $\ln w$ along the loading series of the Polish grid for the best
single-gain law $2\gamma^*\theta$ (with $\gamma^*$ fitted to each state), the per-line linear law
$\Phi_{\rm lin}$ and the self-consistent potential $\hat\Phi$. Weighted-potential errors are shown
above the accuracy gate, with finite evaluations of Eq.~\eqref{eq:certified} dashed.}
\end{figure*}

\paragraph*{Two grids.}
We took the IEEE 118-bus test case and the Polish transmission grid at winter peak (2383 buses), both as
distributed with MATPOWER~\cite{Zimmerman2011}, in the copies of the Netzschleuder
repository~\cite{Peixoto2020}. The models reproduce the case-file power flows; Table~\ref{tab:networks} summarizes
the quantitative checks. The line-gain distributions, operating-point checks and exclusion
of IEEE~300 are documented in the Supplemental Material~\cite{SM}. With the locking procedure of the protocol the
Polish case did not lock at its operating point, so grid states are located by continuation from the lossless grid,
a logged deviation; on IEEE~118 the two procedures give the same states to $10^{-12}$\,rad.

The per-line law $\Phi_{\rm lin}$ is about $56$ and $41$ times more accurate than the best single-gain
law $2\gamma^*\theta$, even with $\gamma^*$ fitted to each state. Its log-weight errors are
$1.5\times10^{-3}$ and $3.3\times10^{-3}$, against $0.082$ and $0.13$. The self-consistent
potential $\hat\Phi$ has measured errors of $9.92\times10^{-7}$ and $1.12\times10^{-6}$.
At the operating points the double-precision bound is at most $7.2$ times the measured error
[Fig.~\ref{fig:networks}(c) shows the loading series]. The stationary-current fraction
$\kappa=\sum_e\abs{\jmath_e}/\sum_e(w_iW_{ij}+w_jW_{ji})$ is $105$ times the surrogate value
on IEEE~118 and $73$ times that value on the Polish grid. The surrogate sets $r/x$ to the
grid median on inductive branches and sets transformer phase shifts to zero, preserving
impedance magnitudes and tap ratios. The absolute
fractions are small ($\kappa=1.3\times10^{-3}$ and $2.9\times10^{-5}$), as expected at small lags.

The memory is nearly uniform at nominal injections. The share relative to uniform memory, $Nw_k$, ranges from $0.87$
to $1.09$ on IEEE 118 and from $0.86$ to $1.23$ on the Polish grid. At the five most and
five least influential buses of each grid, finite differences of the common frequency with
respect to injection agree with $w_k$ to relative error $8\times10^{-8}$. Power balance relates
this sensitivity to losses through Eq.~\eqref{eq:mlf}. The lossless control is uniform to $2\times10^{-12}$ in
log-weight oscillation. The lagging buses hold more memory in these two cases
[Fig.~\ref{fig:networks}(a) and Table~\ref{tab:networks}].

Scaling all injections from $0.1\lambda_c$ to $0.8\lambda_c$, with $\lambda_c$ the operational loading reference defined below, raises $w_{\max}/w_{\min}$ from
$1.1$ to $3.1$ (IEEE 118) and from $1.2$ to $5.9$ (Polish) [Fig.~\ref{fig:networks}(b)].
All $18$
sampled grid states (four loadings, four loss multipliers and one uniform-$r/x$ surrogate per grid)
lock with positive rates. In the $16$ accuracy-eligible states,
$\hat\Phi$ stays within about $1.89\times10^{-4}$ of $\ln w$ in oscillation, and
$\osc(\ln w-\hat\Phi)/\eta^3\le0.027$.
At the last converged states on IEEE~118 ($\lambda=5.030$) and the Polish grid
($\lambda=5.034$), every transverse mode examined in the stability screen decayed.
These continuation endpoints are
consistent with nearby folds. The reference $\lambda_c$ is found by doubling the loading and
bisecting to $2\%$, which returns $5.0$ for both endpoints; it is not an imposed cap. The fixed-voltage model scales all injections together
and omits operating constraints and reactive-power dynamics.

\paragraph*{\emph{C.~elegans} connectome.}
The 460-cell \emph{C.~elegans} gap-junction graph~\cite{Cook2019} supplies a further test
with assigned electrical couplings and synthetic frequencies. All $85$
finite double-precision bound evaluations exceed the measured error, and all
$14$ qualifying small-$\eta$ slopes lie in $[2.7,3.3]$
(Table~\ref{tab:networks}). The Supplemental Material gives the complete connectome analysis and numerical details~\cite{SM}.

\begin{table}[tbp]
\caption{\label{tab:networks}Checks on empirical networks. The criteria in the second column were fixed before these runs; the populations of the bound, constant and third-order rows were selected afterwards. Frequency sensitivity compares $\partial\Omega/\partial P_k$ with $w_k$. The frequency-sensitivity, lossless, sign-rule, per-line law, and current rows report values at the nominal operating point $\lambda=1$ (the per-line and sign-rule criteria also pass across the lossy loading sweep); the bound, constant, and third-order rows use the precision-resolved stable subset ($16$ of $18$ grid states and $107$ of $120$ connectome states). The bound row gives passing comparisons among finite double-precision evaluations of Eq.~\eqref{eq:certified}. For $22$ further eligible connectome states, the path-length test gives no finite bound. The sign-rule row gives the Spearman correlation of $\ln w$ with $\theta$ and the excess median memory of loads over generators. $e=\osc(\ln w-\hat\Phi)$; $e_{\rm lin}$ and $e_{\rm uni}$ are the errors of $\Phi_{\rm lin}$ and of the best single-gain law $2\gamma^*\theta$; $\kappa_{\rm s}$ is the stationary-current fraction of the uniform-$r/x$ surrogate.}
\footnotesize
\begin{ruledtabular}
\begin{tabular}{lllll}
Check & criterion & IEEE 118 & Polish & \emph{C.~elegans}\\
\colrule
Freq. sens. & rel.\ err.\ $\le10^{-5}$ & $8\times10^{-9}$ & $8\times10^{-8}$ & --\\
Lossless & $\osc\ln w\le10^{-10}$ & $10^{-14}$ & $2\times10^{-12}$ & --\\
Bound & bound $\ge$ error & $8/8$ & $8/8$ & $85/85$\\
Constant & $e/\eta^3<1/24$ & $\le0.015$ & $\le0.027$ & $\le0.032$\\
Sign rule & $\le-0.5$; $>0$ & $-0.90$; $7\%$ & $-0.97$; $4\%$ & --\\
Per-line law & $e_{\rm lin}/e_{\rm uni}\le0.5$ & $0.018$ & $0.024$ & --\\
Current & $\kappa/\kappa_{\rm s}\ge10$ & $105$ & $73$ & --\\
Third order & slope $\in[2.7,3.3]$ & -- & -- & $14/14$ pairs\\
\end{tabular}
\end{ruledtabular}
\end{table}

\section{Consequences and conclusions}\label{sec:outlook}

For positive effective rates, a synchronized state's phase differences act as an emergent drift field for the linearized random walk that governs phase relaxation. The potential component of this drift shapes the stationary distribution $w$, fixing how much of a local kick the collective phase remembers, while its cycle affinities drive persistent probability currents. For untwisted states with small edge lags and positive coupling slope at zero, the lifted phases directly predict the log weights at leading order, while an exact self-consistent solution at fixed graph size extends this prediction to cubic accuracy in the edge residual alongside a computable stationarity-residual bound. On the studied chain and product-lattice branches where open boundaries induce a constant asymptotic phase gradient, the same drift localizes source influence at chain ends and box corners, creating a central conductance bottleneck that separates the halves of a long even chain by an exponentially long relaxation time. Conversely, when phase winding sustains a circulation around a cycle, the stationary weights can remain uniform while two-kick reciprocity breaks down.
The stationary-current representation of Proposition~\ref{thm:hodge} extends to asymmetric
edge weights $a_{ij}\neq a_{ji}$ when both directional rates are positive, with
$\ln(a_{ij}/a_{ji})$ adding to the gauge field $f_{ij}$.

\paragraph*{Noise, finite kicks, and signed rates.}
Localization also changes collective noise. For weak independent phase noise of intensity
$D$ at each oscillator, $D_{\rm coll}=D\sum_iw_i^2$~\cite{Masuda2010}. With $\Phi=\ln w$,
\begin{equation}\label{eq:noise}
  D_{\rm coll}=D\,\frac{\sum_ie^{2\Phi_i}}{\big(\sum_ie^{\Phi_i}\big)^2}=\frac{D}{N}\cdot\frac{N}{\rm PR},
\end{equation}
where ${\rm PR}=1/\sum_iw_i^2$ is the participation ratio. The penalty $N/{\rm PR}$ is small on the sampled small-world
graphs and grows with size on two-dimensional families~\cite{SM}. The linear-response weights
need not rank finite-kick shifts, as the supplementary scans show. Negative effective rates
can also give signed influence; their sup-norm growth exponent is bounded by twice the
largest row sum of the negative part of $W$ \proofin{prop:signed}.

\paragraph*{Experimental test.}
In electrochemical oscillator arrays~\cite{Kiss2002}, calibrated weak pulses at two nodes,
with source and receiver interchanged, could measure both the permanent shares and the
reciprocity defect. A reference trajectory separates transient from permanent shifts, and
decreasing pulse amplitudes verifies linear response. On a uniform twisted ring, measurements on
the plateau of Sec.~\ref{sec:topology}, where wrap contributions are small and both responses clear
the measurement floor, probe the cycle affinity directly, with an open reversible chain serving as
a zero control (whereas an untwisted array with a common coupling can still retain a cubic defect at
small lags).

\paragraph*{Open questions.}
It remains open whether stationary currents circulate around spiral cores with nonzero cycle
affinity. In the three full oscillator models with non-odd coupling, the reduced model
predicts whether phase leaders or laggers carry more memory in every eligible comparison,
but the error in the predicted shares varies strongly across graphs. What properties of the
oscillators and their coupling graph determine the first correction to those shares?

\section*{Data availability}

The reproducibility archive contains the simulation and analysis code, numerical datasets,
and per-file SHA-256 manifests. Supplemental Material~\cite{SM} documents the numerical
methods, model definitions and validation checks. Power-grid test cases~\cite{Zimmerman2011} and
the \emph{C.~elegans} gap-junction network~\cite{Cook2019} were obtained from the Netzschleuder
repository~\cite{Peixoto2020}, with download hashes recorded.
The reproducibility archive is publicly deposited on Zenodo at
\href{https://doi.org/10.5281/zenodo.24280237}{10.5281/zenodo.24280237}.

\begin{acknowledgments}
The author is a member of the Gruppo Nazionale per il Calcolo Scientifico (GNCS) of the Istituto Nazionale di Alta Matematica (INdAM).
\end{acknowledgments}

\section*{AI Use Disclosure}
The author used large language models (Claude, Anthropic, Fable 5.2; Gemini, Google DeepMind,
Neon) throughout this work. They wrote and reviewed most of the simulation, analysis and
figure code, which the author then ran on a laptop and on an HPC cluster, and they drafted the
protocols and the logs of later changes under the author's instructions. They were also used to
check the mathematical proofs, to compare the reported numbers with the stored results, and to draft and revise parts of the text to improve the presentation.
The research question, the main results and the proof strategies are the author's. The author
directed each step, ran and inspected the numerical results, read and checked every proof, and
takes full responsibility for the content. The reproducibility archive contains the code, the
dated protocols and the logs.

\appendix

\section{Setting and notation}\label{a:sec:setting}

The appendices prove the statements of the main text. They are self-contained, and Appendix~\ref{a:sec:setting}
fixes the notation. Numerical checks follow the statements. The labels in brackets (C1--C9, L1--L5, S1--S5, T34,
K2 and K3) are those of the corresponding scripts and records in the reproducibility archive.

Let $G=(V,E)$ be a finite connected undirected graph with $N=\abs V\ge2$ nodes and symmetric weights
$a_{ij}=a_{ji}>0$ for $\{i,j\}\in E$ ($a_{ij}=0$ otherwise). Coupling strengths are absorbed into
$a$. The phase model is
\begin{equation}\label{a:eq:model}
  \dot\theta_i=\omega_i+\sum_{j}a_{ij}\,\Gamma(\theta_j-\theta_i),\qquad i\in V,
\end{equation}
with $\Gamma\in C^\infty(\R)$ $2\pi$-periodic. A \emph{locked state} is a solution
$\theta(t)=\theta^*+\Omega t\,\one$. Write $\Delta_{ij}=\theta^*_j-\theta^*_i$ (principal values in
$(-\pi,\pi)$ when we speak of winding; see below). Linearizing \eqref{a:eq:model} about the locked state gives
$\dot{\delta\theta}=J\,\delta\theta$ with
\begin{equation}\label{a:eq:J}
  J=W-\diag(W\one),\qquad W_{ij}=a_{ij}\Gamma'(\Delta_{ij})\quad(i\neq j).
\end{equation}
$J$ has zero row sums. If $\Gamma$ is odd up to an additive constant, $\Gamma'$ is even and $W$ is
symmetric.

\begin{appdefn}[Gauge field]\label[appdefn]{a:def:gauge}
Let $F(x):=\ln\Gamma'(x)-\ln\Gamma'(-x)$, defined where $\Gamma'(\pm x)>0$. It is odd on this domain.
When $\Gamma'(0)>0$, $F'(0)=2\gamma$ with
\begin{equation}\label{a:eq:gamma}
  \gamma:=\frac{\Gamma''(0)}{\Gamma'(0)}.
\end{equation}
For a locked state, the antisymmetric edge function $f_{ij}:=\ln(W_{ij}/W_{ji})=F(\Delta_{ij})$ is the
\emph{gauge field}. The weight $a_{ij}$ cancels. Its \emph{circulation} around an oriented cycle
$C=(i_0,i_1,\dots,i_L=i_0)$ is $\circu_C f=\sum_{k}f_{i_ki_{k+1}}$. With $g:=\Gamma''/\Gamma'$ one
has $F'(x)=g(x)+g(-x)$ and $F'''(0)=2g''(0)$.
\end{appdefn}

For Sakaguchi coupling $\Gamma(x)=\sin(x-\alpha)$ the following hold:
\begin{itemize}[nosep]
\item $\Gamma'(x)=\cos(x-\alpha)$ and $\gamma=\tan\alpha$;
\item $F(x)=\ln\cos(x-\alpha)-\ln\cos(x+\alpha)$;
\item $F'''(0)=4\sec^2\!\alpha\tan\alpha$.
\end{itemize}

For the diffusive-type family $\Gamma(x)=\sin x+c(1-\cos x)$, one has $\gamma=c$.

\begin{appassum}[Positive locked state]\label[appassum]{a:ass:P}
$W_{ij}>0$ for every edge $\{i,j\}\in E$ (in both directions).
\end{appassum}

\begin{appdefn}[Influence]\label[appdefn]{a:def:w}
Under \cref{a:ass:P}, $w$ denotes the unique vector with $w^\trans J=0$, $w>0$, $\one^\trans w=1$
(\cref{a:prop:generator}).
\end{appdefn}

For winding we exclude antipodal edge lags and use principal values $\Delta\in(-\pi,\pi)$, so
$\Delta_{ji}=-\Delta_{ij}$. 

The \emph{winding} of an oriented cycle is
$q_C:=\frac1{2\pi}\sum_k\Delta_{i_ki_{k+1}}\in\Z$. A locked state is \emph{untwisted} if $q_C=0$
for every cycle. Equivalently, there is a real lift $\tilde\theta$ with
$\tilde\theta_j-\tilde\theta_i=\Delta_{ij}$ on every edge.

\section{Generator structure, influence and the maximum principle}\label{a:sec:generator}

\begin{appprop}\label[appprop]{a:prop:generator}
Under \cref{a:ass:P}, $J$ is the generator of an irreducible continuous-time Markov chain on $V$ that
jumps $i\to j$ at rate $W_{ij}$. Hence:
\begin{enumerate}[label=\textup{(\alph*)}]
\item $0$ is a simple eigenvalue of $J$, and every other eigenvalue has negative real part. The
locked state is therefore linearly stable modulo the phase shift.
\item There is a unique $w>0$ with $w^\trans J=0$ and $\one^\trans w=1$.
\item $P(t):=e^{Jt}$ is entrywise nonnegative and row-stochastic for $t\ge0$. Consequently, for
$\delta\theta(t)=P(t)\delta\theta(0)$, $\max_i\delta\theta_i(t)$ is nonincreasing,
$\min_i\delta\theta_i(t)$ is nondecreasing and $\norm{\delta\theta(t)}_\infty\le\norm{\delta\theta(0)}_\infty$
(maximum principle).
\item $w^\trans\delta\theta(t)$ is conserved, and $\delta\theta(t)\to(w^\trans\delta\theta(0))\,\one$.
\end{enumerate}
\end{appprop}

\begin{proof}
The off-diagonal entries of $J$ are $W_{ij}\ge0$, positive exactly on the edges of the connected
graph $G$, and $J\one=0$.

(a) By Gershgorin, every eigenvalue lies in a disk centered at $J_{ii}=-d_i$ with radius $d_i$, where
$d_i=\sum_jW_{ij}>0$. These disks lie in $\{\operatorname{Re}\lambda\le0\}$ and meet the imaginary
axis only at $0$. For $s>\max_id_i$, the matrix $sI+J$ is nonnegative and irreducible, with
eigenvector $\one$ for the eigenvalue $s$. By Perron--Frobenius, $s$ is its spectral radius and is
simple. Hence $0$ is simple for $J$, and the same theorem gives a positive left eigenvector, which
is (b).

(c) $e^{Jt}=e^{-st}e^{(sI+J)t}\ge0$ entrywise, and $e^{Jt}\one=\one$. Each $\delta\theta_i(t)$ is
therefore a convex combination of the entries of $\delta\theta(0)$.

(d) $\frac{d}{dt}w^\trans\delta\theta=w^\trans J\delta\theta=0$. By (a), $e^{Jt}\to\one w^\trans$.
\end{proof}

\noindent Part (a) recovers the positive-rate case of Ermentrout's stability
theorem~\cite{Ermentrout1992,RenErmentrout1998}. The proof also covers nonnegative rates when the directed
support of the positive rates is strongly connected.

\begin{appprop}\label[appprop]{a:prop:shift}
Assume \cref{a:ass:P}. Kick node $s$ of the locked state by $\delta$. The trajectory converges to the
locked orbit shifted by $w_s\delta+O(\delta^2)$.
\end{appprop}

\begin{proof}
In the co-rotating frame, $\mathcal C=\{\theta^*+c\one\}$ is a circle of equilibria. By
\cref{a:prop:generator}(a) it is normally attracting, so every nearby point $\theta$ has a unique
asymptotic phase $\Psi(\theta)$, and $\Psi$ is smooth locally (stable foliation of the normally
attracting circle of equilibria; isochrons
\cite{Guckenheimer1975}). Because $\Psi$ is constant along trajectories, $D\Psi(\theta^*)J=0$.
Because $\Psi(\theta^*+c\one)=c$, $D\Psi(\theta^*)\one=1$. Hence $D\Psi(\theta^*)=w^\trans$. This is
the weight vector of \cite{Kori2009}.
\end{proof}

\begin{appprop}\label[appprop]{a:prop:signed}
Without \cref{a:ass:P}, write $W=W^+-W^-$ with $W^\pm\ge0$, and set
$\nu:=\max_i\sum_jW^-_{ij}$. Then $\norm{e^{Jt}}_\infty\le e^{2\nu t}$.
\end{appprop}

\begin{proof}
Write $J=J^++B$, where $J^+:=W^+-\diag(W^+\one)$ is a generator, so $\norm{e^{J^+t}}_\infty=1$, and
$B:=-W^-+\diag(W^-\one)$, so $\norm B_\infty=2\nu$. By Duhamel's formula,
$e^{Jt}=e^{J^+t}+\int_0^te^{J^+(t-s)}Be^{Js}\,ds$. Gronwall's inequality applied to
$\norm{e^{Js}}_\infty$ gives the bound.
\end{proof}
\noindent\numcheck[C9]{On a random signed generator, $\norm{e^{Jt}}_\infty\le e^{2.54\,t}$ at all 26
sampled $t\in[0,5]$; the largest norm is $1.48$.}

\section{The gauge law}\label{a:sec:gauge}

Throughout this section \cref{a:ass:P} holds. For each edge put $S_{ij}:=\sqrt{W_{ij}W_{ji}}$, which is
symmetric, so that $W_{ij}=S_{ij}e^{f_{ij}/2}$. For $\Phi\in\R^V$ define the \emph{residual}
$h_{ij}(\Phi):=f_{ij}-(\Phi_j-\Phi_i)$ and the \emph{conductance} $c_{ij}(\Phi):=S_{ij}e^{(\Phi_i+\Phi_j)/2}$.

\begin{appprop}\label[appprop]{a:thm:hodge}
The potential $\Phi=\ln w$ satisfies
\begin{equation}\label{a:eq:H}\tag{H}
  \sum_{i\sim j}c_{ij}(\Phi)\,\sinh \big(\tfrac12h_{ij}(\Phi)\big)=0\qquad\text{for all }j\in V.
\end{equation}
Every solution of \eqref{a:eq:H} equals $\ln w$ up to an additive constant; the stationary
probability current on each edge is
\begin{equation}\label{a:eq:current}
  \jmath_{ij}:=w_iW_{ij}-w_jW_{ji}=2\,c_{ij}(\ln w)\,\sinh\!\big(\tfrac12h_{ij}(\ln w)\big).
\end{equation}
\end{appprop}

\begin{proof}
Let $v=e^\Phi$ with $\Phi$ arbitrary. Then
\begin{align*}
v_iW_{ij}-v_jW_{ji}&=S_{ij}\big(e^{\Phi_i+f_{ij}/2}-e^{\Phi_j-f_{ij}/2}\big)\\
&=S_{ij}e^{(\Phi_i+\Phi_j)/2}\big(e^{h_{ij}/2}-e^{-h_{ij}/2}\big)\\
&=2c_{ij}\sinh(h_{ij}/2).
\end{align*}
Moreover $(v^\trans J)_j=\sum_i(v_iW_{ij}-v_jW_{ji})$. So $v^\trans J=0$ if and only if \eqref{a:eq:H}
holds. By \cref{a:prop:generator}(b), $v^\trans J=0$ with $v>0$ forces $v\propto w$. This proves both
the characterization and \eqref{a:eq:current}.
\end{proof}
\noindent\numcheck[C1]{The residual in \eqref{a:eq:H} is $\le2\times10^{-12}$ relative on BA, RGG and IEEE 39-bus states.}

\begin{appcor}\label[appcor]{a:cor:kolmogorov}
The following are equivalent:
\begin{enumerate}[label=\textup{(\roman*)}]
\item $f$ is exact, i.e.\ $f_{ij}=\Phi_j-\Phi_i$ for some $\Phi$;
\item $\circu_Cf=0$ for every cycle $C$;
\item detailed balance holds: $w_iW_{ij}=w_jW_{ji}$ on every edge.
\end{enumerate}
In this case $w=e^\Phi/\sum_ke^{\Phi_k}$. On a tree, (i) always holds, and
$w_j/w_i=\exp\sum_{e\in\mathrm{path}(i,j)}f_e$.
\end{appcor}

\begin{proof}
(i)$\Leftrightarrow$(ii): define $\Phi$ by path sums from a root. This is well defined exactly when
all circulations vanish.

(i)$\Rightarrow$(iii): with this $\Phi$, $h\equiv0$, so \eqref{a:eq:H} holds and $w\propto e^\Phi$. Then
\eqref{a:eq:current} gives $\jmath\equiv0$.

(iii)$\Rightarrow$(i): $\jmath\equiv0$ in \eqref{a:eq:current} forces $h(\ln w)\equiv0$; see
\cite{Kolmogorov1936,Kelly1979}.
\end{proof}

\begin{apprem}
\Cref{a:thm:hodge} writes the gauge field as a potential difference plus the edge affinity $h$.
The affinity gives the stationary current through \eqref{a:eq:current}. It is that current which
has zero divergence. On a uniform twisted ring (Appendix~\ref{a:sec:twist}), $w$ is uniform, $h=f$, and the
current is nonzero when the two directional rates differ.
\end{apprem}

\begin{appthm}[A posteriori error bound]\label{a:thm:weighted}
Let $\hat\Phi\in\R^V$ be arbitrary. Put $\hat c:=c(\hat\Phi)$, $\hat h:=h(\hat\Phi)$,
$\eta:=\max_e\abs{\hat h_e}$, $s(x):=2\sinh(x/2)$, and
\[
  g_j:=\sum_{i\sim j}\hat c_{ij}\,s(\hat h_{ij})\qquad(j\in V).
\]
By the proof of \cref{a:thm:hodge}, $g^\trans=(e^{\hat\Phi})^\trans J$ is the stationarity residual of
$e^{\hat\Phi}$. Let $\hat L$ be the Laplacian with conductances $\hat c$, and define
\begin{gather*}
  \hat R:=\max_{i,j}R_{\rm eff}(i,j),\qquad Q:=\hat R\sum_e\hat c_e,\\
  \norm g_*:=(g^\trans\hat L^+g)^{1/2},\qquad \vartheta:=\tfrac12\hat R\,\norm g_1 .
\end{gather*}
The combinations $\sqrt{\hat R}\norm g_*$, $\vartheta$, $Q$ and $\hat L^+g$ are unchanged by
$\hat\Phi\mapsto\hat\Phi+\text{const}$. Also $Q\ge N-1$.
\begin{enumerate}[label=\textup{(\alph*)}]
\item \emph{A posteriori bound.} If $\vartheta<1$ and
$\mathcal B:=\sqrt{\hat R}\,\norm g_*/(1-\vartheta)<1$, then
\[
  \osc\big(\ln w-\hat\Phi\big)\le-\ln(1-\mathcal B).
\]
\item \emph{Fourth-order correction.} Under the same conditions,
\[
  \osc\big(\ln w-\hat\Phi-\hat L^+g\big)\le\mathcal B\Big(\sqrt Q\,\eta\,e^{\eta/2}+\frac{\mathcal B}{1-\mathcal B}\Big).
\]
\item \emph{Third order with an explicit constant.} Let $\hat\Phi$ satisfy the self-consistent
condition
\begin{equation}\label{a:eq:Hhat}
  \sum_{i\sim j}c_{ij}(\hat\Phi)\,\hat h_{ij}=0\quad\text{for all }j .
\end{equation}
That is, $\hat\Phi$ is the weighted least-squares potential of $f$ with conductances $c(\hat\Phi)$.
Then
\begin{gather*}
  \norm g_*\le\frac{\cosh(\eta/2)}{24}\Big(\sum_e\hat c_e\hat h_e^6\Big)^{1/2},\\
  \norm g_1\le\frac{\cosh(\eta/2)}{12}\sum_e\hat c_e\abs{\hat h_e}^3 .
\end{gather*}
Put $\varepsilon:=\cosh(\eta/2)\,\eta^3/24$. If $\varepsilon Q\le\frac14$, then
\[
  \osc\big(\ln w-\hat\Phi\big)\le\frac{\sqrt Q}{12}\cosh(\eta/2)\,\eta^3,
\]
and the bound in (b) is at most $\frac43\varepsilon Q\,(\eta e^{\eta/2}+2\varepsilon)=O(Q\eta^4)$.
\end{enumerate}
\end{appthm}

\begin{proof}
\emph{Gauge.} Under $\hat\Phi\mapsto\hat\Phi+\kappa$, $\hat c$ and $g$ scale by $e^\kappa$, while
$\hat R$ and $\hat L^+$ scale by $e^{-\kappa}$. So $\sqrt{\hat R}\norm g_*$, $\vartheta$, $Q$ and
$\hat L^+g$ are invariant.

\emph{Step 1: an exact linear equation.} Let $v:=w\,e^{-\hat\Phi}$. As in the proof of
\cref{a:thm:hodge},
\[
  w_iW_{ij}-w_jW_{ji}=\hat c_{ij}\big(v_ie^{\hat h_{ij}/2}-v_je^{-\hat h_{ij}/2}\big).
\]
Split $\hat c_{ij}e^{\pm\hat h_{ij}/2}=\tilde c_{ij}\pm b_{ij}$, with $\tilde c_{ij}:=\hat c_{ij}\cosh(\hat h_{ij}/2)$
symmetric and $b_{ij}:=\hat c_{ij}\sinh(\hat h_{ij}/2)$ antisymmetric. Then $w^\trans J=0$ reads
$(L_{\tilde c}v)_j=\sum_ib_{ij}(v_i+v_j)$ for every $j$, where $L_{\tilde c}$ is the Laplacian with conductances
$\tilde c\ge\hat c$.

Write $v=m(1+\phi)$ with $m$ the mean of $v$, so that $\sum_j\phi_j=0$. Since
$2\sum_ib_{ij}=g_j$,
\begin{equation}\label{a:eq:phi}
  L_{\tilde c}\phi=g+B\phi,\qquad(B\phi)_j:=\sum_ib_{ij}(\phi_i+\phi_j).
\end{equation}

\emph{Step 2: the drift drops out of the energy.} Because $b$ is antisymmetric,
\[
  \langle\phi,B\phi\rangle=\sum_{i,j}b_{ij}\phi_i\phi_j+\sum_j\phi_j^2\sum_ib_{ij}=\tfrac12\sum_jg_j\phi_j^2 .
\]
Pairing \eqref{a:eq:phi} with $\phi$ therefore gives
\begin{equation}\label{a:eq:energy}
  \mathcal E_{\tilde c}(\phi)=\langle g,\phi\rangle+\tfrac12\sum_jg_j\phi_j^2,
\end{equation}
where $\mathcal E_{\tilde c}(x):=\langle x,L_{\tilde c}x\rangle$. The operator $B$ is first order in $\eta$, but it
no longer appears.

\emph{Step 3: four standard facts.} Let $L_c$ be the Laplacian with conductances $c$ on a connected
graph.
\begin{enumerate}[label=(F\arabic*),nosep]
\item $(x_i-x_j)^2\le R^c_{\rm eff}(i,j)\,\mathcal E_c(x)$ for all $x$. This is the variational
definition of effective resistance.
\item For $y\perp\one$, $\langle y,x\rangle\le\mathcal E_c(x)^{1/2}\langle y,L_c^+y\rangle^{1/2}$, by
Cauchy--Schwarz applied to $L_c^{+1/2}y$ and $L_c^{1/2}x$.
\item If $c\le c'$ edgewise, then $R^{c'}_{\rm eff}\le R^c_{\rm eff}$ and $L_{c'}^+\le L_c^+$ on
$\one^\perp$ (Rayleigh monotonicity).
\item \emph{Thomson.} If $y_j=\sum_iF_{ij}$ for an antisymmetric edge flow $F$, then
$\langle y,L_c^+y\rangle\le\sum_eF_e^2/c_e$.
\end{enumerate}
For (F4), let $x:=L_c^+y$. Then
$\langle y,x\rangle=\sum_eF_{ij}(x_j-x_i)\le(\sum_eF_e^2/c_e)^{1/2}\langle y,x\rangle^{1/2}$.

\emph{Step 4: proof of (a).} Since $\tilde c\ge\hat c$, (F1) and (F3) give
$\osc(\phi)^2\le\hat R\,\mathcal E_{\tilde c}(\phi)$. Because $\phi$ has mean zero, $\abs{\phi_j}\le\osc(\phi)$,
so $\frac12\abs{\sum_jg_j\phi_j^2}\le\vartheta\,\mathcal E_{\tilde c}(\phi)$. By (F2) and (F3),
$\langle g,\phi\rangle\le\mathcal E_{\tilde c}(\phi)^{1/2}\norm g_*$. Inserting both bounds into
\eqref{a:eq:energy} gives $(1-\vartheta)\,\mathcal E_{\tilde c}(\phi)^{1/2}\le\norm g_*$, and hence
$\osc(\phi)\le\mathcal B$.

Finally, $\ln w-\hat\Phi=\ln m+\ln(1+\phi)$ and $\min\phi\le0\le\max\phi$. So
\begin{align*}
  \osc\ln(1+\phi)&=\ln\Big(1+\frac{\osc\phi}{1+\min\phi}\Big)\\
  &\le\ln\Big(1+\frac{\mathcal B}{1-\mathcal B}\Big)=-\ln(1-\mathcal B).
\end{align*}

\emph{Step 5: proof of (b).} Let $\chi:=\phi-\hat L^+g$, which has mean zero. By \eqref{a:eq:phi},
$\hat L\chi=(\hat L-L_{\tilde c})\phi+B\phi$. The right-hand side is the divergence of the antisymmetric
flow
\[
  Y_{ij}:=(\hat c_{ij}-\tilde c_{ij})(\phi_j-\phi_i)+b_{ij}(\phi_i+\phi_j).
\]
Using $\abs{\phi}\le\mathcal B$ and $\osc\phi\le\mathcal B$, together with
$\cosh(y)-1+2\sinh(y)\le2y\,e^{y}$, we get $\abs{Y_e}\le\hat c_e\,\mathcal B\,\eta\,e^{\eta/2}$.

By (F4), $\mathcal E_{\hat c}(\chi)=\langle\hat L\chi,\hat L^+\hat L\chi\rangle\le\mathcal B^2\eta^2e^\eta\sum_e\hat c_e$.
By (F1), $\osc\chi\le\mathcal B\,\eta\,e^{\eta/2}\sqrt Q$. Moreover
$\abs{\ln(1+\phi)-\phi}\le\phi^2/(2(1-\abs\phi))$, so $\osc(\ln(1+\phi)-\phi)\le\mathcal B^2/(1-\mathcal B)$.
Adding the two estimates gives (b).

\emph{Step 6: proof of (c).} Under \eqref{a:eq:Hhat}, $g_j=\sum_iF_{ij}$ with
$F_e:=\hat c_e\big(s(\hat h_e)-\hat h_e\big)$. Comparing power series,
$\abs{s(x)-x}\le\frac{\abs x^3}{24}\cosh(x/2)$. Then (F4) gives the bound on $\norm g_*$, and
$\norm g_1\le2\sum_e\abs{F_e}$ gives the bound on $\norm g_1$.

Consequently $\vartheta\le\varepsilon Q\le\frac14$ and $\sqrt{\hat R}\norm g_*\le\varepsilon\sqrt Q$.
Since $Q\ge1$, $\varepsilon\sqrt Q\le\frac14$, so $\mathcal B\le\frac43\varepsilon\sqrt Q\le\frac13$. Hence
$-\ln(1-\mathcal B)\le\frac32\mathcal B\le2\varepsilon\sqrt Q$. The bound on (b) follows the same way.

\emph{$Q\ge N-1$.} This is Foster's identity $\sum_e\hat c_eR_{\rm eff}(e)=N-1$.
\end{proof}

\noindent\numcheck[T34]{The
test set is 37 locked states: BA, RGG and square grids with $N\approx100$--$1600$, diffusive $c=1$ at
$0.25$ and $0.5\,\sigma_{\rm ref}$, where $\sigma_{\rm ref}$ is the numerical continuation limit of the frequency
spread on each graph (not a proved locking threshold), Sakaguchi $\alpha=0.3$, and an order sweep on
BA(400).}
\begin{itemize}[nosep]
\item (a), (b) and (c) hold in every case whose measured error is above the double-precision floor.
  The floor is the spread of $\ln w$ between two solvers: $3\times10^{-13}$ at $N=400$ and
  $2\times10^{-12}$ at $N=1600$.
\item (a) applies in 31 of 34 cases; $\vartheta\ge1$ occurs only on grids near a numerical continuation limit.
  Where it applies, it exceeds the true error by a factor of $1.6$--$22$ (median $3.2$). The bound of
  (b) is much looser: it exceeds the corrected error by a factor of $280$--$5\times10^4$.
\item Order sweep: $\osc(\ln w-\hat\Phi)\propto\eta^{2.97}$, and after the correction of (b),
  $\propto\eta^{3.85}$.
\item The condition of (c) holds in 27 of 34 cases; there the bound exceeds the true error by a
  factor of $62$--$3000$. Foster's identity holds to $10^{-7}$.
\end{itemize}

\begin{appprop}\label[appprop]{a:prop:globalpotential}
On a finite connected graph with positive rates in both directions on every edge,
\eqref{a:eq:Hhat} has exactly one solution up to an additive constant.
\end{appprop}
\begin{proof}
Fix $\Phi\in\one^\perp$ to remove the additive constant. Keep $S_{ij}$ fixed and replace $f$ by
$tf$, $0\le t\le1$. Write $h_{ij}^t=tf_{ij}-(\Phi_j-\Phi_i)$ and
$G_j(\Phi,t)=\sum_{i\sim j}c_{ij}(\Phi)h_{ij}^t$. Antisymmetry gives $G\in\one^\perp$.

At a root, sort the vertices by increasing potential. For each consecutive pair, split the
ordered vertices at that pair and sum $G_j=0$ over the upper set. Internal edges cancel.
Orient the crossing edges from the lower set to the upper set. The $c$-weighted average of
$\Phi_j-\Phi_i$ on these edges equals that of $tf_{ij}$. Each difference is at least the gap
between the pair, and connectedness gives positive total weight across the cut. Thus each
gap is at most $t\norm f_\infty$, and $\osc\Phi\le (N-1)t\norm f_\infty$.

At a root put $b_{ij}=c_{ij}h_{ij}^t/2$, so $b_{ij}=-b_{ji}$ and $\sum_i b_{ij}=0$
(with $b_{ij}=0$ off the edges). Differentiating the conductance and the residual gives
\[
  (D_\Phi G\,y)_j=-(L_c y)_j+\sum_i b_{ij}y_i .
\]
The second term is skew-symmetric. Hence $y^\trans D_\Phi G\,y=-y^\trans L_c y<0$
for nonzero $y\in\one^\perp$, and $D_\Phi G$ is invertible on this subspace.

The oscillation bound and the mean-zero gauge make the closed set of roots over $[0,1]$
compact. At each $t$ its roots are isolated and therefore finite. The implicit function theorem
continues each along a unique local branch. Any additional roots at parameters tending to $t$
would have a subsequence converging to one of these same roots, contradicting local uniqueness.
The number of roots is thus constant on $[0,1]$. At $t=0$ the oscillation bound forces $\Phi=0$,
which solves the equation. At $t=1$ there is consequently exactly one mean-zero solution.
\end{proof}

\begin{appprop}\label[appprop]{a:prop:treeweighted}
Let $\hat\Phi$ solve \eqref{a:eq:Hhat} exactly. Set $b=e^{\eta/2}\cosh(\eta/2)\eta^3/48$.
If $b<1$, then
\begin{equation}\label{a:eq:treeweighted}
  \osc(\ln w-\hat\Phi)\le 2(N-1)[-\ln(1-b)].
\end{equation}
At fixed $N$ the right side is $(N-1)\eta^3(1+O(\eta))/24$.
\end{appprop}

\begin{proof}
With $x_{ij}=\hat h_{ij}/2$, the rates are
$W_{ij}=e^{-\hat\Phi_i}\hat c_{ij}e^{x_{ij}}$. Replace $e^{x_{ij}}$ by
$\cosh x_{ij}+x_{ij}$ to define comparison rates $\widetilde W_{ij}$. They are positive
because $\cosh x>|x|$, and the current of $e^{\hat\Phi}$ on an oriented edge for those rates
is $2\hat c_{ij}x_{ij}$. Its divergence vanishes by \eqref{a:eq:Hhat}.
The relative rate change is bounded by
\[
  e^{|x_{ij}|}|\sinh x_{ij}-x_{ij}|
  \le e^{\eta/2}\cosh(\eta/2)\eta^3/48=b.
\]
The rest is the entrywise perturbation principle of O'Cinneide; his Thm.~1~\cite{OCinneide1993},
proved row by row rather than with spanning trees, gives the same bound with the number of states $N$
in place of $N-1$. We spell it out. For each rooted spanning tree, the ratio of its new to its old rate product lies between
$(1-b)^{N-1}$ and $(1+b)^{N-1}$. By the Markov chain tree theorem, the unnormalized stationary
weights have ratios in the same interval. Normalization changes all log ratios by a common
constant, so their oscillation is bounded by
$(N-1)\ln[(1+b)/(1-b)]\le2(N-1)[-\ln(1-b)]$.
\end{proof}

\begin{apprem}
Since $Q\ge N-1$, the a priori constant $\sqrt Q/12$ of (c) grows at least like $\sqrt N$.

Measured, $Q/N\approx3$--$4$ on BA and $30$--$150$ on RGG. On grids $Q$ also tracks the
influence contrast. Indeed $Q\le e^{\osc\hat\Phi}\,R_S\sum_eS_e$, where $R_S$ is the resistance
diameter with conductances $S$.

In the development states, $\osc(\ln w-\hat\Phi)/\eta^3$ lies in $[0.008,0.027]$. Its fitted
$N$-exponents are $-0.07$ (BA), $0.08$ (RGG) and $0.08$ (grid), against $0.49$, $0.52$ and $0.73$
for $\sqrt Q/12$.

A heuristic explanation of its typical size: by (b) the leading error is $\hat L^+g$, with
$g\approx\operatorname{div}(\hat c\hat h^3)/24$. A cubic flow on a single edge $e$ moves the potential
by $\hat c_eR_{\rm eff}(e)\abs{\hat h_e}^3/24\le\abs{\hat h_e}^3/24$, and the constants measured in the
development set are below $1/24\approx0.042$. Multiple edges contribute to the actual error, and the
atlas includes states with resolved coefficients above $1/24$, as reported in the main text.

Part (a) can certify a numerical candidate directly from the locked-state rates. The spanning-tree
bound assumes an exact solution of \eqref{a:eq:Hhat} and is independent of resistance.
\end{apprem}

\noindent\numcheck[C2]{On BA(1000) with diffusive $c=1$:}
\begin{itemize}[nosep]
\item At $\eta=2.7\times10^{-3}$, $9.1\times10^{-3}$ and $4.9\times10^{-2}$, the error
  $\frac12\osc(\ln w-\hat\Phi)$ is
  $1.9\times10^{-10}$, $7.5\times10^{-9}$ and $1.3\times10^{-6}$. The local slope is $3.0$.
\item The \emph{unweighted} least-squares potential has error $2.1\times10^{-4}$, $1.1\times10^{-3}$ and
  $1.1\times10^{-2}$ at the same points. The slope is about $1.36$, consistent with
  $O(\eta\,\osc\ln c)$, since the conductances themselves vary with the heterogeneity.
\end{itemize}

\begin{appprop}\label[appprop]{a:prop:circ}
Let $C$ be a cycle with principal differences $\abs{\Delta_e}\le d<\pi$ and winding $q_C$.
Assume $\Gamma'(x)>0$ for every $x\in[-d,d]$. Then
\begin{gather*}
\circu_Cf=4\pi\gamma\,q_C+\rho_C,\qquad \abs{\rho_C}\le\frac{M_3}{6}\sum_{e\in C}\abs{\Delta_e}^3,\\
M_3:=\sup_{\abs x\le d}\abs{F'''(x)}.
\end{gather*}
In particular, at fixed graph, untwisted states have $\circu_Cf=O(d^3)$ on every cycle.
A $q$-twisted ring in this domain has circulation $4\pi\gamma q+O(|q|^3/N^2)$.
\end{appprop}

\begin{proof}
$F$ is odd, so $F(x)=2\gamma x+r(x)$ with $\abs{r(x)}\le M_3\abs x^3/6$. Summing over $C$ gives
$\sum_e2\gamma\Delta_e=2\gamma\cdot2\pi q_C$.
\end{proof}

\begin{appthm}\label{a:thm:firstorder}
Let the state be untwisted, with lift $\tilde\theta$. Assume $\abs{\Delta_e}\le d<\pi$ and
$\Gamma'>0$ throughout $[-d,d]$. Set $r_e:=f_e-2\gamma(\tilde\theta_j-\tilde\theta_i)$
and $\rho^*:=\max_T\sum_{e\in T}\abs{r_e}$, the maximum over spanning trees $T$ of $G$. Then
\begin{align*}
  &\max_i\Big|\ln w_i-2\gamma\tilde\theta_i+\ln\textstyle\sum_ke^{2\gamma\tilde\theta_k}\Big|\\
  &\qquad\le\rho^*
  \le\frac{M_3}{6}\max_T\sum_{e\in T}\abs{\Delta_e}^3 ,
\end{align*}
and in oscillation $\osc(\ln w-2\gamma\tilde\theta)\le\rho^*$.
\end{appthm}

\begin{proof}
Write $W_{kl}=W^0_{kl}e^{r_{kl}/2}$ with $W^0_{kl}:=S_{kl}e^{\gamma(\tilde\theta_l-\tilde\theta_k)}$. The
rates $W^0$ satisfy detailed balance with respect to $w^0\propto e^{2\gamma\tilde\theta}$.

By the Markov chain tree theorem \cite{LeightonRivest1986,AnantharamTsoucas1989} (equivalently, the
cofactor formula for $w$ in \cite{Kori2009}),
$w_i\propto\sum_{T\in\mathcal T_i}\prod_{(k\to l)\in T}W_{kl}$, where $\mathcal T_i$ is the set of
spanning trees oriented towards $i$; the same holds for $w^0$ with $W^0$. Each oriented tree weight
changes by the factor $\exp(\frac12\sum_{(k\to l)\in T}r_{kl})\in[e^{-\rho^*/2},e^{\rho^*/2}]$. Hence
the unnormalized ratios $\tilde w_i/\tilde w^0_i$ lie in $[e^{-\rho^*/2},e^{\rho^*/2}]$, which directly
gives $\osc(\ln w-2\gamma\tilde\theta)=\osc\ln(\tilde w/\tilde w^0)\le\rho^*$; after
normalization $w_i/w^0_i\in[e^{-\rho^*},e^{\rho^*}]$. The last inequality follows from
\cref{a:prop:circ}.
\end{proof}

\begin{apprem}
The bound is rigorous but pessimistic for large graphs, because it sums $N-1$ edge errors with no
cancellation. For BA(1000) at $0.5\sigma_{\rm ref}$ the true error is $0.0195$ against $\rho^*=1.56$
(check C2). \Cref{a:thm:weighted} is the sharp statement.

\Cref{a:thm:weighted}(a) also applies to $\hat\Phi=2\gamma\tilde\theta$ itself, which does not satisfy
\eqref{a:eq:Hhat}. On the same BA(1000) states it bounds $\osc(\ln w-2\gamma\tilde\theta)$ by
$8.9\times10^{-4}$, $8.5\times10^{-3}$ and $0.24$ at $0.1$, $0.2$ and $0.35\,\sigma_{\rm ref}$. The sharper
oscillation tree bound $\rho^*$ gives $0.011$, $0.095$ and $0.52$, and the true values are $2.7\times10^{-4}$,
$2.2\times10^{-3}$ and $0.012$. At $0.5\sigma_{\rm ref}$, $\vartheta=2.8>1$ and (a) does not apply. The
one-step correction $2\gamma\tilde\theta+\hat L^+g$ of (b) lowers the true error by a factor of
$74$--$10^5$ over $0.1$--$0.75\,\sigma_{\rm ref}$.

For continuum phase diffusion with a global phase lift and periodic or reflecting boundary
conditions, the adjoint in \cite[App.~D]{Masuda2010} is, up to the diffusion coefficient,
$J^\dagger p=\nabla\cdot[e^{2\gamma\tilde\theta}\nabla(pe^{-2\gamma\tilde\theta})]$. This makes
$p\propto e^{2\gamma\tilde\theta}$ stationary, consistently with the Cole--Hopf formulation \cite{TonjesBlasius2009a}.
\end{apprem}

\section{Linear-response corollaries}\label{a:sec:linear}

\begin{appprop}[Reflection symmetry]\label[appprop]{a:prop:reflect}
Let $\tilde\Gamma(x):=-\Gamma(-x)$. If $\theta^*$ is a locked state of $(\omega,\Gamma)$ with frequency
$\Omega$, then $-\theta^*$ is a locked state of $(-\omega,\tilde\Gamma)$ with frequency $-\Omega$, and
it has the \emph{same} Jacobian. Hence $w$ and the whole spectrum coincide. In particular:
\begin{itemize}
\item \emph{Sakaguchi:} $\tilde\Gamma_\alpha=\Gamma_{-\alpha}$. For identical oscillators the continued
branches satisfy $\theta(-\alpha)=-\theta(\alpha)$ and $w(-\alpha)=w(\alpha)$ exactly.
\item \emph{Diffusive type:} $\tilde\Gamma_c=\Gamma_{-c}$, so $w(-c;-\omega)=w(c;\omega)$.
\end{itemize}
\end{appprop}

\begin{proof}
First, $\tilde W_{ij}=a_{ij}\tilde\Gamma'(-\Delta_{ij})=a_{ij}\Gamma'(\Delta_{ij})=W_{ij}$.
Second, 
\[
-\omega_i+\sum_ja_{ij}\tilde\Gamma(-\Delta_{ij})=-\omega_i-\sum_ja_{ij}\Gamma(\Delta_{ij})=-\Omega.
\]
For identical oscillators, $\omega$ is a frame shift. The branch statement follows from uniqueness
of the analytic continuation from $\alpha=0$ (\cref{a:thm:degree}).
\end{proof}
\noindent\numcheck[C3]{On BA(1000), $\alpha=\pm0.3$, $w$, $\theta$ and $J$ agree to the last bit.}

\begin{appthm}[Degree law]\label{a:thm:degree}
Take identical oscillators with $\Gamma(x)=\sin(x-\alpha)$ and $a_{ij}=KA_{ij}$. Let
$k_i=\sum_jA_{ij}$, let $L=\diag(k)-A$, and let $\beta:=\tan\alpha$. For $\alpha$ near zero with
$\abs\beta\le\beta_0(G)$, the
in-phase branch exists, is analytic and odd in $\beta$, and satisfies \cref{a:ass:P}. Moreover
\begin{align}
\theta(\beta)&=-\beta\,L^+k+\beta^3\theta_3+O(\beta^5),\label{a:eq:theta-exp}\\
\theta_3&=L^+\Big[\tfrac12\textstyle\sum_jA_{ij}\Delta_{1,ij}^2-\tfrac16\sum_jA_{ij}\Delta_{1,ij}^3\Big]_i,\notag\\
\ln w&=-2\beta^2\,L^+k+O(\beta^4)\label{a:eq:degree}\\
&\quad(\text{up to an additive constant}),\notag
\end{align}
where $(dx)_{ij}:=x_j-x_i$ and $\Delta_1=d(-L^+k)$. The function $\ln w$ is even in $\beta$.
\end{appthm}

\begin{proof}
Divide the locking equations by $K\cos\alpha$ and set $\mu:=(\Omega-\omega)/(K\cos\alpha)$. They
become
\begin{align*}
E_i(\theta,\mu,\beta)&:=\textstyle\sum_jA_{ij}\sin\Delta_{ij}-\beta\sum_jA_{ij}\cos\Delta_{ij}-\mu\\
&=0,
\qquad\one^\trans\theta=0 .
\end{align*}
At $\beta=0$, $(\theta,\mu)=(0,0)$, and $D_{(\theta,\mu)}E=(-L,-\one)$ is an isomorphism
$\one^\perp\times\R\to\R^V$. The implicit function theorem gives a unique analytic branch.

Since $E(-\theta,-\mu,-\beta)=-E(\theta,\mu,\beta)$, uniqueness forces the branch to be odd. Expand order by order:
\begin{itemize}[nosep]
\item Order $\beta$: $-L\theta_1-k-\mu_1\one=0$, so $\theta_1=-L^+k$.
\item Order $\beta^2$: vanishes by oddness.
\item Order $\beta^3$: expanding $\sin\Delta=\Delta-\Delta^3/6$ and $\cos\Delta=1-\Delta^2/2$ gives
  \eqref{a:eq:theta-exp}.
\end{itemize}

Since $W_{ij}=KA_{ij}\cos(\Delta_{ij}-\alpha)>0$ for small $\beta$, \cref{a:ass:P} holds. The left null
vector of a simple eigenvalue is analytic, and by \cref{a:prop:reflect} $\ln w$ is even.
For the leading term, use \cref{a:thm:hodge} with $f=F_\alpha(\Delta)=2\beta\Delta+O(\beta\Delta^3)=2\beta^2\,d\theta_1+O(\beta^4)$.
Write $\ln w=\beta^2\Phi_2+O(\beta^4)$. The order-$\beta^2$ part of \eqref{a:eq:H} reads
$\sum_iKA_{ij}(2(d\theta_1)_{ij}-(d\Phi_2)_{ij})=0$, i.e.\ $L(\Phi_2-2\theta_1)=0$. Hence
$\Phi_2=2\theta_1$ up to a constant.
\end{proof}

\noindent\numcheck[C4]{On BA(1000) and WS(1000), fitting against $\tan\alpha$:}
\begin{itemize}[nosep]
\item the error of \eqref{a:eq:degree} scales with slope $4.02$ and $4.03$;
\item the error of $\theta\approx-\beta L^+k$ scales with slope $3.01$ and $3.02$;
\item the leading spread scales with slope $2.00$.
\end{itemize}

\begin{apprem}\label[apprem]{a:rem:freqlaw}
The \emph{phase} formula $\theta\approx(\Gamma(0)/\Gamma'(0))L^+k$ is in \cite[Eq.~4]{Skardal2015}.
What is new here is the influence law \eqref{a:eq:degree} and its exact parity. For heterogeneous
frequencies with $\Gamma(0)=0$ and $\Gamma'(0)>0$, the same argument (IFT in the frequency
spread $\sigma$) gives $\theta=L^+\omega/(K\Gamma'(0))+O(\sigma^2)$ and
$\ln w=2\gamma L^+\omega/(K\Gamma'(0))+O(\sigma^2)$.
\end{apprem}

\section{Exact skin length and the soft sink}\label{a:sec:chains}

Consider a chain $1-2-\cdots-N$ with uniform weight $a$, identical oscillators, and
$\Delta_i:=\theta_{i+1}-\theta_i$. By \cref{a:cor:kolmogorov}, every positive locked state of a tree
satisfies
\begin{equation}\label{a:eq:treeratio}
  \frac{w_{i+1}}{w_i}=\rho(\Delta_i):=\frac{\Gamma'(\Delta_i)}{\Gamma'(-\Delta_i)}\quad\text{exactly.}
\end{equation}

\begin{appprop}\label[appprop]{a:prop:skin}
Let $x_0$ be a zero of $\Gamma$ with $\Gamma'(x_0)>0$ and $\Gamma'(-x_0)>0$. On the semi-infinite
chain $\{1,2,\dots\}$, the profile $\theta_i=-(i-1)x_0$ is an exact locked state with
$\Omega=\omega+a\Gamma(-x_0)$. Its generator is a uniform Hatano--Nelson chain
\cite{HatanoNelson1996} with rates $a\Gamma'(-x_0)$ (rightwards) and $a\Gamma'(x_0)$ (leftwards).
If $r:=\Gamma'(-x_0)/\Gamma'(x_0)<1$, the chain is positive recurrent and
$w_i=(1-r)r^{i-1}$, so the skin length is $\xi=-1/\ln r$.

For Sakaguchi coupling, $x_0=\alpha$, $r=\cos2\alpha$ and $\xi=-1/\ln\cos2\alpha$, valid for
$0<\abs\alpha<\pi/4$.
\end{appprop}

\begin{proof}
Node $1$: $\omega+a\Gamma(-x_0)$. Node $i\ge2$: $\omega+a\Gamma(x_0)+a\Gamma(-x_0)=\omega+a\Gamma(-x_0)$.
The rates follow from \eqref{a:eq:J}, and $w$ follows from \eqref{a:eq:treeratio}.
\end{proof}

\begin{apprem}
The free-end lag $-x_0$ is the free-end boundary condition of the continuum limit of
\cite[(2.5)]{KopellErmentrout1986}, whose $H^-$ is $a\Gamma$ here. In finite chains, it strictly bounds every lag
\cite[Thms.~2.3--2.4]{RenErmentrout1998}. That it holds exactly on the semi-infinite chain, the exact
influence ratio and $\xi$ are new. On the in-phase branch with diffusive coupling and
$\Gamma'(0)>0$, every edge has zero lag and $r=1$. The uniform stationary measure
normalizes to $w_i=1/N$ on a finite chain and has infinite mass on the semi-infinite one.
\end{apprem}

Existence, monotone lags, antisymmetry and stability of the locked state of a finite chain of identical
oscillators are known for a broad class of couplings \cite[Thms.~2.3--2.6]{RenErmentrout1998}. For
Sakaguchi coupling the next lemma makes them explicit, through the center orbit, with the bounds that
the soft-sink theorem needs.

\begin{applem}[Two-source state]\label[applem]{a:lem:twosource}
Take Sakaguchi coupling with $0<\alpha<\pi/4$, uniform weight $a>0$, identical oscillators and even
$N=2n$. Set
\begin{gather*}
  r:=\cos2\alpha,\qquad \mu:=\frac{\Omega-\omega}{a},\qquad \mu_0:=-\sin2\alpha,\\
  \beta:=\arcsin(\sin2\alpha-\sin\alpha),\qquad q:=\frac{\cos\alpha}{\cos\beta},\\
  \kappa_\alpha:=\tan2\alpha+\tan\alpha .
\end{gather*}
Work with shifted lags $u:=\Delta+\alpha$, and for $\delta\in\R$ define
\begin{gather*}
  T_\delta(u):=\arcsin\bigl(\sin2\alpha-\sin(2\alpha-u)-\delta\bigr),\\
  T:=T_0,\qquad
  \rho(u):=\frac{\cos(2\alpha-u)}{\cos u}.
\end{gather*}
So $\rho(u)$ is the ratio $\rho(\Delta)$ of \eqref{a:eq:treeratio}, and
$T_\delta(u)=B_\mu(u-\alpha)+\alpha$ with $\mu=\mu_0+\delta$, where
$B_\mu(\Delta):=-\alpha-\arcsin(\mu-\sin(\Delta-\alpha))$ is the backward lag recursion. Let
$u^*_0:=\alpha$, $u^*_{m+1}:=T(u^*_m)$ be the \emph{center orbit}, and $P:=\prod_{j\ge1}r/\rho(u^*_j)$.
\begin{enumerate}[label=\textup{(\alph*)}]
\item \emph{Center orbit.} On $(0,\alpha]$, $T$ is increasing, $\rho$ is increasing with $r\le\rho\le1$, and
\begin{equation}\label{a:eq:Tbounds}
  r\,u \le T(u) < u\,\rho(u) \le u .
\end{equation}
Hence $u^*_m\downarrow0$, $0<P\le1$, and $\alpha r^m\le u^*_m\le\alpha r^{m-1}/P$ for $m\ge1$.
Equivalently, $B^m_{\mu_0}(0)\to-\alpha$; i.e., the center lag $0$ lies in the basin of the fixed point
$-\alpha$ for every $\alpha\in(0,\pi/4)$.
\item \emph{Existence and uniqueness.} There is exactly one positive (\cref{a:ass:P}) locked state with
$\theta_i=\theta_{N+1-i}$. It has $\Delta_n=0$, $\Delta_{N-k}=-\Delta_k$ and $\Delta_{n-m}=u_m-\alpha$ for
$0\le m<n$. Here $u_0=\alpha$, $u_{m+1}=T_{\delta_N}(u_m)$ and $\mu=\mu_0+\delta_N$, where
$\delta_N\in(0,u^*_n]$ is the unique root of $\Phi_n(\delta):=T_\delta^n(\alpha)$ along which the orbit
stays positive. The state is linearly stable.
\item \emph{Frequency.} $(1-q)\cos\beta \, u^*_n\le\delta_N\le u^*_n$. Hence $\mu-\mu_0\asymp r^{N/2}$.
\item \emph{Lags and influence ratios.} $0<u_m\le u^*_m$ for $0\le m<n$. Hence, for every left-half edge
$k=n-m$ with $1\le m<n$,
\begin{gather*}
  0<\Delta_k+\alpha\le\frac{\alpha}{P}\,r^{\,n-k-1},\\
  r\le\frac{w_{k+1}}{w_k}=\rho(u_m)\le\min\bigl\{1,\ r\,e^{\kappa_\alpha u_m}\bigr\}.
\end{gather*}
In particular $w$ is nonincreasing on $\{1,\dots,n\}$, and $w_{k+1}/w_k=r\,(1+O(r^{N/2-k}))$.
\end{enumerate}
\end{applem}

\begin{proof}
For $1\le k\le n$, node $k$ reads $\mu=\sin(\Delta_k-\alpha)-\sin(\Delta_{k-1}+\alpha)$.
Node $1$ is included by a ghost edge with $\Delta_0:=-\alpha$, so that $\sin(\Delta_0+\alpha)=0$.
Symmetry gives $\Delta_{N-k}=-\Delta_k$, hence $\Delta_n=0$, and the equation of node $N+1-k$ is the
same as that of node $k$.

Under \cref{a:ass:P}, $\cos(\Delta_{k-1}+\alpha)>0$. So node $k$'s equation is equivalent to
$\Delta_{k-1}+\alpha=\arcsin(\sin(\Delta_k-\alpha)-\mu)$ on the principal branch. Put
$u_m:=\Delta_{n-m}+\alpha$ for $0\le m\le n$, so $u_n=\Delta_0+\alpha$. Then the node equations read
$u_{m+1}=T_\delta(u_m)$, with $u_0=\alpha$ and the end condition $u_n=0$.

On edge $k$, positivity of both rates means $\cos(u-2\alpha)>0$ and $\cos u>0$ at $u=u_{n-k}$. On that
range
\begin{equation}\label{a:eq:Tderiv}
\begin{gathered}
  \partial_uT_\delta(u)=\frac{\cos(2\alpha-u)}{\cos T_\delta(u)}>0,\\
  \partial_\delta T_\delta(u)=-\frac1{\cos T_\delta(u)}\le-1 .
\end{gathered}
\end{equation}

(a) Let $u\in(0,\alpha]$. Then $2\alpha-u\in[\alpha,2\alpha)\subset(0,\pi/2)$, so the argument of the
arcsine lies in $[0,\sin2\alpha-\sin\alpha]\subset[0,1)$ and $T(u)\ge0$. Also
\begin{align*}
  \sin T(u)-\sin u&=\sin2\alpha-2\sin\alpha\cos(\alpha-u)\\
  &=2\sin\alpha\,\bigl(\cos\alpha-\cos(\alpha-u)\bigr)<0,
\end{align*}
so $0\le T(u)<u$.

By \eqref{a:eq:Tderiv}, $T'(u)\ge\cos(2\alpha-u)\ge\cos2\alpha=r$. Since $0\le T(u)<u<\pi/2$, also
$T'(u)\le\cos(2\alpha-u)/\cos u=\rho(u)$.

Next, $(\ln\rho)'(u)=\tan(2\alpha-u)+\tan u\in(0,\kappa_\alpha]$. So $\rho$ increases from $\rho(0)=r$
to $\rho(\alpha)=1$, and $\ln(\rho(u)/r)\le\kappa_\alpha u$. Integrating $r\le T'\le\rho$ from
$T(0)=0$, and using that $\rho$ is strictly increasing for the upper bound, gives \eqref{a:eq:Tbounds};
the last inequality is an equality only at $u=\alpha$, where $\rho=1$.

The center orbit is therefore positive and strictly decreasing. Its limit is a fixed point of $T$
in $[0,\alpha)$, which must be $0$.
Induction on \eqref{a:eq:Tbounds} gives $u^*_m\ge\alpha r^m$. In the other direction,
\[
  u^*_m\le\alpha\prod_{j=0}^{m-1}\rho(u^*_j)=\alpha r^{m-1}\prod_{j=1}^{m-1}\frac{\rho(u^*_j)}r\le\frac{\alpha r^{m-1}}P,
\]
because $\rho(u^*_0)=\rho(\alpha)=1$ and every factor $\rho/r$ is $\ge 1$.

For $j\ge1$, $u^*_{j+1}\le\rho(u^*_1)\,u^*_j$ with $\rho(u^*_1)<1$. So $\sum_ju^*_j<\infty$ and
$\ln(1/P)\le\kappa_\alpha\sum_{j\ge1}u^*_j<\infty$. Finally, this yields the identity $B_{\mu_0}^m(0)=u^*_m-\alpha$.

\smallskip
(b) First let $\delta\le0$. For $u\in(0,\pi/2)$ we have $\sin(2\alpha-u)<\sin2\alpha$, so the argument of
the arcsine exceeds $-\delta\ge0$ and $T_\delta(u)>0$ wherever it is defined. By induction $u_n>0$,
so there is no solution.

Now let $\delta>0$. If $u_m\le0$ for some $m<n$, then
$2\alpha-u_m\in[2\alpha,\pi/2)$ by the positive-rate window. The argument of the next
arcsine is at most $-\delta<0$, so whenever defined the next iterate is negative.
Induction contradicts $u_n=0$. Thus $u_m>0$ for $m<n$. Since $T_\delta<T$ and $T$ is increasing,
$u_m\le u^*_m$. The last equation implies $\delta\le\sin u^*_n\le u^*_n$.

\emph{Existence and uniqueness of $\delta_N$.} For $n=1$ the boundary equation gives
$\delta_N=\sin2\alpha-\sin\alpha$. For $n\ge2$, vary $\delta$ in $[0,u^*_n]$ while
$u_m(\delta)>0$ for all $m<n$. In this interval every such input lies in $(0,\alpha]$.
The next arcsine argument lies between $-u^*_n$ and $\sin2\alpha-\sin\alpha$,
strictly inside $(-1,1)$. Consequently each intermediate iterate is continuous
and strictly decreasing in $\delta$ until one reaches zero. The chain rule gives
\begin{equation}\label{a:eq:dPhi}
\begin{aligned}
  \Phi_n'(\delta)&=-\sum_{k=1}^{n}\frac1{\cos u_k}\prod_{j=k}^{n-1}\frac{\cos(2\alpha-u_j)}{\cos u_{j+1}}\\
  &\;\le\;-\frac1{\cos u_n}\;\le\;-1,
\end{aligned}
\end{equation}
since every term is negative. Also $\Phi_n(0)=u^*_n>0$.

If the intermediate iterates stay positive through $\delta=u^*_n$, then
$\Phi_n(u^*_n)\le0$. Otherwise the last intermediate iterate reaches zero first, since
$T_\delta(u)\le T(u)<u$ for $0<u\le\alpha$. At this first parameter $\bar\delta>0$,
$\Phi_n$ has limit $-\arcsin\bar\delta<0$. In both cases continuity gives a root
$\delta_N\in(0,u^*_n]$. The derivative bound and the necessary orbit conditions above give uniqueness.

On the left half every lag has $u\in(0,\alpha]$, so both rates of every edge are positive. Hence the
state satisfies \cref{a:ass:P} and is stable by \cref{a:prop:generator}.

(c) The upper bound was shown in (b). For $\delta\in[0,\delta_N]$ the orbit satisfies
$u_{j+1}\le T(\alpha)=\beta<\alpha$ for $j\ge0$. So in \eqref{a:eq:dPhi}, each factor
$\cos(2\alpha-u_j)/\cos u_{j+1}$ is at most $\cos\alpha/\cos\beta=q<1$, and $1/\cos u_k\le\sec\beta$.
Hence $\abs{\Phi_n'}\le\sec\beta/(1-q)$. Integrating from $0$ to $\delta_N$ gives
$u^*_n=\Phi_n(0)-\Phi_n(\delta_N)\le\delta_N\sec\beta/(1-q)$. Finally $u^*_n\asymp r^n$ by (a).

(d) The bound $u_m\le u^*_m$ is from (b), and (a) bounds $u^*_m$; note $m-1=n-k-1$. The ratio identity
is \eqref{a:eq:treeratio}, since $\cos(\Delta_k-\alpha)/\cos(\Delta_k+\alpha)=\rho(u_m)$. The bounds on
$\rho$ are from (a).
\end{proof}
\noindent\numcheck[C5, L1--L5]{Computed with 160 digits for
$\alpha\in\{0.05,\dots,0.78\}$ and $n\le60$ ($n\le150$ for $\alpha=0.2,0.3$). The bounds
\eqref{a:eq:Tbounds}, $\alpha r^m\le u^*_m\le\alpha r^{m-1}/P$, (c) and $1\le\abs{\Phi_n'}\le\sec\beta/(1-q)$
hold in all 110 cases. Every one of the $N$ node equations holds to $10^{-120}$ (L5). An independent
double-precision Newton continuation of the full chain agrees to $\le6\times10^{-11}$ in the lags in
15 of 16 runs. The 16th ($\alpha=0.7$, $N=40$) does not converge, because $\abs{\lambda_2}=5\times10^{-15}$
makes the Jacobian singular in double precision. The fitted $d\ln\abs{\mu-\mu_0}/dN$ is $-0.178$
against $\frac12\ln r=-0.181$ ($\alpha=0.4$).}

\begin{apprem}[Sharpness]
The limit $c_*:=\lim_mr^{-m}u^*_m$ exists, because
$r^{-m}u^*_m=\alpha\prod_{j<m}T(u^*_j)/(ru^*_j)$ and each factor lies in
$[1,e^{\kappa_\alpha u^*_j}]$ by \eqref{a:eq:Tbounds}. In \eqref{a:eq:dPhi} the last steps of the orbit
sit near $u=0$, where every factor is close to $r$; this suggests $\Phi_n'\to-1/(1-r)$. Hence
$\delta_N/u^*_n\to1-r$ and $\mu-\mu_0\sim(1-r)\,c_*\,r^{N/2}$. This limit is only sketched and is not used below.

\numcheck[L3]{The relative deviation of
$\delta_N/u^*_n$ from $1-r$ decays geometrically in $n$, and is $\le3\times10^{-4}$ in every case with
$r^n\le10^{-5}$.}
\end{apprem}

\begin{appthm}[Soft sink]\label{a:thm:sink}
In the setting of \cref{a:lem:twosource}, let $c_k:=w_kW_{k,k+1}=w_{k+1}W_{k+1,k}$ for $1\le k<N$ be the
edge conductances (the two expressions agree by detailed balance). Let $R_N:=\sum_{k=1}^{N-1}c_k^{-1}$
be the end-to-end resistance and $h_N(i):=R_N^{-1}\sum_{j<i}c_j^{-1}$ the equilibrium potential. Set
\[
  V_N:=\sum_{i=1}^Nw_i\bigl(1-2h_N(i)\bigr)^2,\quad
  \varepsilon_N:=\frac{2n(n-1)\,r^{\,n-1}}{P}.
\]
The eigenvalues $0=\lambda_1>\lambda_2\ge\cdots$ of $J$ are real, because a chain is
reversible. It holds:
\begin{enumerate}[label=\textup{(\roman*)}]
\item \emph{Two-sided bound, every $N$.} The slow eigenvector is antisymmetric, $v_{N+1-i}=-v_i$, and
\begin{equation}\label{a:eq:sinkbounds}
\begin{gathered}
  \frac4{R_N}\;\le\;\abs{\lambda_2}\;\le\;\frac4{R_N\,V_N},\\
  0\le1-V_N\le\frac{4n(n-1)}{a\cos\alpha\;R_N}\le\varepsilon_N .
\end{gathered}
\end{equation}
In particular $0\le\abs{\lambda_2}R_N/4-1\le\varepsilon_N/(1-\varepsilon_N)=O(N^2r^{N/2})$ once
$\varepsilon_N<1$.
\item \emph{Time-scale separation.} For $n\ge2$,
\begin{equation}\label{a:eq:thirdrate}
  \abs{\lambda_3}\ge aP\cos\alpha\,(1-\sqrt r)^2>0.
\end{equation}
For $\alpha$ fixed, this bound is uniform in $N$; in addition, $\abs{\lambda_3}/\abs{\lambda_2}\to\infty$.
\item \emph{Resistance asymptotics.} $a\,r^{n}R_N\to A(\alpha)$, where
\begin{equation}\label{a:eq:Aalpha}
\begin{aligned}
A(\alpha)&=\frac{2rP}{1-r}\Big[\sec\alpha\\
&\qquad+2\sum_{m\ge1}\Big(\prod_{j=1}^m\rho_j\Big)\sec(\Delta^*_m-\alpha)\Big],\\
&\qquad \rho_j:=\rho(u^*_j),\quad \Delta^*_m:=u^*_m-\alpha .
\end{aligned}
\end{equation}
\end{enumerate}
Consequently $\lambda_2=-a\,C(\alpha)\,r^{N/2}\,(1+o(1))$ with $C(\alpha)=4/A(\alpha)$.
\end{appthm}

\begin{proof}
By \cref{a:lem:twosource}(b), $\Delta_{N-k}=-\Delta_k$. This means that
$W_{N+1-k,N-k}=a\cos(\Delta_k-\alpha)=W_{k,k+1}$, so the reflection $\sigma(i)=N+1-i$ commutes with $J$.
It also gives $w\circ\sigma=w$, $c_{N-k}=c_k$ and $h_N\circ\sigma=1-h_N$.

Now notice that $J$ is self-adjoint on $\ell^2(w)$, with Dirichlet form
$\langle-Ju,u\rangle_w=\mathcal E(u):=\sum_kc_k(u_{k+1}-u_k)^2$. Its spectrum is the union of its
spectra on odd and on even functions. Write $v:=u|_{\{1,\dots,n\}}$. Then:
\begin{itemize}[nosep]
\item $u$ odd: $\mathcal E(u)=2\mathcal E_o(v)$ and $\norm u_w^2=2\sum_{i\le n}w_iv_i^2$, where
$\mathcal E_o(v):=\sum_{k<n}c_k(v_{k+1}-v_k)^2+2c_nv_n^2$.
\item $u$ even: same with $\mathcal E_e(v):=\sum_{k<n}c_k(v_{k+1}-v_k)^2$.
\end{itemize}
In the variable $\diag(w)^{1/2}v$ these forms are represented by symmetric tridiagonal matrices with
$M_o=M_e+(2c_n/w_n)\,e_ne_n^\trans$.

This is a positive rank-one update. When $n=1$ the even sector consists only of the neutral
mode. When $n\ge2$, every eigenvector of the unreduced tridiagonal $M_e$ has nonzero
last component, and its eigenvalues are simple. So the eigenvalues interlace strictly:
$0=\lambda_1(M_e)<\lambda_1(M_o)<\lambda_2(M_e)<\lambda_2(M_o)<\cdots$. In both cases
$\abs{\lambda_2}=\lambda_1(M_o)=:\lambda_{\rm odd}$ and the slow mode is odd.
For $n\ge2$, $\abs{\lambda_3}=\lambda_2(M_e)=:\lambda_{\rm even}$.

(i) \emph{Lower bound.} $\mathcal E_o$ is the Dirichlet form of the half chain on $\{1,\dots,n\}$
with measure $w$, killed from node $n$ through conductance $c'_n:=2c_n$. Put $c'_k:=c_k$ for $k<n$, and
let $L_o$ be the generator of this half chain. Let $t(i)$ be the expected killing time from $i$,
so $-L_ot=\one$. Let $\phi>0$ be the Perron eigenvector, $-L_o\phi=\lambda_{\rm odd}\phi$. Then
\begin{align*}
  \textstyle\sum_iw_i\phi_i&=\langle\phi,-L_ot\rangle_w=\lambda_{\rm odd}\langle\phi,t\rangle_w\\
  &\le\lambda_{\rm odd}\max_it(i)\,\sum_iw_i\phi_i ,
\end{align*}
so $\lambda_{\rm odd}\ge1/\max_it(i)$.

Let $s_k$ be the expected time to cross edge $k$ starting from node $k$. These satisfy
$c'_ks_k=w_k+c_{k-1}s_{k-1}$ with $c_0:=0$, hence $c'_ks_k=w(\{1,\dots,k\})$. So $t$ is maximal at
node $1$, and
\begin{align*}
  t(1)&=\sum_{k<n}\frac{w(\{1,\dots,k\})}{c_k}+\frac{w(\{1,\dots,n\})}{2c_n}\\
  &\le\frac12\Big(\sum_{k<n}\frac1{c_k}+\frac1{2c_n}\Big)=\frac{R_N}4 .
\end{align*}
This uses $w(\{1,\dots,n\})=\tfrac12$ and $R_N=2\sum_{k<n}c_k^{-1}+c_n^{-1}$.

\emph{Upper bound.} The function $1-2h_N$ is odd, with $\mathcal E(1-2h_N)=4\,\mathcal E(h_N)=4/R_N$
and $\norm{1-2h_N}_w^2=V_N$. The Rayleigh quotient gives the bound.

\emph{The bound on $V_N$.}
\begin{align*}
  1-V_N&=4\sum_iw_ih_N(1-h_N)\le8\sum_{i\le n}w_ih_N(i)\\
  &=\frac8{R_N}\sum_{j<n}\frac{w(\{j+1,\dots,n\})}{c_j}.
\end{align*}
By \cref{a:lem:twosource}(d), $w$ is nonincreasing on $\{1,\dots,n\}$. Also
$c_j=w_{j+1}\,a\cos(\Delta_j+\alpha)\ge a\cos\alpha\,w_{j+1}$. So each summand is at most
$(n-j)/(a\cos\alpha)$, which gives the middle inequality of \eqref{a:eq:sinkbounds}.

For the last inequality, $R_N\ge c_n^{-1}=(a\cos\alpha\,w_n)^{-1}$. By \cref{a:lem:twosource}(d) and the
monotonicity of $\rho$,
$w_n=w_1\prod_{m=1}^{n-1}\rho(u_m)\le\tfrac12\prod_{m=1}^{n-1}\rho(u^*_m)\le r^{n-1}/(2P)$.

\smallskip
(ii) Let $m_i=w_1r^{i-1}$ on the half chain. By \cref{a:lem:twosource}(d),
$w_i/m_i$ is a product of factors $\rho(u_j)/r$ lying between $1$ and $\rho(u_j^*)/r$,
so $m_i\le w_i\le m_i/P$. Also $c_i\ge a\cos\alpha\,w_{i+1}\ge a\cos\alpha\,m_{i+1}$
for $i<n$. For the weighted variance
$\operatorname{Var}_w(v):=\min_{z\in\R}\sum_{i\le n}w_i(v_i-z)^2$, the upper bound on $w_i$ gives
$\operatorname{Var}_w(v)\le P^{-1}\operatorname{Var}_m(v)$. Thus for nonconstant $v$
\[
  \frac{\mathcal E_e(v)}{\operatorname{Var}_w(v)}
  \ge aP\cos\alpha\,
  \frac{\sum_{i<n}m_{i+1}(v_{i+1}-v_i)^2}{\operatorname{Var}_m(v)}.
\]
The right-hand quotient, before the prefactor, is that of a reflecting birth--death chain with
right and left rates $r$ and $1$. In coordinates $x_i=\sqrt{m_i}\,v_i$ its energy is
$\norm{Bx}^2$, with $(Bx)_i=x_{i+1}-\sqrt r\,x_i$ for $i<n$. The positive eigenvalues of
$B^\trans B$ are the eigenvalues of $BB^\trans$, which has diagonal $1+r$ and adjacent entries
$-\sqrt r$. Its smallest eigenvalue is at least $(1-\sqrt r)^2$ by Gershgorin.
Taking the infimum over nonconstant $v$ proves \eqref{a:eq:thirdrate}. Part (i) gives
$\abs{\lambda_2}\to0$ since $R_N\ge1/(a\cos\alpha\,w_n)$, $w_n\le r^{n-1}/(2P)$ and $V_N\to1$.

(iii) Index the left half from the center: edge $n-m$ carries $u_m=u_m^{(N)}$. By
\eqref{a:eq:treeratio}, $w_{n-m}=w_1\prod_{j=m+1}^{n-1}\rho(u_j)$, and
$c_{n-m}=a\,w_{n-m}\cos(2\alpha-u_m)$. Hence
\begin{gather*}
  \frac{r^n}{c_{n-m}}=\frac{r^{m+1}}{a\,w_1\cos(2\alpha-u_m)\prod_{j=m+1}^{n-1}(\rho(u_j)/r)}\\(0\le m<n).
\end{gather*}
We pass to the limit $N\to\infty$ in four steps.
\begin{enumerate}[label=(\arabic*),nosep]
\item For fixed $j$, $u^{(N)}_j\to u^*_j$. This holds because $\delta_N\to0$ by
\cref{a:lem:twosource}(c) and $T_\delta$ is continuous in $\delta$.
\item The logarithms $\ln(\rho(u_j)/r)$ lie in $[0,\kappa_\alpha u^*_j]$, which is summable. So $\prod_{j=m+1}^{n-1}(\rho(u_j)/r)\to P^{-1}\prod_{j\le m}(r/\rho_j)$ by
dominated convergence.
\item $\tfrac12=\sum_{k\le n}w_k=w_1\sum_{k\le n}\prod_{i<k}\rho(u_{n-i})$. Each term tends to $r^{k-1}$,
since $u_{n-i}\le u^*_{n-i}\to0$, and each is dominated by $r^{k-1}/P$. Hence $w_1\to(1-r)/2$, and
$w_1\ge(1-r)P/2$ for all $N$.
\item By (3), $\cos(2\alpha-u_m)\ge r$ and $\prod(\rho/r)\ge1$, we get
$r^n/c_{n-m}\le2r^m/(a(1-r)P)$, summable in $m$.
\end{enumerate}
Then $r^nR_N=r^n/c_n+2\sum_{m=1}^{n-1}r^n/c_{n-m}$ with dominated convergence gives
\eqref{a:eq:Aalpha}. The term $m=0$ has $u_0=\alpha$, so $\Delta^*_0=0$ and the factor is $\sec\alpha$.

The consequence follows from (i) and (iii).
\end{proof}
\noindent\numcheck[S1--S5]{The checks use multiprecision
eigenvalues (Sturm bisection, 160 digits) of the full chain and of both sectors, for
$\alpha\in\{0.05,\dots,0.78\}$ and $N\le120$ ($N\le300$ for $\alpha=0.2,0.3$). They confirm
\eqref{a:eq:sinkbounds}, odd slow mode, $\abs{\lambda_3}=\lambda_{\rm even}$, and the sector
decomposition of the full spectrum to $10^{-150}$ (Table~\ref{a:tab:sink}).}

\begin{table*}[t]
\caption{\label{a:tab:sink}Soft-sink constant $C(\alpha)$ against 160-digit eigenvalues of the exact two-source state
(checks C6 and S5).
The step-profile estimate $(1-r)^2$, which ignores the central kink (and equals half the $\alpha\to0$ weak-skin
limit $2(1-r)^2$), is too small by a factor of $2.2$--$6.1$ over this range.}
\begin{ruledtabular}
\begin{tabular}{cccccc}
$\alpha$ & $C(\alpha)$ from \eqref{a:eq:Aalpha} & $\abs{\lambda_2}r^{-N/2}$ (numerics) & $N$ & $\abs{\lambda_2}R_N/4-1$ & $P$\\
\colrule
0.2 & 0.0134524 & 0.0134523 & 300 & $2.0\times10^{-4}$ & 0.2522\\
0.3 & 0.0731730 & 0.0731730 & 300 & $3.7\times10^{-11}$ & 0.2541\\
0.4 & 0.2616139 & 0.2616139 & 120 & $3.7\times10^{-8}$ & 0.2546\\
0.5 & 0.7909598 & 0.7909598 & 120 & $1.8\times10^{-14}$ & 0.2495\\
0.6 & 2.4633657 & 2.4633657 & 120 & $<10^{-16}$ & 0.2287\\
\end{tabular}
\end{ruledtabular}
\end{table*}

\begin{apprem}[Relation to metastability and sharpness of the bounds]
\Cref{a:thm:sink}(i) has the capacity interpretation of reversible metastability
\cite[Thm.~1.3]{BEGK2002} (see also Ref.~\cite[Ch.~7]{BovierdenHollander2015}). The capacity between the
endpoint nodes is $1/R_N$ and each half has mass $1/2$. Adding the two exit rates gives
$\abs{\lambda_2}\sim 2(1/R_N)/(1/2)=4/R_N$ at fixed $\alpha$.
The explicit $\varepsilon_N$ can be loose at
small $\alpha$. For example, at $\alpha=0.05$ the slow mode is odd for every even $N$, while
$\varepsilon_N<1$ needs $N$ in the thousands.
\end{apprem}

\begin{apprem}\label[apprem]{a:rem:subband}
Set $a=1$ and choose $\alpha$ so that $r=1/100$. In the symmetric representation of $-J$ take
the two endpoint coordinate vectors and the indicator of the four central sites. For large even
$N$ the supports are disjoint and nonadjacent, so min--max bounds $\abs{\lambda_3}$ by the
largest of their three Rayleigh quotients. Each endpoint quotient tends to $r$. At fixed distance
from the center the lags converge to the center orbit by \cref{a:lem:twosource}, as in (iii) above.
Write $b_1=\cos(2\alpha-u^*_1)$ and $d_j=\cos u^*_j$. In the central quotient the middle
edge cancels, and the two edges adjacent to it each contribute $b_1+d_1-2\sqrt{b_1d_1}$.
The two outgoing edges from the support each contribute $d_2$. The limiting quotient is therefore
$(b_1+d_1+d_2-2\sqrt{b_1d_1})/2$.
Here $\sin u^*_1=\sqrt{1-r^2}-\sqrt{(1-r)/2}$ lies between $0.28$ and $0.30$.
Consequently $d_1>0.9$ and $1/4<b_1<1/2$, using
$b_1=r d_1+\sqrt{1-r^2}\sin u^*_1$. With $d_1,d_2\le1$ the quotient is less than
$1.25-\sqrt{0.9/4}<0.78$. Thus the limit superior is below the uniform chain's band edge
$(1-\sqrt r)^2=0.81$.
\end{apprem}

\begin{appprop}\label[appprop]{a:prop:continuum}
In the setting of \cref{a:lem:twosource}, set $h:=\alpha^2$, $D_m:=(T^m)'(\alpha)$ (the derivative of the
$m$-th iterate of $T$ at the center value $u^*_0=\alpha$), and
\[
  \kappa'(\alpha):=\lim_{m\to\infty}r^{-m}D_m,\qquad S(\alpha):=\sum_{m\ge0}D_m .
\]
\begin{enumerate}[label=\textup{(\alph*)}]
\item \emph{Exact representation, every $\alpha\in(0,\pi/4)$.} Both quantities are finite and
positive, and
\begin{gather*}
  P=\frac{\cos\alpha}{r\,\kappa'(\alpha)},\qquad
  A(\alpha)=\frac{2\,\bigl(2S(\alpha)-1\bigr)}{(1-r)\,\kappa'(\alpha)},\\
  C(\alpha)=\frac{2(1-r)\,\kappa'(\alpha)}{2S(\alpha)-1}.
\end{gather*}
\item \emph{Continuum limit.} There are $K,\alpha_0>0$ such that for $0<\alpha\le\alpha_0$ and all $m\ge0$
\begin{equation}\label{a:eq:kink}
  \bigl|\Delta^*_m+\alpha\tanh(\alpha^2m)\bigr|\le K\alpha^3e^{-\alpha^2m},
\end{equation}
and, with all $O$-constants depending only on $\alpha_0$,
\begin{gather*}
  \kappa'(\alpha)=4\bigl(1+O(\alpha^2)\bigr),\qquad \alpha^2S(\alpha)=1+O(\alpha^2),\\
  c_*=2\alpha\bigl(1+O(\alpha^2)\bigr).
\end{gather*}
Consequently
\begin{align*}
  P&=\tfrac14+O(\alpha^2),\\
  C(\alpha)&=8\alpha^4\bigl(1+O(\alpha^2)\bigr)=8\sin^4\!\alpha\,\bigl(1+O(\alpha^2)\bigr)\\
  &=2(1-r)^2\bigl(1+O(\alpha^2)\bigr).
\end{align*}
\end{enumerate}
\end{appprop}

\begin{proof}
(a) By \eqref{a:eq:Tderiv} with $\delta=0$, $T'(u)=\cos(2\alpha-u)/\cos T(u)$. Hence
\[
  \rho(u^*_j)=T'(u^*_j)\,\frac{\cos u^*_{j+1}}{\cos u^*_j},\qquad
  T'(\alpha)=\frac{\cos\alpha}{\cos u^*_1}.
\]
The product telescopes:
\begin{equation}\label{a:eq:telescope}
  \prod_{j=1}^m\rho(u^*_j)=\frac{D_{m+1}}{T'(\alpha)}\cdot\frac{\cos u^*_{m+1}}{\cos u^*_1}
  =\frac{D_{m+1}\cos u^*_{m+1}}{\cos\alpha}.
\end{equation}

\emph{Finiteness.} By the proof of \cref{a:lem:twosource}(a), $r\le T'\le\rho$ on $(0,\alpha]$ and
$\ln(\rho(u)/r)\le\kappa_\alpha u$. So $r^{-m}D_m=\prod_{j<m}T'(u^*_j)/r$ increases to a finite
limit $\kappa'(\alpha)\ge1$. Moreover $D_m\le\kappa'r^m$, so $S$ is finite.

\emph{Formula for $P$.} Divide \eqref{a:eq:telescope} by $r^m$ and let $m\to\infty$. Since
$u^*_{m+1}\to0$, we get $1/P=r\kappa'/\cos\alpha$.

\emph{Formula for $A$.} Use \eqref{a:eq:telescope} and $D_{m+1}=T'(u^*_m)D_m$ with
$\sec(\Delta^*_m-\alpha)=1/\cos(2\alpha-u^*_m)$. Every summand of \eqref{a:eq:Aalpha} then becomes
\[
  \Big(\prod_{j\le m}\rho_j\Big)\sec(\Delta^*_m-\alpha)=\frac{D_m}{\cos\alpha},
\]
so the bracket in \eqref{a:eq:Aalpha} equals $(1+2\sum_{m\ge1}D_m)/\cos\alpha=(2S-1)/\cos\alpha$.
Inserting $P$ gives $A$, and $C=4/A$.

(b) Write $u=\alpha y$ with $y\in[0,1]$, and $F(y):=T(\alpha y)/\alpha$. Let $K$ denote constants
depending only on $\alpha_0$. All functions involved are analytic, so Taylor expansion gives three
local estimates, uniformly on $[0,1]$.
\begin{enumerate}[label=(E\arabic*),nosep]
\item $F(y)=y-hy(2-y)+R_1(y)$ with $\abs{R_1}\le Kh^2y$.
\item $\ln\bigl(T'(\alpha y)/r\bigr)=2hy+R_2(y)$ with $\abs{R_2}\le Kh^2y$.
\item $\ln\bigl(F(y)/(ry)\bigr)=hy+R_3(y)$ with $\abs{R_3}\le Kh^2y$, and $\ln r=-2h+O(h^2)$.
\end{enumerate}
Each remainder vanishes at $y=0$ and has $O(h^2)$ $y$-derivative, which gives the factor $y$. For
example, $\frac{d}{dy}\ln T'(\alpha y)=\alpha\bigl[\tan(2\alpha-u)+\tan T(u)\,T'(u)\bigr]=2h+O(h^2)$,
because $\abs{T(u)-u}\le(1-r)u=O(\alpha^3)$ by \eqref{a:eq:Tbounds}.

\emph{The orbit.} Let $y_m:=u^*_m/\alpha$, $\tau_m:=hm$, and $Y(\tau):=1-\tanh\tau$. This is the
solution of the logistic equation $\dot Y=-Y(2-Y)$ with $Y(0)=1$. Taylor's theorem gives
$Y(\tau+h)=Y-hY(2-Y)+Q$ with $\abs Q\le\frac{h^2}2\sup\abs{\ddot Y}\le2h^2Y(\tau)$, because
$\ddot Y=2(1-Y)Y(2-Y)$ and $Y$ is decreasing. Subtracting from (E1), the error
$e_m:=y_m-Y(\tau_m)$ satisfies
\begin{gather*}
  e_{m+1}=\bigl(1-2h+h(y_m+Y(\tau_m))\bigr)e_m+R_1(y_m)-Q_m,\\\abs{R_1(y_m)-Q_m}\le Kh^2\bigl(Y(\tau_m)+\abs{e_m}\bigr).
\end{gather*}
The multiplier lies in $[0,1]$. So for $\tau_m\le\tau_1:=\operatorname{artanh}\frac34$ we get
$\abs{e_m}\le Kh$.

For $\tau_m\ge\tau_1$ we have $Y\le\frac14$. Proceeding by induction on $m\ge m_1$ with the
inductive hypothesis $\abs{e_m}\le K_2he^{-\tau_m}\le\frac18$ (valid at $m_1$ for $K_2$ large
enough and $h$ small enough), we have $y_m+Y(\tau_m)\le\frac58$, so the multiplier is at most
$1-\frac{11}8h+Kh^2\le 1-5h/4$ for small $h$. Since $Y(\tau)\le2e^{-2\tau}$,
$\abs{e_{m+1}}\le (1-5h/4)\abs{e_m}+2Kh^2e^{-2\tau_m}$. Taking $K_2\ge8K$
bounds the right-hand side by $K_2he^{-\tau_m}(1-h)\le K_2he^{-\tau_{m+1}}$. This gives
$\abs{e_m}\le K_2he^{-\tau_m}$ for all $m$. As $\Delta^*_m=\alpha(y_m-1)$ and $Y-1=-\tanh$,
this is \eqref{a:eq:kink}.

\emph{Sums along the orbit.} $Y$ is decreasing with $\int_0^\infty Y=\ln2$. So
$h\sum_{j\ge0}Y(\tau_j)=\ln2+O(h)$, and by \eqref{a:eq:kink}, $h\sum_{j\ge0}y_j=\ln2+O(h)$. Likewise,
uniformly in $m$,
\[
2h\sum_{j<m}y_j=2\int_0^{\tau_m}Y+O(h)=2\tau_m-2\ln\cosh\tau_m+O(h).
\]
By (E2), $\ln\kappa'=\sum_{j\ge0}\bigl(2hy_j+O(h^2y_j)\bigr)=2\ln2+O(h)$.
By (E3), $\ln(c_*/\alpha)=\sum_{j\ge0}\ln\bigl(y_{j+1}/(ry_j)\bigr)=\ln2+O(h)$.

By (E2), (E3) and the last display,
\begin{align*}
  D_m&=r^m\exp\Big(\sum_{j<m}\bigl(2hy_j+O(h^2y_j)\bigr)\Big)\\&=\operatorname{sech}^2(\tau_m)\,e^{O(h(1+\tau_m))}.
\end{align*}
For small $h$, $\operatorname{sech}^2(\tau)\,(1+\tau)\,e^{Kh\tau}\le4(1+\tau)e^{-\tau}$, which is
integrable. Since $\operatorname{sech}^2$ is decreasing with integral $1$, it follows that
$hS=h\sum_m\operatorname{sech}^2(\tau_m)+O(h)=1+O(h)$.

\emph{Conclusion.} Insert these into (a), using $1-r=2\sin^2\alpha=2h(1+O(h))$:
$P=\cos\alpha/(r\kappa')=\frac14(1+O(h))$ and
$C=2(1-r)\kappa'/(2S-1)=8h^2(1+O(h))$.
\end{proof}
\noindent\numcheck{Computed with 40 digits for $\alpha=0.00625$--$0.7$. The exact
representation (a) agrees with the direct formulas \eqref{a:eq:Aalpha} to $\le10^{-36}$. The error
in \eqref{a:eq:kink}, scaled as $e^{\alpha^2m}/\alpha^3$, has maximum over $m$ tending to $0.51$ as
$\alpha\to0$, and at most $0.553$ up to $\alpha=0.7$. Each quantity minus its limit, divided by $\alpha^2$,
converges as $\alpha\to0$.}

\begin{apprem}
Richardson extrapolation over $\alpha=0.00625$--$0.05$ gives
\begin{gather*}
  P=\tfrac14+\tfrac1{16}\alpha^2,\qquad \kappa'=4+5\alpha^2,\qquad \alpha^2S=1+\tfrac14\alpha^2,\\
  c_*=2\alpha\bigl(1+\tfrac56\alpha^2\bigr),\qquad C=8\alpha^4\bigl(1+\tfrac76\alpha^2\bigr),
\end{gather*}
each up to $O(\alpha^4)$; the fitted coefficients agree with these fractions to at least four digits. The five
expansions are mutually consistent through (a), but they are not proved here; a proof would need
the second-order modified equation of $T$.

The lag profile is a $\tanh$ kink of width $1/\alpha^2\approx2\xi$. At finite $N=2n$,
\eqref{a:eq:treeratio} uses the lags of the finite-chain shooting orbit. For fixed $m$,
that orbit converges to the center orbit as $n\to\infty$, so \eqref{a:eq:telescope} gives
\[
  \lim_{n\to\infty}\frac{w_{n-m}}{w_n}
  =\frac{\cos\alpha}{D_{m+1}\cos u^*_{m+1}} .
\]
The weak-skin profile of this ratio is $\cosh^2(\alpha^2m)$.
Also $\mu-\mu_0\sim(1-r)c_*r^{N/2}\approx4\alpha^3r^{N/2}$.
\end{apprem}

\emph{General coupling and odd $N$.}
The finite-chain existence and spectral arguments in \cref{a:lem:twosource,a:thm:sink}
extend beyond the Sakaguchi family under the following assumptions.

\begin{appassum}\label[appassum]{a:ass:G}
$\Gamma$ has a zero $x_0>0$ such that:
\begin{enumerate}[label=(G\arabic*),nosep]
\item $\Gamma'>0$ on $[-x_0,x_0]$;
\item $\Gamma_e(x)<\Gamma_e(x_0)$ for $0\le x<x_0$, where $\Gamma_e(x):=\frac12(\Gamma(x)+\Gamma(-x))$ is the
even part;
\item $r:=\Gamma'(-x_0)/\Gamma'(x_0)<1$.
\end{enumerate}
Let $(-W,W)$, with $W>x_0$, be the largest symmetric interval on which $\Gamma'>0$.
\end{appassum}

The sign $x_0>0$ is no restriction since the reflection of \cref{a:prop:reflect} maps $x_0$ to $-x_0$ and
leaves $r$ unchanged. For Sakaguchi coupling, $x_0=\alpha$ and $\Gamma_e(x)=-\sin\alpha\cos x$, so (G2)
holds. (G1) is equivalent to $\alpha<\pi/4$, and $W=\frac\pi2-\alpha$.

(G2) is what makes the lag recursion approach $-x_0$ monotonically: by step 1 of the proof below,
it is equivalent to $T(u)<u$ on $(0,x_0]$. If it fails, $T(u)\ge u$ somewhere in $(0,x_0]$, and the
orbit from the center need not reach $-x_0$.

\emph{Notation.} Let $\mu_0:=\Gamma(-x_0)=2\Gamma_e(x_0)$, let $\Gamma_+^{-1}$ be the inverse of $\Gamma$ on
$(-W,W)$, and set
\begin{gather*}
  T_\delta(u):=x_0-\Gamma_+^{-1}\bigl(\mu_0+\delta-\Gamma(u-x_0)\bigr),\\T:=T_0,\qquad
  \rho(\Delta):=\Gamma'(\Delta)/\Gamma'(-\Delta).
\end{gather*}
For Sakaguchi coupling this is the map of \cref{a:lem:twosource}. The center value of the orbit
depends on the parity of $N$:
\begin{gather*}
  u_0(\delta):=x_0\quad(N=2n),\\u_0(\delta):=x_0-\Gamma_+^{-1}\bigl(\tfrac12(\mu_0+\delta)\bigr)\quad(N=2n+1).
\end{gather*}
The \emph{center orbit} is $u^*_m:=T^m(u_0(0))$, with $\Delta^*_m:=u^*_m-x_0$. Write $D_m:=(T^m)'(u^*_0)$,
$\kappa':=\lim_mr^{-m}D_m$, $S:=\sum_{m\ge0}D_m$ and $P:=\prod_{j\ge1}r/\rho(\Delta^*_j)$, and
$\lfloor N/2\rfloor=n$ in both cases.

\begin{appthm}[General coupling, both parities]\label{a:thm:general}
Let the chain have uniform weight $a$, identical oscillators, coupling satisfying \cref{a:ass:G}, and
$N\ge2$.
\begin{enumerate}[label=(\alph*)]
\item \emph{(\Cref{a:lem:twosource}.)} There is exactly one symmetric locked state with all lags in
$(-W,W)$. It is positive and stable.

Its left-half lags are $\Delta_{n-m}=u_m-x_0$, where $u_0=u_0(\delta_N)$ and
$u_{m+1}=T_{\delta_N}(u_m)$, with $0<u_m\le u^*_m$ for $m<n$. For odd $N$ the center node is $n+1$, and its two edges carry
$\pm\Delta_n$.

The center orbit decreases to $0$ with $u^*_m\asymp r^m$, and the frequency shift satisfies
\[
  K^{-1}u^*_n\le\delta_N\le\max_{[0,x_0]}\Gamma'\cdot u^*_n,
\]
with $K$ depending only on $\Gamma$. Moreover $w_{k+1}/w_k=\rho(\Delta_k)=r\,(1+O(r^{n-k}))$.
\item \emph{(\Cref{a:thm:sink}(i), (ii).)} The slow mode is antisymmetric, and for every $N$
\begin{gather*}
  \frac4{R_N}\le\abs{\lambda_2}\le\frac4{R_NV_N},\qquad 1-V_N=O(Nr^{n}),\\\abs{\lambda_3}\ge\frac cN,
\end{gather*}
where the last bound is for $N\ge4$, with $c>0$ depending only on $a$ and $\Gamma$.
\item \emph{(\Cref{a:thm:sink}(iii), \cref{a:prop:continuum}(a).)} $a\,r^nR_N\to A$, where, with the center
orbit of the matching parity,
\begin{alignat*}{3}
  A&=\frac{2(2S-1)}{(1-r)\,\kappa'\,\Gamma'(x_0)},&\qquad P&=\frac{\Gamma'(0)}{r\kappa'\Gamma'(x_0)}&\qquad&(N\text{ even}),\\
  A&=\frac{4S}{(1-r)\,\kappa'\,\Gamma'(x_0)},&\qquad P&=\frac{\Gamma'(\Delta^*_0)}{r\kappa'\Gamma'(x_0)}&\qquad&(N\text{ odd}).
\end{alignat*}
Hence $\lambda_2=-(4a/A)\,r^{n}\,(1+o(1))$. Equivalently, in the form \eqref{a:eq:Aalpha},
$A=\frac{2rP}{1-r}\bigl[\frac1{\Gamma'(0)}+2\sum_{m\ge1}\frac{\prod_{j\le m}\rho_j}{\Gamma'(\Delta^*_m)}\bigr]$
for even $N$, and $A=\frac{4rP}{1-r}\sum_{m\ge0}\frac{\prod_{1\le j\le m}\rho_j}{\Gamma'(\Delta^*_m)}$ for odd $N$.
\end{enumerate}
\end{appthm}

\begin{proof}
We list each place where the Sakaguchi proofs used the explicit form, and what replaces it.

\emph{1. The map $T$ (replaces \cref{a:lem:twosource}(a)).} Let $u\in[0,x_0]$ and $\Delta=u-x_0$. By (G1),
$\Gamma$ is increasing on $[-x_0,x_0]$. By (G2) at $x=0$, $\Gamma(0)\le\mu_0/2$. Hence
$\mu_0-\Gamma(\Delta)\in[\mu_0-\Gamma(0),0]\subset[\Gamma(0),\Gamma(x_0)]$, so $T$ maps $[0,x_0]$ into itself.
$T$ is smooth and increasing, with $T(0)=0$ and
\[
  T'(u)=\Gamma'(\Delta)/\Gamma'(x_0-T(u)),\qquad T'(0)=r .
\]
Since $\Gamma$ is increasing, $T(u)<u$ is equivalent to $\Gamma(-\Delta')>\Gamma(-\Delta)$ for the image
lag $\Delta'$, that is, to $\Gamma_e(\Delta)<\Gamma_e(x_0)$, which is (G2). So the center orbit decreases to
the unique fixed point $0$.

For odd $N$, $u_0(0)\in(0,x_0)$, because $\Gamma(0)<\mu_0/2<0=\Gamma(x_0)$.

Near $0$, $T(u)=ru+O(u^2)$. Hence $\ln(u^*_{m+1}/(ru^*_m))=O(u^*_m)$ is summable, and
$u^*_m\asymp r^m$. Likewise $\abs{\ln(\rho(\Delta^*_m)/r)}\le Lu^*_m$, so $P$, $\kappa'$ and $S$ are finite
and positive. The explicit bounds \eqref{a:eq:Tbounds} are replaced by constants depending on $\Gamma$.

\emph{2. Shooting (replaces \cref{a:lem:twosource}(b), (c)).}
Write $H_\delta(u):=\mu_0+\delta-\Gamma(u-x_0)$ for the inverse argument.
In a solution with lags in $(-W,W)$, if $\delta\le0$ then the center input is positive
and $H_\delta(u)<0$ at every positive input, so the orbit cannot end at zero.
If $\delta>0$ and an intermediate input is nonpositive, then $H_\delta(u)>0$,
and the next iterate is negative whenever defined. Thus in a solution $\delta>0$
and all intermediate inputs are positive. Since $T_\delta<T$ wherever both are defined
on $[0,x_0]$, and the odd center input decreases with $\delta$, all intermediate inputs
lie below the center orbit. The terminal condition $H_\delta(u_{n-1})=0$ then implies
$\delta<-\mu_0$, since $\Gamma(u_{n-1}-x_0)\le\Gamma(0)<0$.

For existence on $[0,-\mu_0]$, iterate
$x_0-\Gamma_+^{-1}(\min\{0,H_\delta(u)\})$ up to step $n-1$. On $[0,x_0]$ the
clipped inverse argument belongs to $[\Gamma(0),0]$ by (G2).
The iterates are continuous and nonincreasing in $\delta$. This holds also for
the odd center input, which decreases to zero. Hence $H_\delta(u_{n-1})$ increases
strictly from a negative value at $\delta=0$ to a positive one at $-\mu_0$.
At its unique zero, no intermediate input was clipped to zero. Such an input would
force all following inputs to vanish and give $H_\delta(0)=\delta>0$.
The root therefore solves the original equations.

For $\delta\in[0,\delta_N]$ the actual orbit is in $[0,x_0]$, so its terminal
derivative satisfies $\Phi_n'\le-1/\max_{[0,x_0]}\Gamma'$. This gives the upper bound
on $\delta_N$. For the lower bound, choose $u_q>0$ so that
$\Gamma'(u-x_0)/\Gamma'(x_0-v)\le q:=(1+r)/2$ for all $u,v\in[0,u_q]$.
The finite orbit is bounded above by the center orbit, so both inputs of this ratio
are in $[0,u_q]$ after a bounded number of steps. The derivative sum of
\eqref{a:eq:dPhi}, with $\Gamma'$ replacing cosines, is then bounded uniformly in $n$.
The derivative of the odd center input adds a bounded term multiplied by the same
products. Integrating the resulting $|\Phi_n'|\le K$ gives the lower bound.

\emph{3. Spectral bounds (replaces \cref{a:thm:sink}(i), (ii)).} Symmetry gives
$W_{N+1-k,N-k}=a\Gamma'(-\Delta_{N-k})=a\Gamma'(\Delta_k)=W_{k,k+1}$.

For even $N$ the parity argument and the two-sided $\abs{\lambda_2}$ bound are unchanged.
For odd $N$, the reflection fixes the center node and odd
functions vanish there. The odd sector is therefore the even-sector operator with the center node
deleted, so Cauchy interlacing (strict, for an unreduced tridiagonal) again gives
$\abs{\lambda_2}=\lambda_{\rm odd}<\lambda_{\rm even}=\abs{\lambda_3}$.

In the lower bound, the half chain is now killed through $c_n$ rather than $2c_n$. Then
\[
  t(1)\le\tfrac12\sum_{k\le n}c_k^{-1}=\tfrac14R_N,
\]
because $R_N=2\sum_{k\le n}c_k^{-1}$. The test function $1-2h_N$ is still odd, since
$h_N(n+1)=\frac12$.

For $V_N$ and $\lambda_3$, the Sakaguchi facts ``$w$ is nonincreasing'' and
``$c_j\ge a\cos\alpha\,w_{j+1}$'' are replaced as follows:
\begin{itemize}[nosep]
\item $w_i\le Mr^{i-j}w_j$ for $j\le i\le n$ and also for $i=n+1$ when $N$ is odd,
with $M:=e^{L\sum_mu^*_m}$;
\item $c_j=aw_{j+1}\Gamma'(-\Delta_j)\ge am_\Gamma w_{j+1}$, with $m_\Gamma:=\min_{[0,x_0]}\Gamma'$.
\end{itemize}
Then $w(\{j+1,\dots,n\})/c_j\le M/((1-r)am_\Gamma)$. This gives the even-case bounds
$1-V_N\le8nM/((1-r)am_\Gamma R_N)$ and $\sum_{j<n}w(\{j+1,\dots,n\})/c_j\le nM/((1-r)am_\Gamma)$.
For odd $N$, the folded even sector has center mass $w_{n+1}/2$ when the Dirichlet form on
the half chain is normalized as above. The hitting-time sum acquires at most
$\sum_{j\le n}w_{n+1}/(2c_j)\le nM/(2am_\Gamma)$, using detailed balance and the bound on
the influence ratios. The center node also contributes $w_{n+1}=O(r^n)$ to $1-V_N$.
To convert the hitting-time estimate to a spectral bound, let $\mu_i$ be the folded even-sector
masses and use $v_1$ as a candidate center of the weighted variance. Cauchy--Schwarz along the
path from $1$ to $i$ gives
\begin{align*}
\operatorname{Var}_\mu(v)
 &\le \sum_i\mu_i(v_i-v_1)^2\\
 &\le \mathcal E_e(v)\sum_j\frac{\sum_{i>j}\mu_i}{c_j}.
\end{align*}
Here $\mathcal E_e$ sums the folded edge energies. The sum is the same hitting-time
sum bounded above by a constant times $N$. The variational principle therefore gives
$\abs{\lambda_3}\ge c/N$. Finally, $R_N^{-1}\le c_n\le a\max\Gamma'\,w_n=O(r^n)$.

\emph{4. Asymptotics and representation (replaces \cref{a:thm:sink}(iii), \cref{a:prop:continuum}(a)).}
Repeat step (iii) of \cref{a:thm:sink} with $c_{n-m}=aw_{n-m}\Gamma'(\Delta_{n-m})$. For even $N$ the
center edge has $\Delta=0$ and is counted once. For odd $N$ every edge is counted twice.

The telescoping identity becomes
\begin{gather*}
  \rho(\Delta^*_j)=T'(u^*_j)\,\frac{\Gamma'(-\Delta^*_{j+1})}{\Gamma'(-\Delta^*_j)},\\
  \prod_{j=1}^m\rho(\Delta^*_j)=\frac{D_{m+1}\,\Gamma'(-\Delta^*_{m+1})}{\Gamma'(\Delta^*_0)},\\
  \frac{\prod_{1\le j\le m}\rho_j}{\Gamma'(\Delta^*_m)}=\frac{D_m}{\Gamma'(\Delta^*_0)},
\end{gather*}
using $T'(u^*_0)=\Gamma'(\Delta^*_0)/\Gamma'(-\Delta^*_1)$. Since $\Gamma'(-\Delta^*_{m+1})\to\Gamma'(x_0)$, the
formulas for $P$ and $A$ follow as in \cref{a:prop:continuum}(a).
\end{proof}
\noindent\numcheck{Computed with 50 digits. The test couplings are Sakaguchi
$\alpha=0.4$ ($r=0.70$), $\sin(x-0.4)+0.2\sin2x$ ($r=0.83$) and
$\sin(x-0.35)+0.25(1-\cos2x)$ ($r=0.39$), with $N=4$--$61$ of both parities, 42 states in all.}
\begin{itemize}[nosep]
\item All three couplings satisfy \cref{a:ass:G}, and both center orbits are monotone.
\item Every node equation holds to $3\times10^{-45}$. The shooting root lies in the bracket of (a)
  in every case, and $\delta_N/u^*_n\to(1-r)\Gamma'(x_0)$.
\item The bounds of (b) hold in all 42 cases, with an odd slow mode and $\abs{\lambda_3}=\lambda_{\rm even}$.
\item The representation in (c) agrees with the direct sums to $2\times10^{-45}$. At $N=60,61$ the
  measured constant $\abs{\lambda_2}r^{-n}$ matches $4/A$ with relative deviations
  $10^{-4}$, $8\times10^{-3}$ and $10^{-11}$ for the three couplings.
\item An independent double-precision Newton continuation agrees to $\le6\times10^{-10}$ in the lags
  (24 runs).
\end{itemize}

\begin{apprem}\label[apprem]{a:rem:weakskin}
Consider families with $x_0\to0$ and $\ell:=\Gamma''(0)x_0/\Gamma'(0)\to0^+$.
The following formal calculation assumes, uniformly for $0<x\le x_0$, that
\[
  \frac{\Gamma_e'(x)}{\Gamma''(0)x}=1+o(1),\qquad
  \frac{\Gamma_o'(x)}{\Gamma'(0)}=1+o(1),
\]
where $\Gamma_o$ is the odd part. These conditions give $1-r\sim2\ell$ at $x=x_0$.
Integrating the even-part slope from $\abs{\Delta}$ to $x_0$ gives the leading drift
$2(\Gamma_e(x_0)-\Gamma_e(\Delta))\simeq\Gamma''(0)(x_0^2-\Delta^2)$ in the lag recursion.
The resulting continuum kink is $\Delta^*_m\approx-x_0\tanh(\ell m)$.
Assumptions (G1)--(G3) alone also allow a quartic leading even part, for which this reduction fails.
Repeating the argument of \cref{a:prop:continuum}(b) formally gives $\kappa'\to4$, $\ell S\to1$ and
$P\to\frac14$. For both parities this yields
\[
  \abs{\lambda_2}\approx2a\,\Gamma'(0)\,(1-r)^2\,r^{N/2},
\]
for long chains at fixed coupling followed by the weak-skin limit. For odd $N$ the exact exponent is
$\lfloor N/2\rfloor$; the extra factor $\sqrt r$ tends to $1$ in this limit. For Sakaguchi coupling it
reduces to $C\simeq8\sin^4\alpha$.

\numcheck{On three one-parameter families with $x_0\to0$, the ratio of
$C$ to $2\Gamma'(0)(1-r)^2$ tends to $1$ for both parities: $1.0009$ (even) and $1.0006$ (odd) for
$\sin(x-\varepsilon)+0.2\sin2x$ at $\varepsilon=0.025$; $1.05$ and $1.03$ for
$\sin(x-\varepsilon)+0.25(1-\cos2x)$, which converges more slowly because there
$\ell\propto x_0$; and $1.0015$ for Sakaguchi. In all three, $P\to\frac14$.}
For even Sakaguchi chains, \cref{a:prop:continuum} proves the uniform estimates needed for this prefactor.
\end{apprem}

\section{Product graphs and corner localization}\label{a:sec:product}

\begin{appprop}\label[appprop]{a:thm:product}
Let $G$ and $H$ be finite connected graphs with symmetric weights, and take identical oscillators.
Let $(u,\Omega_G)$ and $(v,\Omega_H)$ be positive locked states on $G$ and $H$, with Jacobians
$J_G$, $J_H$ and influences $w_G$, $w_H$. Then on the Cartesian product $G\square H$:
\begin{enumerate}[label=\textup{(\alph*)}]
\item $\theta_{(x,y)}=u_x+v_y$ is a positive locked state, with $\Omega-\omega=(\Omega_G-\omega)+(\Omega_H-\omega)$.
\item Its Jacobian is the Kronecker sum $J=J_G\otimes I+I\otimes J_H$, so
$\operatorname{spec}J=\{\lambda+\mu\}$.
\item $w=w_G\otimes w_H$.
\item If $u$ and $v$ are the in-phase-continued branches along a common parameter path, then so is
$\theta$.
\end{enumerate}
\end{appprop}

\noindent Parts (a) and (d) extend the reduction of rectangular and cubic arrays to independent
chains by Ren and Ermentrout~\cite[Lemma~3.3, Thm.~3.4 and Remark~4]{RenErmentrout1998} to general
Cartesian products; parts (b) and (c), the Kronecker-sum Jacobian and the factorized influence, are not there.

\begin{proof}
The node $(x,y)$ has horizontal neighbors $(x',y)$, weight $a^G_{xx'}$, and vertical neighbors
$(x,y')$, weight $a^H_{yy'}$.

\smallskip
(a) It satisfies
$\omega+\sum_{x'}a^G_{xx'}\Gamma(u_{x'}-u_x)+\sum_{y'}a^H_{yy'}\Gamma(v_{y'}-v_y)$, which equals
$\omega+(\Omega_G-\omega)+(\Omega_H-\omega)$.

(b) The horizontal weight $a^G_{xx'}\Gamma'(u_{x'}-u_x)=(W_G)_{xx'}$ does not depend on $y$, and
likewise for vertical edges.

(c) $(w_G\otimes w_H)^\trans J=(w_G^\trans J_G)\otimes w_H^\trans+w_G^\trans\otimes(w_H^\trans J_H)=0$.
The vector is positive and normalized, hence it is $w$; refer to \cref{a:prop:generator}.

(d) Along the path the product solves the product equations. Its bordered Jacobian is nonsingular
whenever $0$ is simple for $J$. By (b), $\lambda+\mu=0$ requires $\lambda=\mu=0$, since all other
eigenvalues have negative real part. So the continued branch is unique and equals the product.
For a Sakaguchi box, each chain has the positive symmetric branch of \cref{a:thm:general} at
$0<\alpha<\pi/4$. This family is the in-phase-continued branch. The symmetric positive state is unique (\cref{a:thm:general}(a)) and
has a simple zero eigenvalue (\cref{a:prop:generator}), so the implicit function theorem continues it locally in
$\alpha$. The continued state is again symmetric, because its mirror image is a continuation of the same state
and local continuation is unique, and it is positive, because the rates are continuous. By uniqueness it is the
symmetric positive state at the nearby $\alpha$, so the family is one smooth branch on $(0,\pi/4)$. Its lags lie
in $[-\alpha,\alpha]$, so it converges to the in-phase state as $\alpha\to0$, where continuation from $\alpha=0$
is unique (\cref{a:thm:degree}).
The reflection symmetry covers negative $\alpha$.
\end{proof}
\noindent\numcheck[K2]{The errors in $\theta$ and $w$ are $e_\theta\le1.1\times10^{-12}$ and $e_w\le3.9\times10^{-12}$ on $16^2$--$64^2$, $8^3$ and $12^3$ lattices.}

\begin{appcor}\label[appcor]{a:cor:corner}
Take a box $P_{L_1}\square\cdots\square P_{L_D}$ of chains of any lengths $L_d\ge2$, with a coupling
satisfying \cref{a:ass:G} (e.g., Sakaguchi with $\alpha \in (0,\pi/4)$), in the product of the
symmetric two-source states of \cref{a:thm:general}. Then
$w(x)=\prod_dw^{(d)}(x_d)$. Let $\operatorname{dist}(x)$ be the lattice distance to the boundary, and
$S_d(m):=w^{(d)}(\{x_d:\min\{x_d-1,L_d-x_d\}\ge m\})$. Then
\[
\mathbb E_w[\operatorname{dist}]=\sum_{m\ge1}\prod_{d=1}^DS_d(m)\ \xrightarrow[L_d\to\infty]{}\ \frac{r^D}{1-r^D},
\]
and the share of $w$ in the $2^D$ corner boxes of side $\ell$ tends to $(1-r^\ell)^D$.
\end{appcor}

\begin{proof}
The identity $\mathbb E[X]=\sum_{m\ge1}\mathbb P(X\ge m)$ holds for $X\in\{0,1,\dots\}$. Under the
product measure, $\mathbb P(\min_dX_d\ge m)=\prod_d\mathbb P(X_d\ge m)$.

Each factor is a two-sided skin. By symmetry, $S_d(m)=2\sum_{m<k\le n}w^{(d)}_k$ with
$n=\lfloor L_d/2\rfloor$, plus the center mass $w^{(d)}_{n+1}=O(r^n)$ when $L_d$ is odd and $m\le n$.
As in step (3)
of the proof of \cref{a:thm:sink}(iii), $w^{(d)}_k\to\frac{1-r}2r^{k-1}$ for each $k$, with domination
by $Mr^{k-1}$ (\cref{a:thm:general}, step 3). Hence $S_d(m)\to r^m$, and $S_d(m)\le Kr^m$ uniformly
in $L_d$.
Dominated convergence in $m$ gives both limits.
\end{proof}

\begin{apprem}
There is a single zero mode, and its profile is separable. This is \emph{not} an intrinsic
second-order skin effect in the sense of \cite{Kawabata2020}, where corner modes are
$O(L)$ in number. We refer to this as ``corner localization from product structure''.

On the infinite quadrant the corner wave has frequency shift $\Omega-\omega=2a\Gamma(-x_0)$.
On the half-plane the straight-edge wave has shift $a\Gamma(-x_0)+2a\Gamma(0)$.
The corner is faster when $\Gamma(-x_0)>2\Gamma(0)$, which for Sakaguchi coupling at
$0<\alpha<\pi/4$ follows from $\sin2\alpha<2\sin\alpha$.
\end{apprem}

\section{Twisted states and the reciprocity defect}\label{a:sec:twist}

\begin{appprop}\label[appprop]{a:thm:recip}
Under \cref{a:ass:P}, let $P(t)=e^{Jt}$. Here $P_{ij}(t)$ is the response of node $i$ at time $t$ to a
unit kick at node $j$, and $P_{ij}(t)>0$ for $t>0$. For $t>0$ define the \emph{reciprocity defect}
$D_{ij}(t):=\ln\frac{w_iP_{ij}(t)}{w_jP_{ji}(t)}$. The following are equivalent:
\begin{enumerate}[label=\textup{(\roman*)}]
\item $f$ is exact (zero circulation on every cycle);
\item detailed balance holds;
\item $w_iP_{ij}(t)=w_jP_{ji}(t)$ for all $i,j$ and all $t\ge0$;
\item $D_{ij}(t)=0$ for all $i,j$ and $t>0$.
\end{enumerate}
If these hold, the spectrum of $J$ is real: $\operatorname{spec}J\subset\R$.
\end{appprop}

\begin{proof}
(i)$\Leftrightarrow$(ii) is \cref{a:cor:kolmogorov}.

(ii)$\Rightarrow$(iii): with $\Pi=\diag(w)$, (ii) says $\Pi J=J^\trans\Pi$. Then
$\Pi J^n=(J^\trans)^n\Pi$ for every $n$, so $\Pi e^{Jt}=(e^{Jt})^\trans\Pi$.

(iii)$\Rightarrow$(ii): differentiate at $t=0$.

(iii)$\Leftrightarrow$(iv) by definition for $t>0$ and by continuity for the equality in (iii) at $t=0$.

Finally, $\Pi^{1/2}J\Pi^{-1/2}$ is symmetric under (ii), so the spectrum is real.
\end{proof}

\begin{appprop}\label[appprop]{a:prop:ring}
On the ring $C_N$ with $N\ge3$, identical oscillators and uniform $a$, let $\Delta=2\pi q/N$ with
$|\Delta|<\pi$, and assume $\Gamma'(\pm\Delta)>0$. For Sakaguchi coupling this requires
both $\cos(\Delta-\alpha)>0$ and $\cos(\Delta+\alpha)>0$. On the principal attractive branch
$|\alpha|<\pi/2$, these conditions reduce to $|\Delta|+|\alpha|<\pi/2$. The $q$-twisted state
$\theta_i=i\Delta$ has the following properties.
\begin{enumerate}[label=\textup{(\alph*)}]
\item The rates are $A_\pm=a\Gamma'(\pm\Delta)$ for the jumps $i\to i\pm1$, and $w$ is uniform.
\item The circulation is $NF(\Delta)$ and the current on every oriented edge is
$\jmath=(A_+-A_-)/N$. Both are nonzero precisely when $A_+\ne A_-$. If the Taylor
interval in \cref{a:prop:circ} applies, $NF(\Delta)=4\pi\gamma q+O(|q|^3/N^2)$.
\item The eigenvalues are 
\[
\lambda_k=A_+(e^{2\pi ik/N}-1)+A_-(e^{-2\pi ik/N}-1),
\]
with $\operatorname{Im}\lambda_k=(A_+-A_-)\sin(2\pi k/N)$.
\item The propagator is
\begin{align*}
P_{0m}(t)={}&e^{-(A_++A_-)t}\sum_{n\in\Z}(A_+/A_-)^{(m+nN)/2}\\
&\quad \times I_{m+nN}\big(2\sqrt{A_+A_-}\,t\big).
\end{align*}
\item For $0<m<N/2$, $D_{0m}(t)\to mF(\Delta)$ as $t\downarrow0$, and $D_{0m}(t)\to0$ as $t\to\infty$.
\end{enumerate}
\end{appprop}

\begin{proof}
(a)--(c): $J$ is circulant, so its column sums also vanish and $\one$ is a left null vector.
Fourier modes diagonalize the two cyclic shift matrices. The current and cycle affinity
follow from the constant rates. Apply \cref{a:prop:circ} under its Taylor hypothesis.

(d): the lift to $\Z$ is a walk with $\pm1$ jumps at rates $A_\pm$. Its displacement has the Skellam
law $\mathbb P(X_t=k)=e^{-(A_++A_-)t}(A_+/A_-)^{k/2}I_k(2\sqrt{A_+A_-}t)$; see \cite{Skellam1946}.
Projecting modulo $N$ gives the image sum.

(e): $I_{-m}=I_m$, so the ratio of the two unwrapped kernels at displacements $m$ and $-m$ is
$(A_+/A_-)^m$ for every $t>0$. On the ring, $I_k(x)\sim(x/2)^{\abs k}/\abs k!$ as $x\downarrow0$.
For $0<m<N/2$ the unwrapped terms dominate both $P_{0m}$ and $P_{m0}=P_{0,-m}$ as $t\downarrow0$,
and $w$ is uniform. As $t\to\infty$, $P\to\one w^\trans$.
\end{proof}
\noindent\numcheck[C7, K3]{The image sum matches \texttt{expm} to $2\times10^{-15}$;
$D_{0,15}(0.01)\approx15F(2\pi/60)\approx0.9757204$. Direct nonlinear kicks of size $10^{-3}$
agree to $2\times10^{-4}$; the reversible chain and near-reversible heterogeneous ring give
defects of magnitude at most $1.4\times10^{-4}$ at this kick size.}

\begin{apprem}A comment about prior art and measurement:\leavevmode
\begin{itemize}
\item The complex twisted-state spectrum (c) is known \cite{TonjesBlasius2009a,Mihara2018,KimKori2025}.
The paired-kick measurements use the stationary detailed-balance ratio to probe the circulation.
\item $D$ needs two kick experiments, $j\to i$ and $i\to j$. The weights $w_i,w_j$ come from the
same experiments, as the permanent shifts.
\item It separates \emph{gauge} non-reciprocity, for which $D\equiv0$ even when $P_{ij}/P_{ji}$ is far
from $1$ (skin chains), from \emph{flux} non-reciprocity.
\end{itemize}
\end{apprem}

\section{Existence of locked states: a cut criterion}\label{a:sec:cut}

\begin{appprop}[Cut identity]\label[appprop]{a:prop:cut}
Take identical oscillators with $\Gamma(x)=\sin(x-\alpha)$ and $a_{ij}=KA_{ij}$, $A$ unweighted. Let
$k_i$ be the degrees, $\bar k=2\abs E/N$, and for $S\subset V$ let $\abs{\partial S}$ be the number of
edges leaving $S$. Then every locked state satisfies
\begin{gather}
\begin{aligned}
&\textstyle\sum_{i\in S,\,j\notin S}A_{ij}\sin(\Delta_{ij}-\alpha)\\
&\quad=\abs S\,\frac{\Omega-\omega}K+2\sin\alpha\sum_{\{i,j\}\subset S}A_{ij}\cos\Delta_{ij},
\end{aligned}\label{a:eq:cutid}\\
\textstyle\frac{\Omega-\omega}K=-\frac{2\sin\alpha}N\sum_{\{i,j\}\in E}\cos\Delta_{ij}.\label{a:eq:Omega}
\end{gather}
If moreover the principal lags satisfy $\abs{\Delta_e}\le d\le\pi$ on all edges, then for every $S$,
\begin{equation}\label{a:eq:cutcrit}
\abs{\sin\alpha}\Big(\cos d\,\abs S\bar k-2e_{\mathrm{int}}(S)\Big)
\le\abs{\partial S},
\end{equation}
where $e_{\mathrm{int}}(S)$ is the number of edges within $S$.
\end{appprop}

\begin{proof}
Sum the locking equations over $S$, and use $\Gamma(x)+\Gamma(-x)=-2\sin\alpha\cos x$ on the internal
edges; this gives \eqref{a:eq:cutid}. Summing over $V$ gives \eqref{a:eq:Omega}.

By reflection (\cref{a:prop:reflect}) it suffices to consider $\sin\alpha\ge0$. The left side of
\eqref{a:eq:cutid} is at least $-\abs{\partial S}$. On the right side, use
$\sum_E\cos\Delta\ge\abs E\cos d$ in \eqref{a:eq:Omega}, $\sum_{\mathrm{int}}\cos\Delta\le e_{\mathrm{int}}(S)$,
which yields \eqref{a:eq:cutcrit}.
\end{proof}
\noindent\numcheck[C8]{Equation~\eqref{a:eq:cutid} holds to $7\times10^{-16}$ on 1015 sets of the ER(1000) state at $\alpha=0.09$; \eqref{a:eq:cutcrit} has no violations.}

\begin{appcor}[Leading cut estimate; heuristic]
Dropping the $O(d^2)$ term in \eqref{a:eq:cutcrit} gives the estimate
$|\sin\alpha|\lesssim1/(n(\bar k-2)+2)$ for a dangling tree of $n$ nodes attached by one bridge
when the denominator is positive, and $|\sin\alpha|\lesssim k_i/\bar k$ for a single node of
degree $k_i$.
\end{appcor}

\begin{apprem}\leavevmode
\begin{itemize}
\item The finite-lag condition \eqref{a:eq:cutcrit} is necessary. The leading estimate omits its
$d$-dependent term.
\item A numerical stop need not correspond to a fold. In particular, on free-boundary boxes
the product construction gives positive stable locks for $|\alpha|<\pi/4$.
\item This is the cut-set form of the classical single-node condition
\cite[Lemma~4.2 and Sec.~4.3]{DorflerBullo2014}, applied to the effective frequencies $-K k_i\sin\alpha$.
\end{itemize}
\end{apprem}

\bibliography{refs}

\end{document}